\documentclass[a4paper,11pt]{elsarticle}

\journal{European Journal of Operational Research}

\usepackage{hyperref}
\usepackage{amsthm}
\usepackage{amssymb}
\usepackage{amsmath}
\usepackage{adjustbox}
\usepackage{threeparttable}
\usepackage{multirow}
\usepackage{booktabs}
\usepackage{enumitem}
\setlist{noitemsep}
\usepackage{makecell}
\usepackage{rotating}
\usepackage{caption}
\usepackage{geometry}
\usepackage{xcolor}
\usepackage{soul}
\newtheorem{decisionrule}{Decision Rule}
\newtheorem{definition}{Definition}
\newtheorem{remark}{Remark}

\DeclareMathOperator*{\argmax}{arg\,max}

\renewcommand{\baselinestretch}{1.3}

\biboptions{authoryear}

\definecolor{Mycolor1}{HTML}{EA5A00}
\definecolor{Mycolor2}{HTML}{6E55DE}
\definecolor{Mycolor3}{HTML}{D22378}

\newcommand{\RefA}[2]{\hspace{0pt}\marginpar{\textcolor{Mycolor1}{1.#1}}{\textcolor{Mycolor1}{#2}}}
\newcommand{\RefB}[2]{\hspace{0pt}\marginpar{\textcolor{Mycolor2}{2.#1}}{\textcolor{Mycolor2}{#2}}}
\newcommand{\RefC}[2]{\hspace{0pt}\marginpar{\textcolor{Mycolor3}{3.#1}}{\textcolor{Mycolor3}{#2}}}

\renewcommand{\RefA}[2]{{#2}}
\renewcommand{\RefB}[2]{{#2}}
\renewcommand{\RefC}[2]{{#2}}

\makeatletter
\def\ps@pprintTitle{%
  \let\@oddhead\@empty
  \let\@evenhead\@empty
  \let\@oddfoot\@empty
  \let\@evenfoot\@oddfoot
}
\makeatother

\begin{document}

\begin{frontmatter}

\title{Solving Large-Scale Dynamic Vehicle Routing Problems with\\Stochastic Requests}

\author[tueie]{Jian Zhang}
\ead{j.zhang4@tue.nl}

\author[tuemcs]{Kelin Luo}
\cortext[mycorrespondingauthor]{Corresponding author}
\ead{k.luo@tue.nl}

\author[tueie]{Alexandre M. Florio\corref{mycorrespondingauthor}}
\ead{a.de.macedo.florio@tue.nl}

\author[tueie]{Tom Van Woensel}
\ead{t.v.woensel@tue.nl}

\address[tueie]{Department of Industrial Engineering and Innovation Sciences,\\ Eindhoven University of Technology, 5600MB Eindhoven, Netherlands}
\address[tuemcs]{Department of Mathematics and Computer Science,\\ Eindhoven University of Technology, 5600MB Eindhoven, Netherlands}

\begin{abstract}
Dynamic vehicle routing problems (DVRPs) arise in several applications such as technician routing, meal delivery, and parcel shipping. We consider the DVRP with stochastic customer requests (DVRPSR), in which vehicles must be routed dynamically with the goal of maximizing the number of served requests. We model the DVRPSR as a multi-stage optimization problem, where the first-stage decision defines route plans for serving scheduled requests. Our main contributions are knapsack-based linear models to approximate accurately the expected reward-to-go, measured as the number of accepted requests, at any state of the stochastic system. These approximations are based on representing each vehicle as a knapsack with a capacity given by the remaining service time available along the vehicle's route. We combine these approximations with optimal acceptance and assignment decision rules and derive efficient and high-performing online scheduling policies. We further leverage good predictions of the expected reward-to-go to design initial route plans that facilitate serving dynamic requests. Computational experiments on very large instances based on a real street network demonstrate the effectiveness of the proposed methods in prescribing high-quality offline route plans and online scheduling decisions.
\end{abstract}


\end{frontmatter}

\section{Introduction} \label{SecIntroduction}
\RefA{1}{We consider the dynamic vehicle routing problem with stochastic requests (DVRPSR), a fundamental problem in freight transport. In the multi-vehicle DVRPSR, a set of \textit{static} requests is known in advance, while other requests (hereafter, \textit{dynamic} requests) arrive randomly during the service period. Vehicles depart from a given location, and must return to that location before the end of the service period. The goal is to serve as many customer requests as possible by assigning an initial route plan to each vehicle, and by routing vehicles dynamically in response to new requests.}

Previous research on the DVRPSR is mostly restricted to small- and medium-scale instances on synthetic graphs with Euclidean or Manhattan distances \cite[e.g.,][]{azi2012dynamic,klapp2018dynamic,ulmer2018budgeting,ulmer2019offline,voccia2019same}. This paper scales the DVRPSR up to practical sizes: we address the dynamic and real-time routing of multiple vehicles on a real street network with more than 16,000 nodes (intersections), each corresponding to a potential customer location. Real-time decision-making's tight time frames pose considerable challenges when solving the DVRPSR in our large-scale setting.

A solution to the DVRPSR consists of an initial route plan and an online scheduling policy. The former is computed \emph{offline}, i.e., before the beginning of the service period, and defines the initial vehicle routes covering all static requests. In many recently developed DVRPSR solution methods, initial route plans are computed by myopic heuristics such as cheapest insertion (CI) and savings heuristic \cite[e.g.,][]{ulmer2018budgeting,ulmer2019offline,van2019delivery}, which attempt to minimize the travel time required for serving static requests, but completely ignore the probable realization of dynamic requests. Ideally, however, the initial route plan should also facilitate the acceptance of dynamic requests during the service period. As shown by \cite{bent2004scenario} and \cite{ferrucci2016pro}, exploiting the stochastic knowledge of dynamic requests can improve the quality of the initial routes.

In addition to the initial route plan, the DVRPSR requires an online policy to quickly determine acceptance/rejection decisions, request-vehicle assignments, and updated routes. Value function approximation (VFA) and rollout algorithms are the most widely used approximate dynamic programming (ADP) methods for computing online policies \cite[e.g.,][]{klapp2018dynamic,ulmer2019offline,van2019delivery}. Both methods estimate the expected reward-to-go by simulation and can effectively solve small DVRPSR instances. However, when the problem scale becomes realistically large, the application of many VFA policies (e.g., those developed in \cite{ulmer2018budgeting} and \cite{ulmer2019offline}) is restrained by the high dimensionality of lookup tables. Rollout algorithms are unsuitable for large-scale problems as well, since they perform all simulations online and their solution quality and computation time depend largely on the number of simulated scenarios.

\RefA{1}{Our key methodological contributions are knapsack-based linear models to approximate the expected reward-to-go at any state of the dynamic system. We derive these models by representing each vehicle as a knapsack whose capacity is given by the remaining available service time (or \emph{budget}) along the vehicle's route. We assess the cost (or budget consumption) of assigning future requests to a vehicle by predicting the vehicle's future location along its planned route. We combine these approximations with the decision rules of an optimal policy in order to derive high-performing and efficient scheduling policies. Moreover, we demonstrate how these approximations can be used to select initial route plans that facilitate the acceptance of requests in the dynamic phase.}

\begin{figure}[ht]
    \begin{center}
	    \caption{A Large-scale DVRPSR Instance.\label{FigSimulator}}
        \includegraphics[trim={200 200 50 200},clip,width=0.85\textwidth]{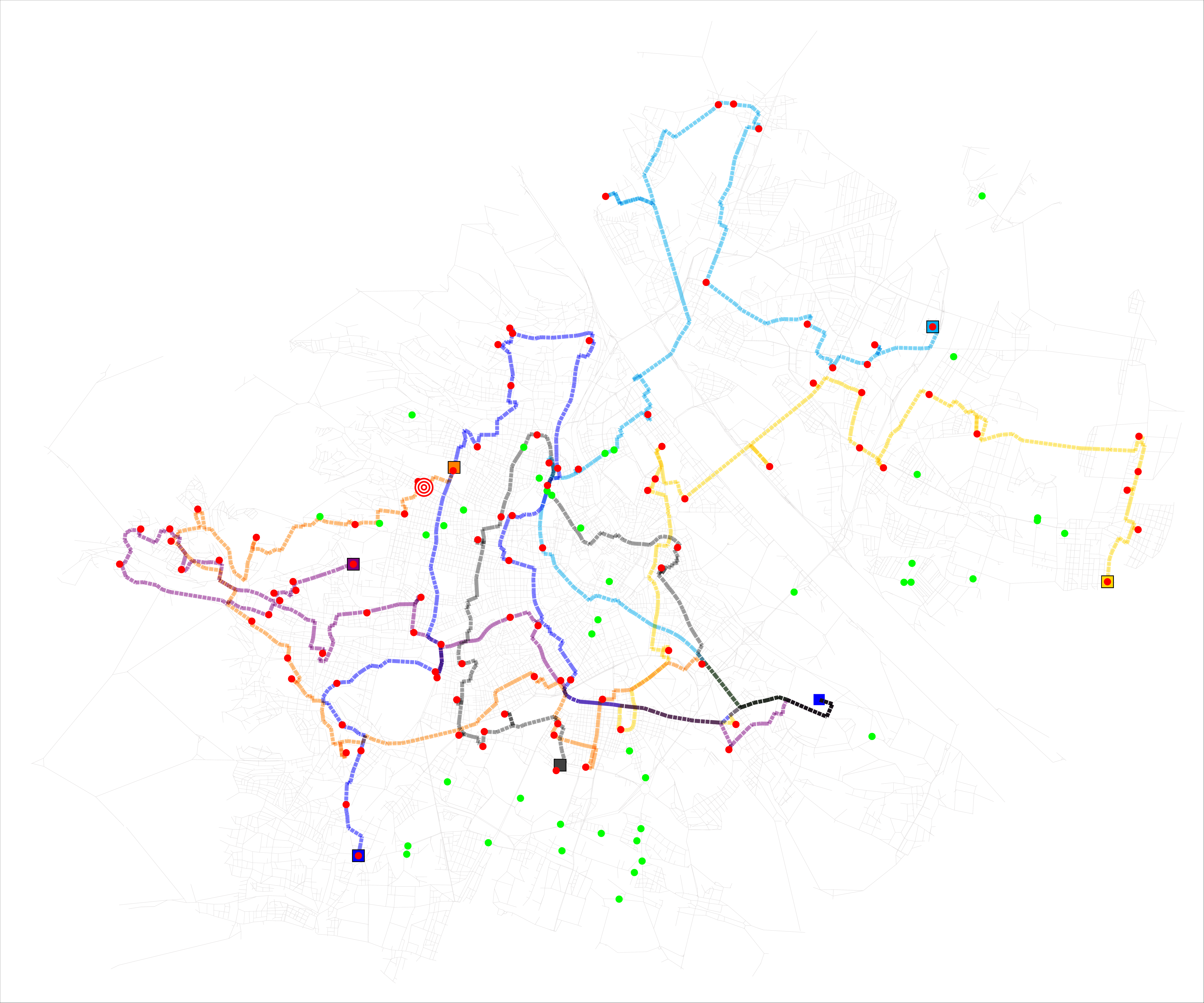}
    \end{center}
	\footnotesize\textit{Notes.} The green dots represent the served requests; the red dots represent the accepted requests that have not yet been served; the red circle represents a dynamic request that has just arrived; the blue square in the southeast represents the depot; the other squares represent vehicles' current locations and the lines represent the remaining planned routes (each vehicle is associated with a specific color).
\end{figure}

To summarize, the main contributions of this paper are the following:

\begin{enumerate}
\item We propose a sequential optimization model for the DVRPSR, which jointly optimizes the offline decision concerning the initial route plan and the online decisions for scheduling dynamic requests. This is the first DVRPSR model that incorporates, at the same time, multiple vehicles, a real street network, and real-time decision-making. We further characterize the decision rules of optimal scheduling policies.
\item We develop interpretable knapsack-based approximations of the expected reward-to-go, which accurately and efficiently predict the expected number of future accepted requests from any state in our DVRPSR model. These approximations improve the service area coverage achieved by the initial route plan and, in combination with the optimal decision rules, lead to high-performing and efficient online scheduling policies.
\item We demonstrate by an extensive computational study on large-scale instances the merits of initial route plans that consider the spatiotemporal distribution of dynamic requests. Further, we show that the scheduling policies based on the knapsack approximations are suitable for real-time decision-making and, at the same time, outperform conventional ADP methods widely applied to solve dynamic vehicle routing problems (DVRPs).
\item We make freely available the implementation of a generic DVRPSR simulator (see Figure~\ref{FigSimulator} for a snapshot and \href{https://youtu.be/D57xNfU73as}{https://youtu.be/D57xNfU73as} for animated simulations) and a set of standard benchmark instances to allow reproducibility and comparison of results.
\end{enumerate}

The remainder of this paper is organized as follows. Section~\ref{SecLiterature} reviews the relevant literature. Section~\ref{SecProblem} introduces a sequential optimization model for the DVRPSR. Section~\ref{SecOnline} and Section~\ref{SecOffline} elaborate on online scheduling policies and offline planning algorithms, respectively. The results of the computational study are presented and discussed in Section~\ref{SecComputation}. Finally, Section~\ref{SecConclusions} draws concluding remarks and proposes directions for further research.

\section{Literature Review}\label{SecLiterature}
A detailed review of contributions to the DVRPSR is presented in Section 2.1. The applications of ADP to DVRPs are discussed in Section 2.2.

\subsection{Dynamic Vehicle Routing with Stochastic Requests}
We first focus on the previous research on the DVRPSR. Table \ref{TabLiterature} provides an overview of the relevant contributions. In terms of problem settings, we indicate, for each paper, whether decisions are made in real-time \RefC{1.2}{(i.e., whether decision epochs are triggered by request arrivals)}, the maximum expected number of requests, the maximum number of vehicles, and the graph type (synthetic graph or street network). A ``n/a'' (not applicable) in the column ``Veh.'' means that fleet size is not constrained. Table \ref{TabLiterature} also summarizes the methodologies employed for computing initial route plans and online scheduling policies. A ``n/a'' in the column ``Initial plan'' indicates that static requests are not considered.
	
\begin{table}[!htb]
    \renewcommand{\baselinestretch}{1.2}
	\caption{DVRPSR: Classification of the Relevant Literature\label{TabLiterature}}
	{
	    \addtolength\tabcolsep{-0.1em}
		\begin{adjustbox}{center}
			\begin{threeparttable}
				\centering
				\footnotesize
				\begin{tabular}{lcccccccc}
					\hline
					\multirow{2}{*}{Literature} & \multicolumn{4}{c}{Problem setting} & \multicolumn{2}{c}{Initial plan} & \multicolumn{2}{c}{Online policy} \\\cmidrule(r){2-5}\cmidrule(lr){6-7}\cmidrule(lr){8-9}
						& RT & Req. & Veh. & Graph & Method & Ant. & Method & Ant. \\\hline
					\cite{bent2004scenario} & $\times$ & 100 & n/a & Synth. 			& MSA & \checkmark & MSA & \checkmark \\
					\cite{chen2006dynamic} & $\times$ & 100 & n/a & Synth. 				& CGBH & $\times$ & CGBH & $\times$ \\ 					
					\cite{hvattum2006solving} & $\times$ & 130 & n/a & Synth. 			& DSHH & \checkmark & DSHH & \checkmark \\
					\cite{ichoua2006exploiting} & \checkmark & 240 & 6 & Synth. 		& TS & \checkmark & TS & \checkmark \\
					\cite{thomas2007waiting} & $\times$ & 50 & 1 & Synth. 				& GRASP & $\times$ & Heuristics & \checkmark \\ 		
					\cite{azi2012dynamic} & \checkmark & 144 & 5 & Synth. 				& \multicolumn{2}{c}{\hspace{0.5cm}n/a} & ALNS & \checkmark \\						
					\cite{ferrucci2015general,ferrucci2016pro} & $\times$ & 150 & 12 & SSSN & TS & \checkmark & TS & \checkmark \\					
					\cite{klapp2018dynamic} & $\times$ & 40 & 1 & Synth.       			& B\&C & \checkmark & Rollout & \checkmark \\			
					\cite{ulmer2018budgeting} & $\times$ & 100 & 1 & Synth. & CI & $\times$ & VFA & \checkmark \\					
					\cite{ulmer2018value} & $\times$ & 100 & 1 & Synth. & CI & $\times$ & VFA & \checkmark \\					
					\cite{ulmer2019offline} & $\times$ & 100 & 1 & Synth.      			& CI & $\times$ & VFA+rollout & \checkmark  \\			
					\cite{van2019delivery} & $\times$ & 400 & n/a & Synth., SSSN 		& CW & $\times$ & VFA, rollout & \checkmark \\					
					\cite{voccia2019same} & $\times$ & 192 & 13 & Synth.         		& \multicolumn{2}{c}{\hspace{0.5cm}n/a} & MSA & \checkmark \\						
					\cite{ulmer2020dynamic} & \checkmark & 180 & 3 & Synth. 			& \multicolumn{2}{c}{\hspace{0.5cm}n/a} & VFA & \checkmark \\						
					\cite{ulmer2020meso} & \checkmark & 50 & 1 & Synth.      			& \multicolumn{2}{c}{\hspace{0.5cm}n/a} & VFA & \checkmark \\						
					\textbf{This paper} & \textbf{\checkmark} & \textbf{947} & \textbf{20} & \textbf{LSSN} & \textbf{PbCGBH} & \checkmark & \textbf{PbPs} & \textbf{\checkmark} \\	
					\hline
				\end{tabular}
			\end{threeparttable}
		\end{adjustbox}\vspace{0.5em}
	}
	\footnotesize\textit{Notes.} RT, real-time; `Req.', requests; `Veh.', vehicles; `Ant.', anticipation; `Synth.', synthetic; SSSN, small-scale street network; LSSN, large-scale street network; MSA, multiple scenario approach; CGBH, column-generation-based heuristic; DSSH, dynamic stochastic hedging heuristic; TS, tabu search metaheuristic; GRASP, greedy randomized adaptive search procedure; ALNS, adaptive large neighborhood search; B\&C, branch-and-cut; CI, cheapest insertion; CW, Clarke-Wright savings heuristic; PbCGBH, potential-based CGBH; PbPs, potential-based policies.
\end{table}

Early research on the DVRPSR proposes and evaluates several policies which follow certain decision rules (e.g., first-come, first-served) and ignore probabilistic knowledge about future requests when prescribing decisions \citep[e.g.,][]{bertsimas1991stochastic,tassiulas1996adaptive,gendreau1999parallel,larsen2002partially}. \RefC{2.1}{\cite{chen2006dynamic} develop a dynamic column generation framework to solve a DVRPSR. Assuming that no probabilistic information of future requests is available, they generate single-vehicle routes (i.e., columns) continuously and solve a set partitioning model periodically with existing columns. \cite{bent2004scenario}, \cite{hvattum2006solving}, and \cite{ichoua2006exploiting} are the earliest works to} show that exploiting stochastic information of dynamic requests leads to significantly better initial route plans and online scheduling policies. The multiple scenario approaches (MSAs) presented in \cite{bent2004scenario} and \cite{voccia2019same} compute and reoptimize route plans for multiple scenarios, each corresponding to a static VRP involving known requests and randomly sampled future requests. From multiple route plans, the optimal one is selected by consensus functions. \cite{ichoua2006exploiting} and \cite{ferrucci2015general,ferrucci2016pro} integrate randomly sampled future requests into tabu search metaheuristics for static VRPs, by which anticipatory DVRPSR policies are derived.

MSAs and anticipatory (meta-)heuristics are applied to DVRPSRs by decomposing the problem into sequential static VRPs \cite[e.g.,][]{hvattum2006solving,ferrucci2015general,ferrucci2016pro,klapp2018dynamic}. In these methods, decision epochs occur at the end of fixed time intervals. Thus, there is usually sufficient computation time for reoptimizing static VRPs, but new requests must wait until the end of the current time interval before being scheduled or rejected. To improve customer responsiveness, the decision epochs in many DVRPSRs are triggered by events, e.g., when a vehicle arrives at a node or when a new request arrives. In these cases, sequential stochastic optimization models are more suitable than sequential static VRP formulations. For example, \cite{thomas2007waiting} and \cite{ulmer2018budgeting} propose Markov decision processes for single-vehicle DVRPSRs, where a decision epoch occurs every time the vehicle stops at a customer location or returns to the depot. At each epoch, a subset of the newly arrived requests is accepted, and the vehicle's route is updated accordingly.

In most DVRPSRs, some dynamic requests may be rejected due to the limited number of vehicles and hard deadlines of service periods \cite[e.g.,][]{klapp2018dynamic,ulmer2018budgeting,voccia2019same}. It is hence necessary to consider customers' expectations for quick responses from the service provider. If a decision epoch occurs only when a vehicle stops or at predetermined time instants, most responses to dynamic requests will be delayed. \RefC{1.2}{Column ``RT'' of} Table \ref{TabLiterature} shows that only 4 models define decision epochs as the moments when new requests arrive. Among these, \cite{azi2012dynamic} spend more than 200 seconds to compute a decision for small instances with 3 vehicles and 144 expected requests, whereas in the remaining works, the computation times are not specified.
	
Concerning problem scales, most DVRPSR instances in the previous research are defined on synthetic graphs with modest numbers of requests and vehicles. In \cite{van2019delivery}, several instances based on the real street network of Copenhagen (Denmark) are solved, but the authors consider only 10 customer locations and hence the problem scale is still limited. \cite{ferrucci2015general} formulate test instances on the real street network of Dortmund (Germany), but the number of nodes (potential customer locations) or the geographical distribution of dynamic requests is not specified. Therefore, no research has addressed solving large-scale DVRPSR instances in real-time.

\subsection{Approximate Dynamic Programming in DVRPs}
Sequential stochastic optimization (also referred to as Markov decision process) is a commonly used mathematical model for DVRPs (including DVRPSRs) \citep{pillac2013review,psaraftis2016dynamic,ritzinger2016survey,soeffker2021stochastic}. However, in practice, solving these models exactly by classical dynamic programming is computationally prohibitive because of the ``curses of dimensionality". Consequently, a variety of ADP methods have been developed for DVRPs. The basic idea of ADP is to make decisions according to the estimated reward-to-go, which is generally obtained by simulation. \RefA{2}{We refer the interested reader to \cite{powell2011approximate} and \cite{powell2012approximate} for tutorials on the different ADP methods and their applications to various transportation and logistics problems.}

\cite{powell2011approximate} classifies ADP into four broad categories: myopic policies, look-ahead policies, value function approximation (VFA), and policy function approximation (PFA). Myopic policies guide the selection of decisions in such a way that immediate rewards are maximized. For DVRPs, these policies tend to perform poorly because they do not take into account any forecast information about dynamic requests. Look-ahead policies, on the contrary, approximate the reward-to-go by explicitly sampling and simulating the future. Rollout policies are a typical look-ahead policy class \RefA{2}{initially proposed by \cite{bertsekas1996neuro}}. For applications to DVRPSRs, \cite{ulmer2019offline} demonstrate that any base policy with deterministic decision rules can be improved by rollout algorithms. However, look-ahead policies have the critical shortcoming that they perform all simulations online, and hence their performance is highly restricted by the tight time frames of DVRPSRs. For example, the rollout policies employed in \cite{ulmer2019offline} simulate 16 scenarios at each decision epoch. For instances with only 100 expected requests on a Euclidean plane, the maximum computation time for selecting a decision exceeds 6 minutes, which is a far too long response time when considering real-time applications.

In contrast to look-ahead policies, VFA policies conduct all simulations offline (i.e., a priori) and hence are more suitable for the applications with limited online computation time. The main component of a VFA policy is an approximation of the value function which captures the expected reward-to-go. \RefA{2}{\cite{powell2012approximate} establish a unified VFA framework to address the very large state spaces of DVRPs. Based on this framework, \cite{ulmer2018budgeting}, \cite{van2019delivery}, and \cite{ulmer2020meso} propose non-parametric VFA (N-VFA), parametric VFA (P-VFA), and meso-parametric VFA (M-VFA) policies, respectively, to approximate value functions for DVRPSRs.} The main drawback of these VFA policies is the high computational complexity. Both N-VFA and M-VFA rely on lookup tables (completely or partially) which suffer from the curses of dimensionality \citep{powell2015tutorial}. P-VFA requires less computational effort, but its linear basis function ignores non-linearity and cannot approximate the value function in every detail \citep{ulmer2020meso}.
	
The last category of ADP policies, PFA, is similar to VFA in the sense that all simulations are conducted offline, but is easier to implement because it usually only requires tuning the value of a parameter rather than approximating a complex value function. \RefA{2}{PFA policies are suitable for the cases where the structure of a good decision policy can be easily observed and captured by an analytical function, which returns a decision given a state \citep{powell2012approximate}.} A typical PFA policy is developed by \cite{ulmer2019same} for a dynamic dispatching problem, in which every vehicle can only perform direct trips between a depot and one of several parcel lockers. At each fixed time point, the policy identifies the locker with the largest amount of parcels to be delivered. If this amount is larger than a threshold value, a vehicle is dispatched to serve that locker; otherwise, all available vehicles wait at the depot. The threshold balances the trade-off between early dispatch for fast deliveries and late dispatch for consolidation, and extensive offline simulations determine its value. \RefA{2}{For a more complex dynamic dispatching problem with both vehicles and drones, \cite{ulmer2018same} propose a PFA policy based on the intuition that it is beneficial to serve remote customers by drones and serve nearby customers by vehicles. With offline simulations, they tune a threshold to measure whether each customer is far enough from the depot so that it should be served by a drone.}

The three categories of anticipatory ADP policies (i.e., look-ahead, VFA, and PFA) require a sufficiently large number of online or offline simulations. As a result, it is hard to balance their computational efficiency and solution quality. \RefC{1.3}{To the best of our knowledge, there is no proven method} to solve accurately large-scale DVRPSR instances while respecting the tight time frames imposed by real-time decision-making. This paper fills this gap by proposing accurate and efficient approximations of the reward-to-go, which lead to practical algorithms to solve large-scale DVRPSRs in real-time. Moreover, the proposed approximations are based on sample scenarios and hence can better cope with the heterogeneity of stochastic information (e.g., time-varying request rates) than VFA \citep{soeffker2021stochastic}. With a set of PFA and rollout policies as benchmarks, we validate through an extensive computational study that the proposed algorithms outperform the conventional ADP policies in terms of both efficiency and quality.

\section{Problem Definition} \label{SecProblem}
In this section, we formally introduce the DVRPSR. Section \ref{gendef} presents general definitions, and Section \ref{ssom} models the DVRPSR as a sequential stochastic optimization problem. \RefA{4}{The notations used in this section are summarized in Table \ref{tab:notation_pro} in \ref{sec:notations}.}

\subsection{General Definitions} \label{gendef}
The DVRPSR is defined on a strongly connected digraph $\mathcal{G}=(\mathcal{V}, \mathcal{A})$, where $\mathcal{V}=\{0,\ldots,V\}$ is the set of nodes and $\mathcal{A}$ is the set of arcs. Graph $\mathcal{G}$ represents the street network of the service area, and sets $\mathcal{V}$ and $\mathcal{A}$ represent, respectively, road intersections and road segments. We assume that all potential customers are located at intersections. The length of each segment $(i,j)\in\mathcal{A}$ is given by $d_{ij}$. The service period is represented by the interval $[0,U]$. We use interchangeably the terms \emph{time} and \emph{instant} to refer to a value $u\in[0,U]$. Before the service period begins, a fleet of $K$ vehicles is stationed at node 0, which corresponds to the location of the depot. Vehicles travel at a constant speed $\overline{s}$, so the travel time along a segment $(i,j)$ is given by $t_{ij}=d_{ij}/\overline{s}$. Further, the duration (under speed $\overline{s}$) of the fastest path from node $i$ to node $j$ is given by $t(i,j)$.

A customer request consists of a triple $(u,i,d)$, where $u\in[0,U]$ is the request arrival time (alternatively, the time when the request is \emph{generated} by the customer), $i\in\mathcal{V}\setminus\{0\}$ is the customer location, and $d\in\mathbb{R}_{>0}$ is the \emph{duration} of the request, i.e., the service time required by the customer to fulfill the request. A set $\mathcal{S}$ of static requests is known before the beginning of the service period. Static requests must be served by the end of the service period. In addition, dynamic requests arrive randomly during the service period according to a non-homogeneous Poisson process with a request rate $\Lambda(u)=\sum_{i\in\mathcal{V}\setminus\{0\}}\lambda_{i}(u)$, where $\lambda_{i}(u)$ are node specific request rates. We denote by $\mathcal{D}(u)$ the ordered (by arrival time) set of dynamic requests generated up to instant $u$. Hence, $\{|\mathcal{D}(u)|\}_{u\in[0,U]}$ is the counting process associated with the request arrival process $\{\Lambda(u)\}_{u\in[0,U]}$, so $|\mathcal{D}(u)|$ is a random variable indicating the number of dynamic requests generated up to instant $u$. Finally, the duration of a dynamic request is known at the moment it arrives, and so are the durations of all static requests. The durations of future (not yet generated) requests follow a known probability distribution.

Once a dynamic request arrives, the scheduling policy in effect decides whether to accept or to reject the request. Accepted requests must be served before the end of the service period, and rejected requests are lost (it is assumed that rejected requests are served by competitors). The scheduling policy must take a decision in an online fashion, that is, almost immediately after a dynamic request arrives. A static or accepted dynamic request $(u,i,d)$ is assigned to a vehicle $k\in\{1,\ldots,K\}$. The request is served by routing vehicle $k$ to location $i$ and keeping the vehicle there for an amount of time $d$. Once a request is served, the vehicle proceeds to the location of another request, or back to the depot if there are no further requests assigned to that vehicle. In particular, waiting is not allowed, except when a vehicle is idle at the depot. The goal in the DVRPSR is to route the $K$ vehicles, so that all requests in $\mathcal{S}$ are served, all vehicles return to the depot no later than instant $U$, and the number of dynamic requests accepted and served is maximized.

\subsection{A Sequential Stochastic Optimization Model for the DVRPSR} \label{ssom}
We now propose a model for the DVRPSR inspired by the unified ADP framework of \cite{powell2012approximate} \RefA{3}{and by the route-based Markov decision processes proposed by \cite{ulmer2020modeling} for modeling DVRPs}.

\subsubsection{Decision Epochs and States}
Decision epochs correspond to the moments when decisions must be taken. In the DVRPSR, in addition to a first decision concerning the initial route plan to serve requests in $\mathcal{S}$, a decision is also made whenever a dynamic request arrives. Hence, there are in total $1+|\mathcal{D}(U)|$ decision epochs, where $|\mathcal{D}(U)|$ is uncertain.

The state of the system at the moment a decision must be made is represented by state variables $S_{t}$, $t\in\{0,\ldots,T\}$, where $T=|\mathcal{D}(U)|$ is the last decision epoch. State $S_{0}$ represents all parameters of the problem, since no uncertainty has been realized before the beginning of the service period. Next, we formally define routes, which are instrumental for characterizing all other states $S_{t}$, $t\neq0$:
\begin{definition}[Route]
A route is a sequence of pairs $\theta=((v_{0},\delta_{0}),\ldots,(v_{f},\delta_{f}))$, where $v_{0}=v_{f}=0$, $(v_{i},v_{i+1})\in\mathcal{A}$, $i\in\{0,\dots,f-1\}$, $\delta_{0}=\delta_{f}=\emptyset$, and $\delta_{i}$, $i\in\{1,\ldots,f-1\}$, is the (possibly empty) set of requests to serve when arriving at the $i$-th node in $\theta$. For notational convenience, $\Delta(\theta)=\cup_{i=1}^{f-1}\delta_{i}$ denotes the set of all requests served by $\theta$, and $\Delta_{u}(\theta)=\{(u',i,d)\in\Delta(\theta):u'\geq u\}$ denotes the ordered (by $u'$) set of all requests served by $\theta$ from instant $u$ onwards.
\end{definition}

Given a route $\theta=((v_{0},\delta_{0}),\ldots,(v_{f},\delta_{f}))$ and a starting time $\tau\in[0,U]$, the \emph{budget} $b(\tau,\theta)$ of route $\theta$ when it starts at instant $\tau$ is defined as the slack time relative to the end of the service period:
\begin{equation} \nonumber
b(\tau,\theta)=U-\tau-\sum_{i=0}^{f-1}t_{v_{i},v_{i+1}}-\sum_{(u,i,d)\in\Delta(\theta)}d.
\end{equation}

A state $S_{t}$, $t\in\{1,\ldots,T\}$, is a $(K+1)$-tuple $S_{t}=(V_{1},\ldots,V_{K},r_{t})$, where $V_{k}$ are vehicle states, and $r_{t}=(u_{t},i_{t},d_{t})$ is the $t$-th element of $\mathcal{D}(U)$. A vehicle may be either in an idle state, in which case we say that $V_{k}=\emptyset$, or in a routing state executing a route $\theta_{k}$, in which case we say that $V_{k}=(\tau_{k},\theta_{k})$, where $\tau_{k}$ is the instant when vehicle $k$ started route $\theta_{k}$ and $b(\tau_{k},\theta_{k})\geq0$. \RefC{2.2}{As $u_{t}$, $d_{t}$, and $\tau_{k}$ are continuous variables, the number of possible states is infinite.} \RefA{3}{The idea of incorporating route information in the state space was originally proposed by \cite{ulmer2020modeling}, and is particularly useful as our solution methods rely heavily on information about planned routes to guide online decisions.} So, when in a routing state, a vehicle must be either in transit (traveling along a road segment) or in service (fulfilling a request at a customer).

Since travel times and service times of accepted requests are deterministic, a state collects sufficient information to define precisely the locations of all vehicles. In particular, given a state it is possible to identify the status of each vehicle (idle, in transit or in service), to compute the remaining time before a vehicle reaches the endpoint of the road segment along which it currently travels (if the vehicle is in transit), or to compute the remaining time that a vehicle will stay at a particular node serving requests (if the vehicle is in service).

\subsubsection{Offline and Online Decisions}
It is convenient to discuss separately the decision at $S_{0}$ and the decisions at other states $S_{t}$, $t\neq0$, since they are structurally different. We denote by $\mathcal{F}$ the feasible decision set at state $S_{0}$. Set $\mathcal{F}$ contains all sets of routes $\{\theta_{1},\ldots,\theta_{k'}\}$, $k'\leq K$, such that $b(0,\theta_{1}),\ldots,b(0,\theta_{k'})\geq0$ and $\Delta(\theta_{1}),\ldots,\Delta(\theta_{k'})$ is a partition of $\mathcal{S}$. We represent generically the decision made at $S_{0}$ by a decision vector $\mathbf{y}\in\mathcal{F}$. This decision corresponds to an initial route plan and is computed offline (i.e., before the start of the delivery period) by, e.g., one of the algorithms described in Section \ref{SecOffline}.

We now consider a state $S_{t}=(V_{1},\ldots,V_{K},r_{t})$, $t\neq0$. That is, a state where the dynamic request $r_{t}=(u_{t},i_{t},d_{t})$ has just arrived and an online decision must be taken, in real-time. Given $S_{t}$, a \emph{scheduling policy} is the set of rules that determines:

\begin{enumerate}[label=(\roman*)]
\item An \emph{acceptance} decision: this is an `accept' or `reject' decision relative to request $r_{t}$.
\item An \emph{assignment} decision: if the request is accepted, this decision allocates the request to a vehicle $k\in\{1,\ldots,K\}$.
\item A \emph{routing} decision: when the request is allocated to vehicle $k$, this decision defines how route $\theta_{k}$ will be adjusted (in case $V_{k}=(\tau_{k},\theta_{k})$) or initialized (in case $V_{k}=\emptyset$) to accommodate the accepted request.
\end{enumerate}

There are $1+K$ combinations of acceptance and assignment decisions, but in general there are many more possibilities concerning routing decisions, especially when considering the road network graph. The set of possible decisions when the system is at state $S_{t}$ is denoted by $\mathcal{X}_{t}$. Each decision $x_{t}\in\mathcal{X}_{t}$, $t\neq0$, specifies acceptance and, when required, also assignment and routing decisions. Formally, a (deterministic) scheduling policy maps every possible state $S_{t}$, $t\neq0$, to a decision $x_{t}\in\mathcal{X}_{t}$. We further denote by $X^{\pi}(S_{t})$ the decision prescribed by policy $\pi$ when at state $S_{t}$. Once a decision is made at state $S_{t}$, $t\neq T$, the system evolves and the next decision epoch occurs when the $(t+1)$-th request arrives. State $S_{t+1}$ is specified by the preceding state $S_{t}$, the decision $x_{t}\in\mathcal{X}_{t}$, and the $(t+1)$-th request $r_{t+1}=(u_{t+1},i_{t+1},d_{t+1})$. This functional relationship is represented by the transition function $S^{M}(\cdot)$:
\begin{equation} \nonumber
S_{t+1}=S^{M}(S_{t},x_{t},r_{t+1}).
\end{equation}

\subsubsection{Rewards and Objective Function}
We denote by $R(S_{t},x_{t})$ the reward associated with decision $x_{t}\in\mathcal{X}_{t}$ when at state $S_{t}$. At state $S_{0}$, any decision $\mathbf{y}\in\mathcal{F}$ has a reward of 0. At any other state $S_{t}$, $t\neq0$, the reward of a decision is 1 if the dynamic request $(u_{t},i_{t},d_{t})$ is accepted, and 0 otherwise.

The previous definitions allow us to propose the following sequential optimization model (SOM) for the DVRPSR:
\begin{equation} \label{model:som}
\max_{\pi\in\Pi,\mathbf{y}\in\mathcal{F}}\mathbb{E}\Bigg[\mathbb{E}\bigg[\sum_{t=1}^{T}R(S_{t},X^{\pi}(S_{t}))\bigg| S_{1}\bigg]\Bigg],\tag{DVRP-SOM}
\end{equation}
where $\Pi$ is the set of all scheduling policies. The double expectation is required because state $S_{1}=S^{M}(S_{0},\mathbf{y},r_{1})$ is unknown, since it depends on the (also unknown) first dynamic request $r_{1}=(u_{1},i_{1},d_{1})$.

Model \eqref{model:som} maximizes the expected total reward (or the expected number of accepted requests) conditional on the state resulting from the initial decision $\mathbf{y}$. In this work, we are concerned with finding a good policy $\pi$ as well as a good initial plan $\mathbf{y}$ to be used in combination with $\pi$. The online and offline algorithms that we propose in Sections \ref{SecOnline} and \ref{SecOffline}, respectively, are based on the idea of the \emph{potential} $\Phi_{\pi}(S_{t})$ of a state $S_{t}$, $t\neq0$, which is defined as the expected reward-to-go once the system occupies state $S_{t}$ and is controlled by policy $\pi$:
\begin{equation} \label{eq:potential}
\Phi_{\pi}(S_{t})=\mathbb{E}\left[\sum_{t'=t}^{T}R(S_{t'},X^{\pi}(S_{t'}))\middle| S_{t}\right],\qquad\qquad t\neq0.
\end{equation}

When policy $\pi$ is relatively simple, \eqref{eq:potential} can be estimated in a reasonable time by simulation and sample average approximation. An example of such a simple policy is a greedy policy that accepts a request whenever it is feasible to serve it, and assigns requests to vehicles according to cheapest insertion. We further discuss this policy (and its rollout version) in Section \ref{sec:bench}. In the following section, we describe computationally efficient approximations of $\Phi_{\pi}(S_{t})$ when policy $\pi$ is \RefB{1}{anticipatory}, that is, when it takes into account the spatiotemporal distribution of future requests and prescribes acceptance and assignment decisions accordingly.

\section{Online Scheduling Policies} \label{SecOnline}
In this section, we develop several scheduling policies. Across all policies, routing decisions are prescribed according to one of the two routing policies defined in Section \ref{ss:routingpcy}. Section \ref{ss:charact} characterizes acceptance and assignment decisions in the optimal scheduling policy. Next, in Section \ref{ss:mka} we derive our potential approximation model following a sequence of interpretable steps. Sections \ref{sec:pbp} and \ref{sec:spbp} introduce two potential-based scheduling policies, and Section \ref{sec:bench} presents greedy and rollout benchmark policies. \RefA{4}{All notations used in this section are summarized in Table \ref{tab:notation_alg} in \ref{sec:notations}.}

\subsection{Routing Policies}\label{ss:routingpcy}
When a dynamic request is accepted and assigned to a vehicle in a routing state, it is reasonable to schedule the new request to be served in a position such that the budget of the adjusted route is maximized. Depending on whether it is allowed to reorder unserved requests, this leads to two main routing policies, which are defined next:
\begin{definition}[Cheapest Insertion (CI) Routing Policy $\rho_{_\mathsf{CI}}$]
Consider a request $r=(u,i,d)$ and a vehicle state $V=(\tau,\theta)$. Under the CI routing policy $\rho_{_\mathsf{CI}}$, route $\theta$ is adjusted into a new route $\rho_{_\mathsf{CI}}(\theta,r)$ in such a way that $\Delta_{u}(\rho_{_\mathsf{CI}}(\theta,r))\setminus\Delta_{u}(\theta)=\{r\}$, $\Delta_{u}(\theta)$ is a subsequence of $\Delta_{u}(\rho_{_\mathsf{CI}}(\theta,r))$ and the budget $b(\tau,\rho_{_\mathsf{CI}}(\theta,r))$ is maximized.
\end{definition}

\begin{definition}[Reoptimization Routing Policy $\rho_{_\mathsf{R}}$]
Consider a request $r=(u,i,d)$ and a vehicle state $V=(\tau,\theta)$. Under the reoptimization routing policy $\rho_{_\mathsf{R}}$, route $\theta$ is adjusted into a new route $\rho_{_\mathsf{R}}(\theta,r)$ in such a way that $\Delta_{u}(\rho_{_\mathsf{R}}(\theta,r))\setminus\Delta_{u}(\theta)=\{r\}$ and the budget $b(\tau,\rho_{_\mathsf{R}}(\theta,r))$ is maximized.
\end{definition}

Given a route $\theta$ and a request $r$, $\rho_{_\mathsf{CI}}(\theta,r)$ can be computed in polynomial time on $|\theta|$, while computing $\rho_{_\mathsf{R}}(\theta,r)$ requires solving a traveling salesman problem with $|\Delta_{u}(\theta)|+1$ nodes.

\subsection{Characterizing the Optimal Scheduling Policy}\label{ss:charact}
Given a first stage decision $\mathbf{y}$, the optimal scheduling policy $\pi^{*}$ is such that:
\begin{equation} \nonumber
\pi^{*}=\argmax_{\pi\in\Pi}\,\mathbb{E}\left[\Phi_{\pi}(S^{M}(S_{0},\mathbf{y},r_{1}))\right].
\end{equation}

We now define a class of potential-based policies $\Pi_{P}\subset\Pi$ such that acceptance and assignment decisions are prescribed according to the following decision rules (DRs):

\begin{decisionrule}[Acceptance Rule]Consider that the system is at state $S_{t}$, and let $x_{t}^{-}$ be the `reject' decision with respect to request $r_{t}$. Then, a policy $\pi\in\Pi_{P}$ accepts request $r_{t}$ if and only if there exists an `accept' decision $x_{t}^{+}\in\mathcal{X}_{t}\setminus\{x_{t}^{-}\}$ such that:
\begin{equation} \label{eq:dr1}
1+\mathbb{E}\left[\Phi_{\pi}(S^{M}(S_{t},x_{t}^{+},r_{t+1}))\right]\geq\mathbb{E}\left[\Phi_{\pi}(S^{M}(S_{t},x_{t}^{-},r_{t+1}))\right].\tag{DR1}
\end{equation}
\end{decisionrule}

\begin{decisionrule}[Assignment Rule]Consider that the system is at state $S_{t}$, and let $\{\mathcal{X}_{t}^{k}\}$ be a partition of $\mathcal{X}_{t}\setminus\{x_{t}^{-}\}$, in which $\mathcal{X}_{t}^{k}$ is the set of decisions where $r_{t}$ is accepted and assigned to vehicle $k$. Further, for all $\mathcal{X}_{t}^{k}\neq\emptyset$, let
\begin{equation}\nonumber
x_{t}^{k}=\argmax_{x_{t}\in\mathcal{X}_{t}^{k}}\mathbb{E}\left[\Phi_{\pi}(S^{M}(S_{t},x_{t},r_{t+1}))\right].
\end{equation}

Then, a policy $\pi\in\Pi_{P}$ assigns an accepted request $r_{t}$ to the vehicle $k^{*}$ such that:
\begin{equation} \label{eq:dr2}
k^{*}=\argmax_{k}\mathbb{E}\left[\Phi_{\pi}(S^{M}(S_{t},x_{t}^{k},r_{t+1}))\right].\tag{DR2}
\end{equation}
\end{decisionrule}

\RefB{1}{Essentially, \eqref{eq:dr1} determines that request $r_{t}$ is accepted only if the immediate reward offsets the expected decrease in potential (due to budget consumption) at the (random) state $S_{t+1}$, and \eqref{eq:dr2} determines that, if $r_{t}$ is accepted, it is assigned to the vehicle such that the expected potential at $S_{t+1}$ is maximized.} Note that $\pi^{*}\in\Pi_{P}$, otherwise $\pi^{*}$ is not optimal. Intuitively, good approximations of $\mathbb{E}[\Phi_{\pi^{*}}(S_{t+1})|x_{t}]$, $x_{t}\in\mathcal{X}_{t}$, can be used to prescribe sensible decisions according to \eqref{eq:dr1} and \eqref{eq:dr2}. Next, we discuss how to efficiently compute such approximations.

\subsection{Multiple-knapsack Approximation of $\mathbb{E}[\Phi_{\pi^{*}}(S_{t+1})|x_{t}]$}\label{ss:mka}
Given a request $r_{t}=(u_{t},i_{t},d_{t})$, our goal is to approximate $\mathbb{E}[\Phi_{\pi^{*}}(S_{t+1})|x_{t}]$ for any decision $x_{t}\in\mathcal{X}_{t}$, where $S_{t+1}=(V_{1},\ldots,V_{K},r_{t+1})$ is random since request $r_{t+1}$ is unknown at the current instant $u_{t}$. To this end, we propose an approximation based on a multiple-knapsack model where each vehicle state $V_{k}=(\tau_{k},\theta_{k})$ is interpreted as a knapsack with capacity $b(\tau_{k},\theta_{k})$. The objective is to maximize the number of future requests accepted. The cost (in terms of budget consumption) of serving a future request depends on each vehicle state (or knapsack). In particular, it depends on the vehicle location and the remaining route at the moment the future request arrives. The procedure to predict the location of a vehicle at a future instant is based on the following two observations:

\begin{remark}[Late Depot Arrival]\label{rem:lda}
Since the objective in the DVRPSR is to maximize the number of accepted requests, vehicles will usually return to the depot at an instant near the end of the service period.
\end{remark}

\begin{remark}[Request Order Preservation]\label{rem:sop}
When following the CI routing policy, the relative order of scheduled requests along a route does not change when a new request is assigned to the route. When following the reoptimization policy, the relative order of already scheduled requests usually changes only slightly, if at all, when a new request is assigned.
\end{remark}

\RefA{5}{Note that our goal of approximating the expected number of requests that can be served from a given state is related, in an inverse or symmetric fashion, to the problem of approximating the length of the minimum-length tour to visit a given set of nodes. The latter is studied in the stream of literature on routing approximations \citep{daganzo1984length,daganzo1984distance}. In our case, the route lengths are given, and we estimate the maximum number of additional requests that can be accommodated along the planned routes.}

Remarks \ref{rem:lda} and \ref{rem:sop} suggest the concept of \emph{effective speed}, which is useful to predict the location of a vehicle at future instants. The effective speed is defined as the traveling speed such that a route finishes exactly at instant $U$, assuming no further requests are added to the route. \RefA{4}{We denote by $\check{s}_{u_{t}}(V_{k})$ the effective speed of a vehicle with state $V_{k}=(\tau_{k},\theta_{k})$ at instant $u_{t}$ (the formula for computing $\check{s}_{u_{t}}(V_{k})$ is given in \ref{sec:algpredict}). Then,} it is reasonable to predict the location of vehicle $k$ at a future instant $u'>u_{t}$ by, from the current location, following route $\theta_{k}$ at speed $\check{s}_{u_{t}}(V_{k})$ for an amount of time $u'-u_{t}$. By predicting at the current instant $u_{t}$ the location of vehicle $k$ at a future instant $u'$, we also estimate the set of nodes along $\theta_{k}$ to be traversed after instant $u'$, which we denote by $\mathcal{V}_{u_{t}}(V_{k},u')$. The procedure to predict the vehicle location and determine $\mathcal{V}_{u_{t}}(V_{k},u')$ requires a simple simulation of $\theta_{k}$ and is detailed in \ref{sec:algpredict}. Set $\mathcal{V}_{u_{t}}(V_{k},u')$ is then used to estimate the cost (or budget consumption) $c_{u_{t}}(V_{k},r')$ when assigning a future request $r'=(u',i',d')$ to a vehicle with state $V_{k}$:
\begin{equation}\nonumber
c_{u_{t}}(V_{k},r')=\min_{j\in\mathcal{V}_{u_{t}}(V_{k},u')}t(j,i')+t(i',j)+d'.
\end{equation}

We are interested, at the current instant $u_{t}$, in approximating $\mathbb{E}[\Phi_{\pi^{*}}(S_{t+1})|x_{t}]$. Note that $S_{t+1}=S^{M}(S_{t},x_{t},r_{t+1})=(V_{1},\ldots,V_{K},r_{t+1})$, so vehicle states $V_{1},\ldots,V_{K}$ are known and already take into account decision $x_{t}$. Let $\omega=\{r_{1}^{\omega},\ldots,r_{T_{\omega}}^{\omega}\}$ be a sample path of the request arrival process $\{\Lambda(u)\}_{u\in(u_{t},U]}$, and $\overline{K}=\{k\in\{1,\ldots,K\}:V_{k}\neq\emptyset\}$. Then, at time $u_{t}$, after decision $x_{t}$ is taken, the number of requests in $\omega$ accepted by the optimal scheduling policy is approximated by the following model with variables $z_{kr}$, $k\in\overline{K}$ and $r\in\omega$:
\begin{align}
&\phi_{\pi^{*}}^{\omega}(S_{t+1}|x_{t})=\max&&\sum_{k\in\overline{K}}\sum_{r\in\omega}z_{kr},\nonumber\\
&\text{s.t.}&&\sum_{r\in\omega}c_{u_{t}}(V_{k},r)z_{kr}\leq b(\tau_{k},\theta_{k}),&{k\in\overline{K}},\label{eq:kpctr}\\
&&&\sum_{k\in\overline{K}}z_{kr}\leq1,&r\in\omega,\label{eq:reqa}\\
&&&0\leq z_{kr}\leq1,&k\in\overline{K},r\in\omega.\label{model:mka}\tag{MKA}
\end{align}

Model \eqref{model:mka} captures the multiple-knapsack structure of the online phase of the DVRPSR. \RefB{4}{The objective function measures the total number of requests in $\omega$ accepted and assigned to a vehicle. For each vehicle, constraints \eqref{eq:kpctr} limit the total budget consumption due to assigned requests. These constraints consider the predicted vehicle locations, as determined by effective speeds and planned routes,} and charge a vehicle-dependent cost for each accepted request, based on scheduling the request in a position such that budget consumption is minimum. \RefB{4}{Constraints \eqref{eq:reqa} forbid assigning a request to more than one vehicle. Idle vehicles are ignored in model \eqref{model:mka}, since it is not possible to predict their future locations with the concept of effective speed. Unless the fleet size is over-dimensioned, all vehicles remain busy during the service period and are considered in the model.}

We consider continuous instead of binary variables in \eqref{model:mka}, so we are able to solve the model fast in a real-time context where computational time is scarce. This is also justified since, first, the coefficients of the formulation have already been approximated, so solving \eqref{model:mka} with integer variables would not necessarily produce a better approximation; second, $|\omega|$ is usually much larger than $|\overline{K}|$, so the difference between the optimal solution values of \eqref{model:mka} and its integer counterpart is relatively small. These observations lead us to hypothesize that the approximation attained by \eqref{model:mka} is sensible, which is confirmed empirically. To achieve stability, when approximating $\mathbb{E}[\Phi_{\pi^{*}}(S_{t+1})|x_{t}]$ we take the average of solution values obtained over a set $\Omega=\{\omega_{1},\ldots,\omega_{H}\}$ of $H$ sample paths:
\begin{equation}\label{eq:mka}
\mathbb{E}[\Phi_{\pi^{*}}(S_{t+1})|x_{t}]\approx\hat{\Phi}_{\pi^{*}}(S_{t+1}|x_{t})=\frac{1}{H}\sum_{\omega\in\Omega}\phi_{\pi^{*}}^{\omega}(S_{t+1}|x_{t}).
\end{equation}

Finally, note that our approximation scheme requires no parameters other than $H$.

\subsection{Potential-based Policy (PbP)}\label{sec:pbp}
The PbP scheduling policy consists in approximating state potentials according to \eqref{eq:mka} and applying decision rules \eqref{eq:dr1} and \eqref{eq:dr2}. Given a routing policy $\rho$ and a state $S_{t}=(V_{1},\ldots,V_{K},r_{t})$, $r_{t}=(u_{t},i_{t},d_{t})$, the policy works as follows:

\begin{enumerate}[label=(\roman*)]
\item If there is an idle vehicle $k$ and it can serve request $r_{t}$, that is, $u_{t}+t(0,i_{t})+d_{t}+t(i_{t},0)\leq U$, then $r_{t}$ is accepted, assigned to vehicle $k$, and the state of vehicle $k$ changes to $V_{k}=(u_{t},\theta_{k})$, where $\theta_{k}$ corresponds to the fastest route from the depot to node $i_{t}$, and back.
\item Otherwise, a set $\tilde{\mathcal{X}}_{t}=\{\tilde{x}_{t}^{-}\}$ of candidate decisions is initialized, where $\tilde{x}_{t}^{-}$ is the `reject' decision related to request $r_{t}$.
\item For each non-idle vehicle $k$, $V_{k}=(\tau_{k},\theta_{k})$, we compute the adjusted route $\rho(\theta_{k},r_{t})$. If $b(\tau_{k},\rho(\theta_{k},r_{t}))\geq0$, then the decision $\tilde{x}_{t}^{k}$ of accepting and assigning $r_{t}$ to $k$ under routing policy $\rho$ is feasible and is added to $\tilde{\mathcal{X}}_{t}$.
\item A decision is selected from $\tilde{\mathcal{X}}_{t}$ by evaluating decision rules \eqref{eq:dr1} and \eqref{eq:dr2} based on approximations $\hat{\Phi}_{\pi^{*}}(S_{t+1}|\tilde{x}_{t})$ of $\mathbb{E}[\Phi_{\pi^{*}}(S_{t+1})|\tilde{x}_{t}]$, $\tilde{x}_{t}\in\tilde{\mathcal{X}}_{t}$, computed by \eqref{eq:mka}.\label{st:evaldr}
\end{enumerate}

Different variants of the PbP arise, depending on the adopted routing policy. We denote by PbP$_{_\textsf{R}}(H)$ the PbP under routing policy $\rho_{_\mathsf{R}}$, when $H$ sample paths are used to approximate state potentials.

\subsection{Simplified Potential-based Policy (S-PbP)}\label{sec:spbp}
\RefB{1}{Step \ref{st:evaldr} of policy PbP requires solving model \eqref{model:mka} under $H$ sample paths for each non-idle vehicle that can serve the new request. Additionally, \eqref{model:mka} is solved $H$ times to evaluate the potential of the state resulting from the reject decision. Hence, in policy PbP, model} \eqref{model:mka} is created and solved up to $H(K+1)$ times before a decision is taken. In very large instances, this may require more computational time than what is available in a real-time setting. An alternative approximation can be computed faster by considering single-knapsack models. Consider a state $S_{t}=(V_{1},\ldots,V_{K},r_{t})$, $r_{t}=(u_{t},i_{t},d_{t})$. Given a non-idle vehicle $k$ with $V_{k}=(\tau_{k},\theta_{k})$ and a sample path $\omega=\{r_{1}^{\omega},\ldots,r_{T_{\omega}}^{\omega}\}$ of process $\{\Lambda(u)\}_{u\in(u_{t},U]}$, we define the following model with variables $z_{r}$, $r\in\omega$:
\begin{align}
&p_{u_{t}}^{\omega}(V_{k})=\max&&\sum_{r\in\omega}z_{r},\nonumber\\*
&\text{s.t.}&&\sum_{r\in\omega}c_{u_{t}}(V_{k},r)z_{r}\leq b(\tau_{k},\theta_{k}),\nonumber\\*
&&&0\leq z_{r}\leq1,&r\in\omega.\label{model:ka}\tag{KA}
\end{align}

Model \eqref{model:ka} overestimates the number of requests from $\omega$ that will be assigned to vehicle $k$, since it does not take into account other vehicles that may compete for the same requests. \RefA{4}{To compensate for this bias, we define a coefficient $\alpha^{\omega}$ as the ratio between the potential estimated by model \eqref{model:mka} and model \eqref{model:ka} before accepting or rejecting request $r_{t}$:}
\begin{equation}\label{eq:alpha}
\alpha^{\omega}=\frac{\phi_{\pi^{*}}^{\omega}(S_{t+1}|\tilde{x}_{t}^{-})}{\sum_{V_{k}\neq\emptyset}p_{u_{t}}^{\omega}(V_{k})}.
\end{equation}

Ratio $\alpha^{\omega}$ falls within the unit interval, and is \RefB{1}{inversely proportional to the amount of} competition among vehicles to serve requests. Then, another approximation of $\mathbb{E}[\Phi_{\pi^{*}}(S_{t+1})|\tilde{x}_{t}^{k}]$, where $\tilde{x}_{t}^{k}$ is an `accept' decision that assigns request $r_{t}$ to vehicle $k$ and changes its state to $V_{k}^{+}$, is given by:
\begin{equation}\label{eq:ka}
\mathbb{E}[\Phi_{\pi^{*}}(S_{t+1})|\tilde{x}_{t}^{k}]\approx\frac{1}{H}\sum_{\omega\in\Omega}\alpha^{\omega}\Big(p_{u_{t}}^{\omega}(V_{k}^{+})+\sum_{V_{l}\neq\emptyset,l\neq k}p_{u_{t}}^{\omega}(V_{l})\Big).
\end{equation}

Policy S-PbP follows by using \eqref{eq:ka} to approximate the expected potential of states in step (iv) of the procedure described in Section \ref{sec:pbp}. \RefB{1}{Policy S-PbP requires solving model \eqref{model:mka} $H$ times for the computation of the numerators of \eqref{eq:alpha}. Additionally, for each non-idle vehicle that can serve $r_{t}$, model \eqref{model:ka} is solved $2H$ times: $H$ times for the `accept' and $H$ times for the `reject' decision.} Since variables are continuous, the single-knapsack model \eqref{model:ka} can be efficiently solved in $\mathcal{O}(|\omega|\log|\omega|)$ time without invoking a linear solver.

Analogously to PbP$_{_\textsf{R}}(H)$, we denote by S-PbP$_{_\textsf{R}}(H)$ the S-PbP when routing policy $\rho_{_\mathsf{R}}$ is adopted and $H$ sample paths are used to approximate potentials.

\subsection{Benchmark Policies}\label{sec:bench}
We discuss a few benchmark policies in this section.

\subsubsection{Greedy Policy}
Given a routing policy $\rho$ and a state $S_{t}=(V_{1},\ldots,V_{K},r_{t})$, $r_{t}=(u_{t},i_{t},d_{t})$, the greedy policy is specified as follows:

\begin{enumerate}[label=(\roman*)]
\item \textit{Same as step (i) from Section \ref{sec:pbp}.}
\item Otherwise, request $r_{t}$ is accepted if and only if there is a non-idle vehicle $k$, $V_{k}=(\tau_{k},\theta_{k})$, such that $b(\tau_{k},\rho(\theta_{k},r_{t}))\geq0$.
\item If $r_{t}$ is accepted, it is assigned to the non-idle vehicle $k$, $V_{k}=(\tau_{k},\theta_{k})$, such that $b(\tau_{k},\theta_{k})-b(\tau_{k},\rho(\theta_{k},r_{t}))$ is minimum.
\end{enumerate}

Again, different variants of the greedy policy arise, depending on the adopted routing policy. We denote by GP$_{_\textsf{CI}}$ and GP$_{_\textsf{R}}$ the greedy policies obtained by considering routing policies $\rho_{_\mathsf{CI}}$ and $\rho_{_\mathsf{R}}$, respectively.

\subsubsection{PFA Policy}\label{SubSecPFA}
PFA is a class of ADP methods that rely on analytical functions and is suitable for the cases where the structure of the optimal policy is known \citep{powell2012approximate}. We design a benchmark PFA policy in which the expected reward-to-go $\mathbb{E}[\Phi_{\pi^*}(S_{t+1})|x_t]$ is estimated by an analytical function $\hat{\Phi}_{_\textsf{PFA}}(S_{t+1}|x_t)$ parameterized by a value $\gamma\in\mathbb{R}_{>0}$. \RefA{4}{The precise mathematical definition of $\hat{\Phi}_{_\textsf{PFA}}(S_{t+1}|x_t)$ and the offline procedure to tune parameter $\gamma$ are detailed in \ref{sec:PFA}.} Given a routing policy $\rho$ and a state $S_{t}=(V_{1},\ldots,V_{K},r_{t})$, $r_{t}=(u_{t},i_{t},d_{t})$, the PFA policy is specified as follows:

\begin{enumerate}[label=(\roman*)]\setcounter{enumi}{3}
\item [(i)-(iii)] \textit{Same as steps (i)-(iii) from Section \ref{sec:pbp}.}
\item A decision is selected from $\tilde{\mathcal{X}}_{t}$ by evaluating decision rules \eqref{eq:dr1} and \eqref{eq:dr2} based on approximations $\hat{\Phi}_{_\textsf{PFA}}(S_{t+1}|\tilde{x}_{t})$ of $\mathbb{E}[\Phi_{\pi^{*}}(S_{t+1})|\tilde{x}_{t}]$, $\tilde{x}_{t}\in\tilde{\mathcal{X}}_{t}$, computed by \eqref{eq:pfa}.
\end{enumerate}

We denote by PFA$_{_\textsf{CI}}$ and PFA$_{_\textsf{R}}$ the PFA policies with routing policies $\rho_{_\mathsf{CI}}$ and $\rho_{_\textsf{R}}$ adopted in the online phase, respectively. \RefB{6}{As the offline procedure to tune parameter $\gamma$ requires simulating the policy under a large number of sample paths,} to guarantee computational efficiency the routing policy $\rho_{_\mathsf{CI}}$ is employed in the offline simulations.

\subsubsection{Rollout Policies}
A rollout algorithm is a general procedure to improve the performance of a base policy. It works by simulating candidate decisions, online and under the base policy, and then prescribing the decision with the maximum simulated reward-to-go. We refer to \cite{bertsekas1996neuro} for more details and a formal introduction to the method. As shown in Table \ref{TabLiterature}, rollout algorithms are frequently employed for solving DVRPs, and for this reason we also define rollout-based benchmark policies.

Given a base policy, a routing policy $\rho$ and a state $S_{t}=(V_{1},\ldots,V_{K},r_{t})$, $r_{t}=(u_{t},i_{t},d_{t})$, the rollout policy is specified as follows:

\begin{enumerate}[label=(\roman*)]\setcounter{enumi}{3}
\item[(i)-(iii)] \textit{Same as steps (i)-(iii) from Section \ref{sec:pbp}.}
\item For each candidate decision $\tilde{x}_{t}\in\tilde{\mathcal{X}}_{t}$, the system is simulated over $H$ sample paths from the current instant $u_{t}$ until the last instant $U$, given that decision $\tilde{x}_{t}$ is taken at state $S_{t}$ and the base policy is in effect. As a result of the simulations, we obtain $\hat{r}(\tilde{x}_{t})$, $\tilde{x}_{t}\in\tilde{\mathcal{X}}_{t}$, the approximated reward-to-go associated with decision $\tilde{x}_{t}$.
\item The policy prescribes decision $\tilde{x}_{t}\in\tilde{\mathcal{X}}_{t}$ such that $\hat{r}(\tilde{x}_{t})$ is maximum.
\end{enumerate}

Because of the reduced computational time available in a real-time context and the large number of online simulations required in step (iv) of the rollout procedure, only very efficient scheduling policies are suitable as base policies. Among the policies previously presented, only GP$_{_\textsf{CI}}$ and PFA$_{_\textsf{CI}}$ are appropriate as base policies, since all other policies require solving multiple integer programs or multiple linear programs before a decision is taken. \RefB{7}{In fact, as we discuss in Section \ref{SubSecExpOnline}, even under the computationally efficient base policies GP$_{_\textsf{CI}}$ and PFA$_{_\textsf{CI}}$ the decision times of rollout policies in large instances are severely compromised due to the costly online simulations.}

Finally, different routing policies can be adopted when defining the set of candidate decisions in step (iii) of the procedure. So, we denote by R$_{_\textsf{CI}}$-GP$_{_\textsf{CI}}(H)$ and R$_{_\textsf{CI}}$-PFA$_{_\textsf{CI}}(H)$ the rollout policies with base policies GP$_{_\textsf{CI}}$ and PFA$_{_\textsf{CI}}$, respectively, where routing policy $\rho_{_\mathsf{CI}}$ is adopted to generate candidate decisions and $H$ samples paths are simulated in step (iv) to approximate rewards-to-go. Likewise, we define rollout policies R$_{_\textsf{R}}$-GP$_{_\textsf{CI}}(H)$ and  R$_{_\textsf{R}}$-PFA$_{_\textsf{CI}}(H)$ when routing policy $\rho_{_\mathsf{R}}$ is used to generate candidate decisions.

\section{Offline Planning Algorithms}\label{SecOffline}
\RefA{6}{We now focus on the decision at the initial state $S_{0}$. In Section \ref{SecOfflineARP}, we briefly discuss previous methods for route planning in DVRPs, and formalize the problem. Section \ref{sec:budplanner} describes a benchmark method based on travel time minimization, and Section \ref{sec:potplanner} presents our potential-based approach for route planning.}

\subsection{Anticipative Route Planning}\label{SecOfflineARP}
\RefA{6}{A good route plan should cover static requests and leave the system in a favorable state for fitting dynamic requests that arrive during the service period. In this context, previous work proposed designing route plans that include static and sampled dynamic requests \citep{bent2004scenario}, waiting strategies in anticipation of the location of dynamic requests \citep{Branke_2005,thomas2007waiting}, combinations between waiting and request sampling \citep{ichoua2006exploiting}, and managing slack times to facilitate acceptance of future requests when time windows are enforced \citep{Mitrovi_Mini__2004}. For a detailed discussion, we refer to \cite{Ichoua12007}. Our potential-based planner is closest to the multiple scenario approach from \cite{bent2004scenario}, as we rely on sample paths from the spatiotemporal request distribution to evaluate the potential of candidate planned routes.}

Formally, the offline planning algorithm defines an initial set of routes $\mathbf{y}=\{\theta_{1},\ldots,\theta_{k'}\}$, $k'\leq K$, such that all static requests are covered. For the purpose of computing an offline plan, it is useful to reformulate \eqref{model:som} in a way more suitable for static VRP optimization:
\begin{align}
&\max_{\{\theta_{1},\ldots,\theta_{k'}\}\subset\Theta,k'\leq K}&&g(S_{1})\equiv\mathbb{E}\left[\Phi_{\pi^{*}}(S_{1})\right],\nonumber\\
&\text{s.t.}&&\sum_{k\in\{1,\ldots,k'\}}\big[(u,i,d)\in\Delta(\theta_{k})\big]=1,&(u,i,d)\in\mathcal{S},\label{eq:spc}\\
&&&V_{k}=(0,\theta_{k}),&k\in\{1,\ldots,k'\},\nonumber\\
&&&V_{k}=\emptyset,&k\in\{k'+1,\ldots,K\},\label{model:spm}\tag{DVRP-SPM}
\end{align}
where $S_{1}=(V_{1},\ldots,V_{K},r_{1})$, $r_{1}=(u_{1},i_{1},d_{1})$, and $\Theta=\{\theta:b(0,\theta)\geq0\}$ is the set of all feasible routes.

\RefB{1}{The feasible solution space of \eqref{model:spm} is defined by the partitioning constraints \eqref{eq:spc}. The states of the vehicles to which a planned route is assigned are initialized accordingly, and it is assumed that these vehicles depart from the depot at instant $0$ (i.e., waiting time at the depot is not allowed).} The cost function is non-linear and also depends on the optimal (but unknown) policy $\pi^{*}$. A reasonable step towards a tractable offline planning algorithm is to replace the cost function $g$ by an additively separable (on $V_{1},\ldots,V_{K}$), mostly order-preserving mapping $\tilde{g}$, and solve \eqref{model:spm} with a standard (static) VRP algorithm or heuristic. In the remainder of this section, we discuss two possible ways of implementing this approach.

\subsection{Budget-based (Myopic) Planner}\label{sec:budplanner}
The budget of a route is a suitable proxy for assessing the capacity of the route to accommodate future requests. Therefore, a sensible mapping is given by:
\begin{equation}\label{eq:myo}
\tilde{g}_{\mathsf{myo}}(S_{1})=\sum_{V_{k}\neq\emptyset}b(0,\theta_{k})+\sum_{V_{k}=\emptyset}U.
\end{equation}

We call $\tilde{g}_{\mathsf{myo}}$ the \emph{myopic} mapping, since it does not anticipate the realization of future dynamic requests. Under cost function \eqref{eq:myo}, model \eqref{model:spm} reduces to a duration-constrained (or distance-constrained) VRP \citep[DCVRP;][]{toth2002vehicle} defined on an auxiliary, complete digraph $\mathcal{G}'=(\mathcal{V}', \mathcal{A}')$, which is created in the following way: first, the set of nodes is defined as $\mathcal{V}'=\{0\}\cup\{i:(u,i,d)\in\mathcal{S}\}$; second, with each arc $(i,j)\in\mathcal{A}'$, where $(u_{i},i,d_{i}),(u_{j},j,d_{j})\in\mathcal{S}$, a cost $t(i,j)+d_{j}$ is associated. Then, a myopic plan is obtained by solving a DCVRP on $\mathcal{G}'$, with node $0$ as the depot, up to $K$ vehicles available, and a duration limit of $U$ on individual routes.

Depending on the number of static requests $|\mathcal{S}|$, finding the optimal solution to the DCVRP may be computationally prohibitive. Moreover, a near-optimal solution is enough to allow us to assess the quality of budget-based first stage DVRPSR decisions, compared to the potential-based plans described in the next section. Therefore, we opt for a heuristic method instead of an exact algorithm to solve the DCVRP. The two-step heuristic works as follows:

\begin{enumerate}[label=(\roman*)]
\item A large set of feasible routes $\Theta'\subset\Theta$ is generated by applying column generation \citep{barnhart1998branch} on the linear relaxation of model \eqref{model:spm} under cost function \eqref{eq:myo}. The details of this procedure, which employs standard techniques commonly applied in (static) VRP algorithms, are described in \ref{sec:colgen}.
\item A myopic plan is determined by finding the set of routes $\{\theta_{1},\ldots,\theta_{k'}\}\subset\Theta'$, $k'\leq K$, \RefB{1}{that serve all static requests with minimum total duration. To this end, we employ an off-the-shelf integer solver to solve exactly model \eqref{model:spm} under cost function \eqref{eq:myo} with the restricted set of variables~$\Theta'$.}\label{st:bnb}
\end{enumerate}

\subsection{Potential-based Planner}\label{sec:potplanner}
\RefB{8}{Instead of maximizing budget, the potential-based planner maximizes the potential of the state resulting from the first stage decision.} Given sample paths $\Omega=\{\omega_{1},\ldots,\omega_{H}\}$ of the request arrival process $\{\Lambda(u)\}_{u\in[0,U]}$, we define the potential-based mapping $\tilde{g}_{\mathsf{pb}}$:
\begin{equation}\label{eq:pb}
\tilde{g}_{\mathsf{pb}}(S_{1})=\frac{1}{H}\sum_{\omega\in\Omega}\sum_{V_{k}\neq\emptyset}p_{0}^{\omega}(V_{k}),
\end{equation}
\RefB{8}{where $V_{k}=(0,\theta_{k})$, $k\in\{1,\ldots,k'\}$, are the initial states of the vehicles to which a route is assigned.}

Mapping $\tilde{g}_{\mathsf{pb}}$ takes into account the distribution of future requests and is also additively separable on $V_{1},\ldots,V_{K}$. An algorithm for generating potential-based plans follows:

\begin{enumerate}[label=(\roman*)]
\item \textit{Same as step (i) from Section \ref{sec:budplanner}.}
\item \RefB{8}{For each route $\theta\in\Theta'$, compute the cost coefficient $\sum_{\omega\in\Omega}p_{0}^{\omega}((0,\theta))/H$, which correlates with the number of dynamic requests that can be accommodated along $\theta$.}
\item A potential-based plan is determined by finding the set of routes $\{\theta_{1},\ldots,\theta_{k'}\}\subset\Theta'$, $k'\leq K$, \RefB{8}{that serve all static requests with maximum (estimated) potential. As in step \ref{st:bnb} of the myopic planner, we invoke an integer solver to solve exactly model \eqref{model:spm} under cost function \eqref{eq:pb} with the restricted set of variables~$\Theta'$.}
\end{enumerate}

Compared to the myopic planner, the potential-based planner favors routes with a total budget that can be more efficiently used to serve dynamic requests. \RefA{6}{The potential-based mapping promotes flexibility to serve dynamic requests as planned routes with more slack time are preferred over low-budget routes.} In addition, whereas the DCVRP-based myopic plan usually consists of few routes, each serving several static requests, the potential-based first stage decision exhibits a more balanced distribution of static requests to vehicles and promotes a better coverage of the service area. \RefA{6}{By considering a large number of sample paths when evaluating \eqref{eq:pb}, we ensure that the initial plan is robust, in the sense that it performs well, on average, for a variety of possible dynamic request trajectories. Also note that, similar to the potential-based policies, the potential-based planner requires no parameter other than the number of sample paths.}

\section{Computational Study} \label{SecComputation}
In this section, we present an extensive computational study on large-scale test instances based on a real street network. Section \ref{SubSecInstances} introduces the detailed configurations of the test instances. Section \ref{SubSecExpOffline} compares the myopic and potential-based offline planning algorithms, and Section \ref{SubSecExpOnline} compares the online scheduling policies. All algorithms and simulations are programmed in C++ and executed on a single core of an Intel\textsuperscript{\textregistered} Xeon\textsuperscript{\textregistered} Gold 6130 (2.1GHz) processor with 8GB of available RAM. IBM\textsuperscript{\textregistered} CPLEX\textsuperscript{\textregistered} version 12.10 is employed for solving linear and integer programs. \RefB{3}{When adopting the reoptimization policy $\rho_{_\mathsf{R}}$, the resulting TSPs are solved with a standard branch-and-cut procedure.} Finally, all software and test instances developed in this project, with which all results presented in this section can be replicated, are available open-source at \href{https://github.com/amflorio/dvrp-stochastic-requests}{https://github.com/amflorio/dvrp-stochastic-requests}.

\subsection{Test Instances}\label{SubSecInstances}
\begin{figure}[!htb]
\centering
\caption{Graph Representing the Street Network of Vienna, Austria.\label{FigMap}}
\includegraphics[trim={20 20 20 20},clip,width=0.75\textwidth]{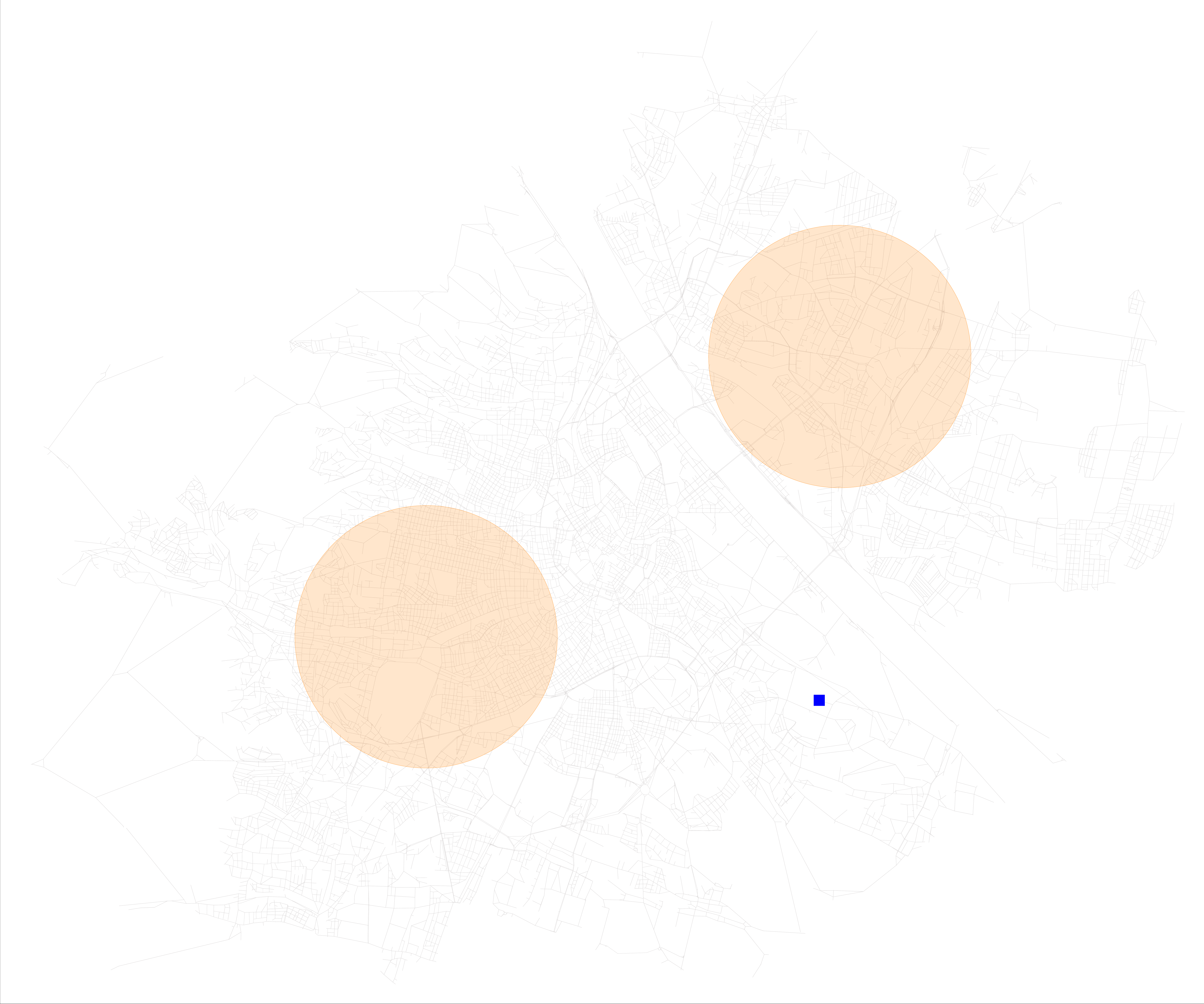}\\
{\footnotesize\textit{Notes.} The square represents the depot, and the two shaded areas represent request clusters.} 
\end{figure}

Using geographical data from the OpenStreetMap (OSM) project, we construct our test instances based on the street network of Vienna, Austria. The street network graph is generated from OSM data in the following way: first, we fetch all available data (nodes and ways in OSM terminology) for the city of Vienna; next, we remove the ways that correspond to walkways or service roads, and remove the nodes that are not intersections of two or more streets. The final graph consists of 16,080 nodes (intersections) and 36,424 arcs (street segments) and is illustrated in Figure \ref{FigMap}. The depot, indicated by the square, is located in the southeast of the city, and two request clusters, depicted as circular shaded areas, are defined in the northeast and southwest areas.

The fleet size $K$ varies from 2 to 20, and each vehicle travels at a constant speed $\bar{s}$=20 km/h. The duration of the service period is set to $U=600$ minutes. The ratio between the expected number of dynamic requests and the expected total number of requests (static and dynamic) is given by the degree of dynamism $\eta=\mathbb{E}[T/(|\mathcal{S}|+T)]$. We assume $\sum_{i\in\mathcal{V}\setminus\{0\}}\lambda_{i}(u)=\Lambda$, $u\in[0,U]$, where $\Lambda\in\{0.2,0.4,0.8,1.5\}$ is the (constant) overall request rate per minute. Note that our smallest instances with $\Lambda=0.2$ and $\eta=0.75$, such that $\mathbb{E}[T]+|\mathcal{S}|=160$, are larger than most synthetic instances considered in the previous works mentioned in Table \ref{TabLiterature}.

Given $\Lambda$ and $\eta$, the expected numbers of static and dynamic requests are fixed. Static requests are uniformly distributed on \RefA{7}{the network nodes. To assess the performance of the proposed methods under distinct dynamic request patterns, we generate dynamic requests }according to three spatiotemporal distributions: uniform and time-independent (UTI), clustered and time-independent (CTI), and clustered and time-dependent (CTD). In the UTI distribution, dynamic requests are also uniformly distributed on the network nodes ($\lambda_{i}(u)=\lambda_{j}(u)$, for all $i,j\in\mathcal{V}\setminus\{0\}$ and $u\in[0,U]$), while in the CTI and CTD distributions, half of the dynamic requests are expected to originate from a node within the two circular areas with a radius of 3 km (see Figure \ref{FigMap}). The request rate at each node does not change over time in the UTI and CTI distributions ($\lambda_{i}(u)=\lambda_{i}(u')$, for all $i\in\mathcal{V}\setminus\{0\}$ and $u,u'\in[0,U]$). In the CTD distribution, the request rates in the northeast and southwest clusters decrease and increase over time, respectively: 80\% (20\%) of the clustered requests are expected to appear in the northeast (southwest) cluster at the beginning of the service period ($u=0$), and this percentage decreases (increases) linearly to 20\% (80\%) at the end of the service period ($u=U$). \RefA{7}{The UTI distribution models more closely a practical large-scale instance, because the high density of nodes in the inner city induces a natural demand cluster.} \RefC{2.8}{Distributions CTI and CTD are somewhat arbitrary but useful nevertheless, as they allow testing the policies on request trajectories largely distinct from those given by the UTI distribution.} Regardless of which request distribution is adopted, each static or dynamic request's service duration is normally distributed with a mean of 10 minutes and a standard deviation of 150 seconds.

Each test instance is characterized by a certain combination of request rate $\Lambda$, degree of dynamism $\eta$, fleet size $K$, and request distribution. For each instance, we independently generate one set of static requests and multiple sets of dynamic requests. Combining the former with one of the latter yields a \textit{request scenario}, which corresponds to a specific problem and is solved individually. Table \ref{TabInstances} summarizes all test instances, the number of request scenarios per instance, as well as the offline and online algorithms employed to solve each scenario. We parameterize both the potential-based planner and policies with $H=50$, as our preliminary experiments indicated that additional sample paths do not lead to improvements in solution quality. Hence, for simplicity, S-PbP$_{_\textsf{R}}(50)$ and PbP$_{_\textsf{R}}(50)$ are from now on referred to as S-PbP and PbP, respectively. When solving a request scenario under a non-deterministic online policy (all policies except GP$_{_\textsf{CI}}$ and GP$_{_\textsf{R}}$), we repeat the solution procedure 15 times to evaluate the policy's performance in expectation. \RefB{10}{All sample paths required by the non-deterministic policies are independently generated.} Considering all instances, request scenarios, and algorithms tested, we simulate a total of 27,964 executions of the DVRPSR. As the total computational time required in this simulation study was quite high (approximately, 248 CPU days), we also make available in the online repository the complete raw simulation results.

\begin{table}[!htb]
    \renewcommand{\baselinestretch}{1.2}
    \caption{Test Instances and Offline and Online Algorithms Simulated\label{TabInstances}}
    \addtolength\tabcolsep{-0.1em}
    \begin{adjustbox}{center}
        \begin{threeparttable}
            \centering
            \footnotesize
                        \begin{tabular}{ccccccccr}
            \hline
                \multicolumn{5}{c}{Instance parameters} && \multicolumn{2}{c}{Algorithms} & \multirow{2}{*}{Simulations} \\\cmidrule{1-5}\cmidrule(r){7-8}
                $\Lambda$ & $\eta$ & $K$ & Dist. & Scen. && Offline & Online &  \\
            \hline\vspace{0.75em}
                0.2 & 0.75 & 2, 3 & \makecell{UTI\\CTI\\CTD} & 5 && \makecell{MY\\PB} & \makecell{GP$_{_\textsf{CI}}$, GP$_{_\textsf{R}}$, R$_{_\textsf{CI}}$-GP$_{_\textsf{CI}}$(10, 25, 50, 100),\\R$_{_\textsf{R}}$-GP$_{_\textsf{CI}}$(10, 25, 50, 100), PFA$_{_\textsf{CI}}$, PFA$_{_\textsf{R}}$,\\R$_{_\textsf{R}}$-PFA$_{_\textsf{CI}}$(25, 50, 100)$^*$, S-PbP, PbP} & 12,270\\\vspace{0.75em}
                0.4 & 0.85 & 3, 5 & \makecell{UTI\\CTI\\CTD} & 5 && \makecell{MY\\PB} & \makecell{GP$_{_\textsf{CI}}$, GP$_{_\textsf{R}}$, R$_{_\textsf{CI}}$-GP$_{_\textsf{CI}}$(10, 25, 50, 100),\\R$_{_\textsf{R}}$-GP$_{_\textsf{CI}}$(10, 25, 50, 100), PFA$_{_\textsf{CI}}$, PFA$_{_\textsf{R}}$,\\R$_{_\textsf{R}}$-PFA$_{_\textsf{CI}}$(25, 50, 100)$^*$, S-PbP, PbP} & 12,270\\\vspace{0.75em}
                0.8 & 0.90 & 6, 12 & \makecell{UTI\\CTI\\CTD} & 5 && PB & \makecell{GP$_{_\textsf{CI}}$, GP$_{_\textsf{R}}$, PFA$_{_\textsf{CI}}$, PFA$_{_\textsf{R}}$,\\R$_{_\textsf{CI}}$-GP$_{_\textsf{CI}}$(10), R$_{_\textsf{R}}$-GP$_{_\textsf{CI}}$(10),\\R$_{_\textsf{R}}$-PFA$_{_\textsf{CI}}$(10), S-PbP, PbP} & 3,210\\
                1.5 & 0.95 & 10, 20 & UTI & 1** && PB & \makecell{GP$_{_\textsf{CI}}$, GP$_{_\textsf{R}}$, PFA$_{_\textsf{CI}}$, PFA$_{_\textsf{R}}$,\\ R$_{_\textsf{CI}}$-GP$_{_\textsf{CI}}$(10), R$_{_\textsf{R}}$-GP$_{_\textsf{CI}}$(10),\\R$_{_\textsf{R}}$-PFA$_{_\textsf{CI}}$(10), S-PbP, PbP} & 214\\
            \hline
            \end{tabular}
            \footnotesize
            \textit{Notes.} Dist., distribution; Scen., scenario; MY, myopic planner; PB, potential-based planner;\\
            * The online algorithms marked by * are tested with the potential-based planner only;\\
            {** Only one scenario is tested for each large instance with $\Lambda=1.5$ due to the long simulation time.}
        \end{threeparttable}
    \end{adjustbox}\vspace{0.5em}
    \end{table}

Throughout each simulation, the implemented framework records the arrival of requests, the decisions taken, and the states' evolution, and it also outputs each state in a visual format. The visualization of a simulation enables a better interpretation of the decisions taken by the scheduling policy. As an example, at \href{https://youtu.be/D57xNfU73as}{https://youtu.be/D57xNfU73as} we can visualize a simulation of policies R$_{_\textsf{R}}$-GP$_{_\textsf{CI}}$(50) and PbP on an instance with $\Lambda=0.4$, $K=5$ and UTI distributed requests. The main indicators used to evaluate the performance of the algorithms are the \emph{acceptance rate}, defined as the number of dynamic requests accepted divided by the total number of dynamic requests generated, and the \emph{maximum decision time}, which is the maximum computation time required for taking a single decision. The maximum decision time is an important metric to assess a scheduling policy's suitability for real-time applications.

\subsection{Evaluation of Offline Planning Algorithms}\label{SubSecExpOffline}
To evaluate the performance of the myopic and potential-based planners, we employ both algorithms to compute initial route plans for all instances with $\Lambda\in\{0.2,0.4\}$. Table \ref{TabOfflineAR} compares the average acceptance rates across eight online policies, where columns ``Impr.'' indicate the percentage improvement achieved by the potential-based plan in relation to the myopic plan when a specific online algorithm is used. \RefA{8}{A more detailed comparison, which shows the results of different request distributions, is provided in Table \ref{TabOfflineARDetail} in \ref{sec:detailed results}.} \RefB{1}{Overall, when combined with the best-performing policies S-PbP and PbP, the potential-based plans improve acceptance rates by $3.6\%$ and $5.0\%$, respectively. Although myopic and potential-based plans perform equally well in the smallest instance with $\Lambda=0.2$ and $K=2$, in all other instances the potential-based plans significantly improve the average acceptance rates, regardless of the request distribution and the adopted online policy.}

\begin{table}[!htb]
\begin{center}
    \renewcommand{\baselinestretch}{1.2}
    \caption{Comparison of Offline Planning Algorithms: Acceptance Rate (\%) and Relative Improvement (\%)\label{TabOfflineAR}}
\footnotesize
\begin{threeparttable}
\addtolength\tabcolsep{-0.05em}
\begin{tabular}{cccrrrrrrrrrrrrrrrr}
	\hline
	\multicolumn{2}{c}{Instance} & & \multicolumn{3}{c}{GP$_{_\textsf{CI}}$} &  & \multicolumn{3}{c}{GP$_{_\textsf{R}}$} &  & \multicolumn{3}{c}{R$_{_\textsf{CI}}$-GP$_{_\textsf{CI}}$(50)} &  & \multicolumn{3}{c}{R$_{_\textsf{R}}$-GP$_{_\textsf{CI}}$(50)} \\\cmidrule(lr){1-2}\cmidrule(lr){4-6}\cmidrule(lr){8-10}\cmidrule(lr){12-14}\cmidrule(lr){16-18}
	$\Lambda$ & $K$ &  & \multicolumn{1}{c}{MY} & \multicolumn{1}{c}{PB} & \multicolumn{1}{c}{Impr.} &  & \multicolumn{1}{c}{MY} & \multicolumn{1}{c}{PB} & \multicolumn{1}{c}{Impr.} &  & \multicolumn{1}{c}{MY} & \multicolumn{1}{c}{PB} & \multicolumn{1}{c}{Impr.} &  & \multicolumn{1}{c}{MY} & \multicolumn{1}{c}{PB} & \multicolumn{1}{c}{Impr.} \\\cmidrule(lr){1-18}
	0.2 & 2                 && 19.5 & 19.5 & 0.0    && 19.8 & 19.8 & 0.0    && 24.2 & 24.3 & +0.4*  && 24.4 & 24.3 & $-$0.1* \\\vspace{0.5em}
        & 3                 && 39.1 & 43.6 & +11.4  && 39.6 & 45.4 & +14.8  && 49.3 & 50.6 & +2.5   && 49.6 & 52.0 & +4.8 \\
    0.4 & 3                 && 20.0 & 23.0 & +15.0  && 20.3 & 23.6 & +16.1  && 26.0 & 28.0 & +7.9   && 26.1 & 28.6 & +9.5 \\\vspace{0.5em}
        & 5                 && 42.6 & 49.2 & +15.6  && 43.1 & 50.1 & +16.4  && 51.4 & 55.2 & +7.2   && 51.6 & 55.7 & +7.9 \\
\multicolumn{2}{c}{Overall} && 30.3 & 33.8 & +11.7  && 30.7 & 34.8 & +13.2  && 37.7 & 39.5 & +4.7   && 37.9 & 40.1 & +5.9 \\
	\cmidrule(lr){1-18}
	\multicolumn{2}{c}{Instance} & & \multicolumn{3}{c}{PFA$_{_\textsf{CI}}$} &  & \multicolumn{3}{c}{PFA$_{_\textsf{R}}$} &  & \multicolumn{3}{c}{S-PbP} &  & \multicolumn{3}{c}{PbP} \\\cmidrule(lr){1-2}\cmidrule(lr){4-6}\cmidrule(lr){8-10}\cmidrule(lr){12-14}\cmidrule(lr){16-18}
	$\Lambda$ & $K$ &  & \multicolumn{1}{c}{MY} & \multicolumn{1}{c}{PB} & \multicolumn{1}{c}{Impr.} &  & \multicolumn{1}{c}{MY} & \multicolumn{1}{c}{PB} & \multicolumn{1}{c}{Impr.} &  & \multicolumn{1}{c}{MY} & \multicolumn{1}{c}{PB} & \multicolumn{1}{c}{Impr.} &  & \multicolumn{1}{c}{MY} & \multicolumn{1}{c}{PB} & \multicolumn{1}{c}{Impr.} \\\cmidrule(lr){1-18}
	0.2 & 2                 && 26.6 & 26.6 & $-$0.2*&& 26.8 & 26.8 & 0.0    && 27.7 & 27.7 & 0.0    && 27.7 & 27.7 & +0.1*   \\\vspace{0.5em}
        & 3                 && 41.6 & 48.5 & +16.6  && 42.8 & 50.2 & +17.5  && 52.0 & 52.7 & +1.3   && 52.5 & 54.3 & +3.5    \\
    0.4 & 3                 && 24.0 & 29.8 & +24.4  && 24.2 & 30.0 & +24.0  && 28.7 & 31.4 & +9.2   && 28.9 & 31.8 & +10.0   \\\vspace{0.5em}
        & 5                 && 45.9 & 54.3 & +18.5  && 46.2 & 54.4 & +17.7  && 53.3 & 55.8 & +4.5   && 54.5 & 58.0 & +6.4    \\
\multicolumn{2}{c}{Overall} && 34.5 & 39.8 & +15.3  && 35.0 & 40.4 & +15.3  && 40.5 & 41.9 & +3.6   && 40.9 & 43.0 & +5.0    \\
	\hline
\end{tabular}
\footnotesize\textit{Notes.} Impr., relative improvement of PB over MY: (PB$-$MY)/MY$\times\%$.\\
{* indicates the improvements that are NOT statistically significant  (i.e., $p>0.05$ in paired two-sample \textit{t}-tests).}
\end{threeparttable}
\end{center}
\end{table}

\begin{figure}[!htb]
\centering
\caption{Comparison of Initial Route Plans.\label{FigInitialPlans}}
\includegraphics[width=0.9\textwidth]{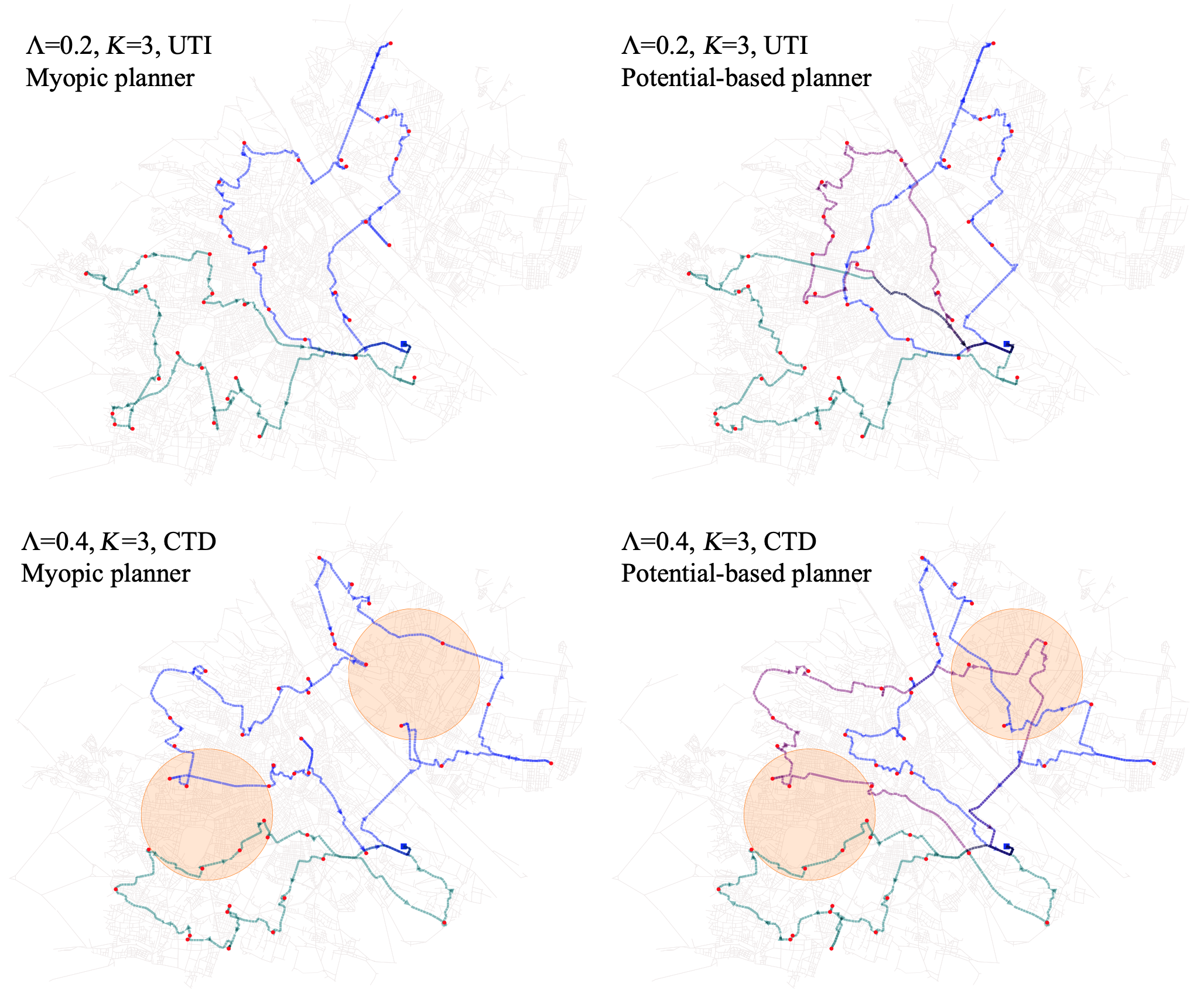}\\
{\footnotesize\textit{Notes.} Each red dot represents a static request. Each colored line represents the initial route of a vehicle.}
\end{figure}

Figure \ref{FigInitialPlans} compares the initial plans generated by the two offline planning algorithms for two instances to develop intuition behind these improvements. In the upper instance, we see that the myopic planner defines two initial routes and leaves one vehicle idle at the depot (which is activated upon arrival of the first dynamic request). This myopic plan maximizes the total budget but ignores the spatiotemporal distribution of dynamic requests. In comparison, the potential-based planner defines one additional route and, by doing so, promotes a better coverage of the dense area in the city center, which constitutes a natural cluster as requests are uniformly distributed over intersections. Similarly, in the lower instance, we observe that the potential-based planner uses the third vehicle to improve the northeast cluster's coverage, where the request rate is high, especially during the early stage of the service period.

Potential-based plans' success is also attributed to a more balanced set of initial routes, which cover similar distances and serve a similar number of static requests. \RefB{5}{This is achieved simply by evaluating routes under the knapsack-based potential approximation, without requiring additional parameters.} The capacity (in terms of budget) of those routes is consumed in a more balanced way along the service period, which gives more flexibility to the scheduling policy when deciding the best vehicle to serve a dynamic request. \RefB{5}{On the other hand, myopic plans are characterized by very long routes and idle vehicles.} While it is possible to enforce the myopic planner to use all vehicles, this alone does not lead to a balanced set of routes as the DCVRP model minimizes the total duration, so additional constraints (and parameters) would be required to achieve the same effect.

Table \ref{TabOfflineCT} presents the runtime requirements of both planners. The fleet size and the request distribution have little effect on the offline runtime, so the results correspond to the averages over all fleet sizes and request distributions for each request rate $\Lambda$. For each instance, the two planners generate the same route set $\Theta'$ because they have the same column generation procedure as their first step. The runtime of the potential-based planner increases with $\Lambda$, as it requires each route to be simulated $H=50$ times under sample paths that grow linearly with $\Lambda$.

\begin{table}[!htb]
    \renewcommand{\baselinestretch}{1.2}
    \caption{Average Runtime of Offline Planning Algorithms\label{TabOfflineCT}}
    \centering
    \footnotesize
    \begin{threeparttable}
    \begin{tabular}{ccccccccc}
    \hline
    \multirow{2}{*}{$\Lambda$} & \multirow{2}{*}{$\mathbb{E}[T]$} & \multirow{2}{*}{$|\mathcal{S}|$} & \multirow{2}{*}{$|\Theta'|$} & \multicolumn{2}{c}{Time (min)} \\
    \cmidrule{5-6}
     &  &  &  & MY & PB \\
    \hline
    0.2 & 120 & 40 & 5,986 & 22 & 26 \\
    0.4 & 240 & 42 & 4,656 & 26 & 31 \\
    0.8 & 480 & 53 & 6,279 & 55 & 66 \\
    1.5 & 900 & 47 & 5,748 & 37 & 57 \\
    \hline
    \end{tabular}
    \end{threeparttable}\\\vspace{0.5em}
    \footnotesize\textit{Notes.} MY, myopic planner; PB, potential-based planner.
    \end{table}

\subsection{Evaluation of Online Scheduling Policies}\label{SubSecExpOnline}
We first compare the online scheduling policies applied to the instances with $\Lambda\in\{0.2,0.4\}$. As Section \ref{SubSecExpOffline} demonstrates that potential-based plans outperform myopic plans, all results presented in this section are based on the former. \RefA{8}{Moreover, since the spatiotemporal distribution of dynamic requests has little impact on the performance of online policies, all results presented in this section are averaged over the three distributions. The interested reader is referred to Tables \ref{TabOnlineARSmallDetail} to \ref{TabOnlineTLargeDetail} in \ref{sec:detailed results} for detailed results.} Taking the results of policy GP$_{_\textsf{CI}}$ as a reference, Table \ref{TabOnlineARSmall} presents the relative improvements in acceptance rate obtained by the other policies (the acceptance rates of GP$_{_\textsf{CI}}$ can be found in Table \ref{TabOfflineAR}). In addition, the maximum decision times are presented in Table \ref{TabOnlineDTSmall}.

\begin{table}[!htb]
\renewcommand{\baselinestretch}{1.2}
    \caption{Comparison of Online Scheduling Policies ($\Lambda\in\{0.2,0.4\}$): Percentage Improvement in the Acceptance Rate Compared to GP$_{_\textsf{CI}}$\label{TabOnlineARSmall}}
    \footnotesize
\addtolength\tabcolsep{0.0em}
\begin{adjustbox}{center}
\begin{threeparttable}
\begin{tabular}{ccrrrrrrrrrrr}
	\hline
	\multicolumn{2}{c}{Instance} & \multicolumn{11}{c}{Online scheduling policy} \\\cmidrule(lr){1-2}\cmidrule(lr){3-13}
	\multirow{2}{*}{$\Lambda$} & \multirow{2}{*}{$K$} & \multicolumn{1}{c}{\multirow{2}{*}{GP$_{_\textsf{R}}$}} & \multicolumn{3}{c}{R$_{_\textsf{R}}$-GP$_{_\textsf{CI}}(H)$} & \multicolumn{1}{c}{\multirow{2}{*}{PFA$_{_\textsf{CI}}$}} & \multicolumn{1}{c}{\multirow{2}{*}{PFA$_{_\textsf{R}}$}} & \multicolumn{3}{c}{R$_{_\textsf{R}}$-PFA$_{_\textsf{CI}}(H)$} & \multicolumn{1}{c}{\multirow{2}{*}{S-PbP}} & \multicolumn{1}{c}{\multirow{2}{*}{PbP}} \\\cmidrule(lr){4-6}\cmidrule(lr){9-11}
	& & & \multicolumn{1}{c}{$H$=25} & \multicolumn{1}{c}{$H$=50} & \multicolumn{1}{c}{$H$=100} & & & \multicolumn{1}{c}{$H$=25} & \multicolumn{1}{c}{$H$=50} & \multicolumn{1}{c}{$H$=100} & & \\\cmidrule(lr){1-13}
	                0.2 & 2 & +1.9* & +23.5 & +24.9 & +25.7 & +36.4 & +37.7 & +39.6 & +40.8 & +40.9 & +42.4 & +42.4 \\\vspace{0.5em}
                        & 3 & +4.3  & +17.3 & +19.3 & +20.8 & +11.4 & +15.4 & +21.1 & +22.5 & +23.6 & +20.9 & +24.7 \\
                    0.4 & 3 & +2.8* & +21.7 & +24.5 & +26.1 & +29.8 & +30.7 & +34.7 & +36.6 & +37.4 & +36.8 & +38.5 \\\vspace{0.5em}
                        & 5 & +1.8* & +10.3 & +13.1 & +15.1 & +10.4 & +10.5 & +14.3 & +16.6 & +17.9 & +13.2 & +17.7 \\
\multicolumn{2}{c}{Overall} & +2.8  & +16.4 & +18.7 & +20.4 & +17.7 & +19.4 & +23.6 & +25.4 & +26.4 & +23.9 & +27.1 \\
	\hline
\end{tabular}
{* indicates the improvements that are NOT statistically significant  (i.e., $p>0.05$ in paired two-sample \textit{t}-tests).}
\end{threeparttable}
\end{adjustbox}
    \end{table}

\begin{table}[!htb]
\renewcommand{\baselinestretch}{1.2}
    \caption{Comparison of Online Scheduling Policies ($\Lambda\in\{0.2,0.4\}$): Maximum Decision Time (seconds)\label{TabOnlineDTSmall}}
    \footnotesize
\addtolength\tabcolsep{-0.1em}
\begin{adjustbox}{center}
\begin{tabular}{ccrrrrrrrrrrrr}
	\hline
	\multicolumn{2}{c}{Instance} & \multicolumn{12}{c}{Online scheduling policy} \\\cmidrule(lr){1-2}\cmidrule(lr){3-14}
	\multirow{2}{*}{$\Lambda$} & \multirow{2}{*}{$K$} & \multicolumn{1}{c}{\multirow{2}{*}{GP$_{_\textsf{CI}}$}} & \multicolumn{1}{c}{\multirow{2}{*}{GP$_{_\textsf{R}}$}} & \multicolumn{3}{c}{R$_{_\textsf{R}}$-GP$_{_\textsf{CI}}(H)$} & \multicolumn{1}{c}{\multirow{2}{*}{PFA$_{_\textsf{CI}}$}} & \multicolumn{1}{c}{\multirow{2}{*}{PFA$_{_\textsf{R}}$}} & \multicolumn{3}{c}{R$_{_\textsf{R}}$-PFA$_{_\textsf{CI}}(H)$} & \multicolumn{1}{c}{\multirow{2}{*}{S-PbP}} & \multicolumn{1}{c}{\multirow{2}{*}{PbP}} \\\cmidrule(lr){5-7}\cmidrule(lr){10-12}
	& & & & \multicolumn{1}{c}{$H$=25} & \multicolumn{1}{c}{$H$=50} & \multicolumn{1}{c}{$H$=100} & & & \multicolumn{1}{c}{$H$=25} & \multicolumn{1}{c}{$H$=50} & \multicolumn{1}{c}{$H$=100} & & \\\cmidrule(lr){1-14}
	                0.2 & 2 & $<$0.1 & 0.3 & 1.0 & 2.0 & 3.1    & $<$0.1 & 0.3 & 1.0 & 2.1 & 3.9 & 0.4 & 0.6    \\\vspace{0.5em}
                        & 3 & $<$0.1 & 0.2 & 1.9 & 3.7 & 6.3    & $<$0.1 & 0.2 & 1.9 & 3.5 & 7.2 & 0.4 & 0.6    \\
                    0.4 & 3 & $<$0.1 & 0.4 & 3.1 & 6.1 & 14.1   & $<$0.1 & 0.3 & 3.2 & 7.3 & 14.5 & 0.6 & 1.1   \\\vspace{0.5em}
                        & 5 & $<$0.1 & 0.3 & 7.4 & 14.6 & 33.2  & $<$0.1 & 0.3 & 8.7 & 17.1 & 34.5 & 0.8 & 2.4  \\
\multicolumn{2}{c}{Overall} & $<$0.1 & 0.4 & 7.4 & 14.6 & 33.2  & $<$0.1 & 0.3 & 8.7 & 17.1 & 34.5 & 0.8 & 2.4  \\
	\hline
\end{tabular}
\end{adjustbox}
    \end{table}

The best-performing policy is PbP, which displays an average acceptance rate higher than those of the PFA and rollout algorithms, with or without reoptimization, and for any number of sample paths. Moreover, PbP can take decisions much faster than rollout policies. For the same number of sample paths $H=50$, PbP promotes a sixfold decrease in decision times, on average. Policy S-PbP is even more computationally efficient, at the cost of somewhat lower acceptance rates compared to PbP (but still better than \RefB{1}{the other policies except R$_{_\textsf{R}}$-PFA$_{_\textsf{CI}}(H\geq50)$}, on average). In these instances, PFA-based policies require no more than 60 seconds of offline computation time for parameter tuning. We also notice that, compared to the CI routing policy $\rho_{_\textsf{CI}}$, the reoptimization routing policy $\rho_{_\textsf{R}}$ improves the acceptance rates slightly in both the greedy and PFA policies, at the cost of short additional computation time per decision. Similar differences caused by $\rho_{_\textsf{CI}}$ and $\rho_{_\textsf{R}}$ can be observed in the rollout policies, as illustrated in Figure \ref{FigEvaluationRCI} in \ref{sec:detailed results}.

\begin{figure}[!htb]
\centering
\caption{Comparison of Online Scheduling Policies ($\Lambda\in\{0.2,0.4\}$).\label{FigEvaluationRH}}
{\includegraphics[width=0.6\textwidth]{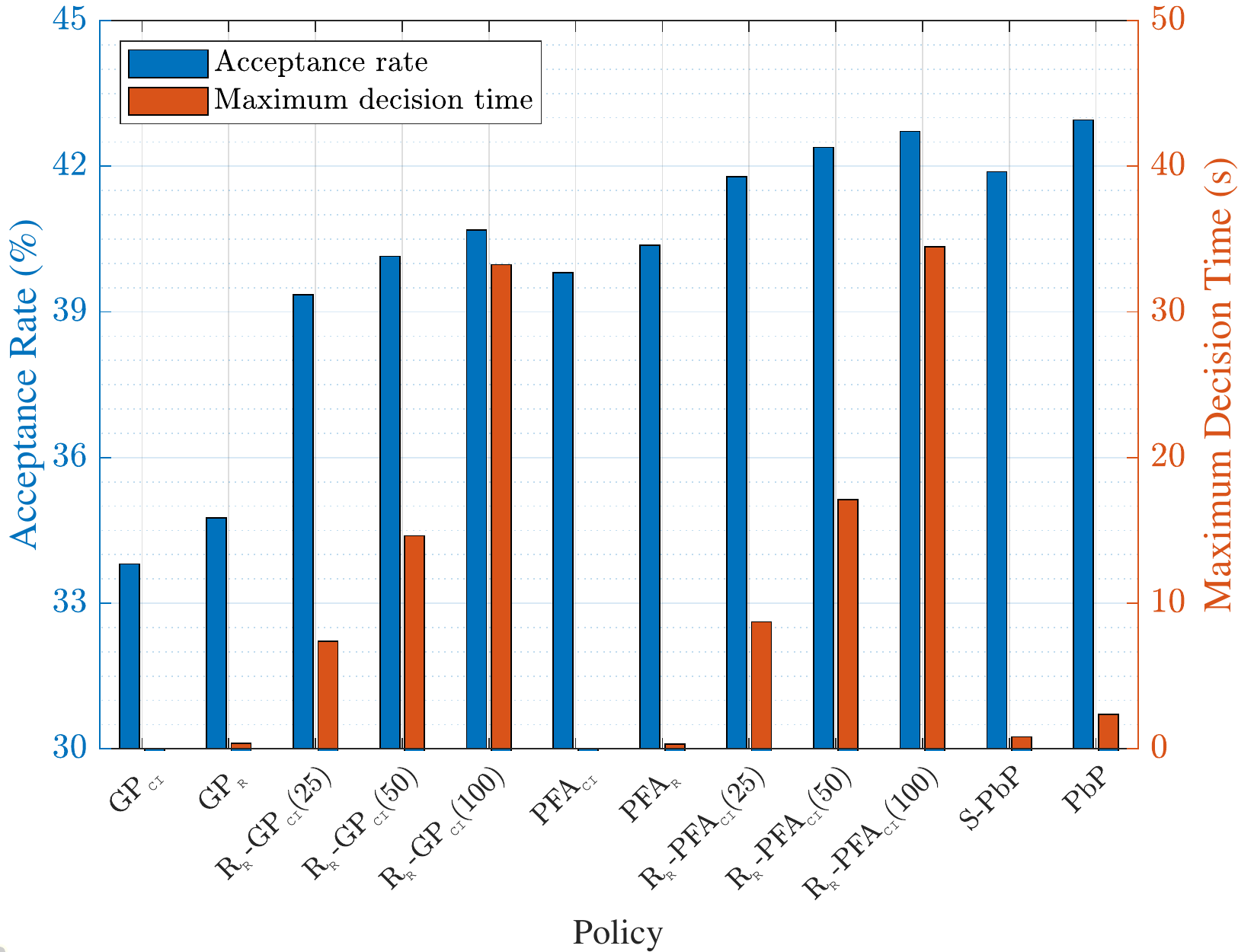}}
\end{figure}

Figure \ref{FigEvaluationRH} compares the overall performance and decision times of all scheduling algorithms (except the rollout algorithms with the CI routing policy $\rho_{_\textsf{CI}}$) graphically. As expected, the performance and decision times of rollout policies increase with the number of sample paths. As also expected, policies GP$_{_\textsf{CI}}$ and PFA$_{_\textsf{CI}}$ are the least computationally expensive, but the former has the worst solution quality. Although the PFA-based policies lead to drastically better results than the GP-based policies, their solution quality does not surpass that of policy S-PbP until the number of sample paths increases to $H=50$, at which point their maximum decision times are 21 and 7 times longer than those of S-PbP and PbP, respectively, which limits their application in a real-time context.

Tables \ref{TabOnlineARLarge} and \ref{TabOnlineTLarge} present the acceptance rates and decision times, respectively, obtained in the instances with $\Lambda\in\{0.8,1.5\}$. In Table \ref{TabOnlineARLarge}, the relative improvements over GP$_{_\textsf{CI}}$ are provided in parentheses. The number of sample paths in rollout algorithms is limited to $H=10$ as rollout decision times increase considerably in these large instances. Furthermore, the offline runtime required by the PFA-based policies grows to 15 minutes in the largest instances with $\Lambda=1.5$ and $K=20$. Policies S-PbP and PbP perform significantly better than the other algorithms. In addition, we see that S-PbP requires at most 11.3 seconds to compute a decision in the largest instances, which reaffirms S-PbP as a high-performing scheduling policy suitable for real-time decision-making in the DVRPSR. Rollout algorithms with $H=10$ sample paths, on the other hand, require much longer runtimes (e.g., more than 150 seconds in the largest instances) and are frequently not effective enough to outperform their base policies. Increasing the number of sample paths would lead to unacceptable online decision times (for example, up to six minutes when $H=25$). The acceptance rates and maximum decision times of the different policies are compared graphically in Figure \ref{FigEvaluationRHBig} in \ref{sec:detailed results}.

\begin{table}[!htb]
\renewcommand{\baselinestretch}{1.2}
    \caption{Comparison of Online Scheduling Policies ($\Lambda\in\{0.8,1.5\}$): Acceptance Rate (\%) and Percentage Improvement Compared to GP$_{_\textsf{CI}}$\label{TabOnlineARLarge}}
    \footnotesize
    \addtolength\tabcolsep{-0.25em}
\begin{adjustbox}{center}
\begin{threeparttable}
\begin{tabular}{cccccccccc}
	\hline
	\multicolumn{2}{c}{Instance} & \multicolumn{8}{c}{Online scheduling policy} \\\cmidrule(lr){1-2}\cmidrule(lr){3-10}
	$\Lambda$ & $K$ & \multicolumn{1}{c}{GP$_{_\textsf{CI}}$} & \multicolumn{1}{c}{GP$_{_\textsf{R}}$} & \multicolumn{1}{c}{R$_{_\textsf{R}}$-GP$_{_\textsf{CI}}(10)$} & \multicolumn{1}{c}{PFA$_{_\textsf{CI}}$} & \multicolumn{1}{c}{PFA$_{_\textsf{R}}$} & \multicolumn{1}{c}{R$_{_\textsf{R}}$-PFA$_{_\textsf{CI}}(10)$} & \multicolumn{1}{c}{S-PbP} & \multicolumn{1}{c}{PbP} \\\cmidrule(lr){1-10}
	0.8 &  6 & 30.7 & 30.8(+0.6*)                  & 32.9(+7.2)    & 37.0(+20.5)  & 36.5(+19.1) & 36.9(+20.4) & 39.7(+29.3)  & 40.5(+31.9) \\\vspace{0.5em}
        & 12 & 63.9 & 63.0($-$1.4*)                & 62.4($-$2.3)  & 65.9 (+3.1)  & 66.1 (+3.5) & 64.7 (+1.3) & 69.4 (+8.7)  & 73.8(+15.6) \\
    1.5 & 10 & 33.4 & 35.7(+7.0)                   & 34.6(+3.8)    & 39.7(+19.1)  & 39.9(+19.7) & 39.9(+19.5) & 43.1(+29.0)  & 44.2(+32.5) \\\vspace{0.5em}
        & 20 & 62.4 & 64.5(+3.4)                   & 61.6($-$1.4)  & 67.1 (+7.4)  & 66.7 (+6.8) & 65.0 (+4.1) & 70.8(+13.4)  & 76.9(+23.2) \\
\multicolumn{2}{c}{Overall} & 47.3 & 47.1($-$0.4*) & 47.7(+0.8)    & 51.5 (+8.9)  & 51.5 (+8.8) & 50.9 (+7.6) & 54.7(+15.6)  & 57.4(+21.2) \\
	\hline
\end{tabular}
{* indicates the improvements that are NOT statistically significant  (i.e., $p>0.05$ in paired two-sample \textit{t}-tests).}
\end{threeparttable}
\end{adjustbox}
    \end{table}

\begin{table}[!htb]
\renewcommand{\baselinestretch}{1.2}
    \caption{Comparison of Online Scheduling Policies ($\Lambda\in\{0.8,1.5\}$): Decision Times (seconds)\label{TabOnlineTLarge}}
    \footnotesize
    \addtolength\tabcolsep{-0.35em}
\begin{adjustbox}{center}
\begin{tabular}{ccrrrrrrrrrrrrrrrr}
	\hline
	\multicolumn{2}{c}{Instance} & \multicolumn{16}{c}{Online scheduling policy} \\\cmidrule(lr){1-2}\cmidrule(lr){3-18}
	\multirow{2}{*}{$\Lambda$} & \multirow{2}{*}{$K$} & \multicolumn{2}{c}{GP$_{_\textsf{CI}}$} & \multicolumn{2}{c}{GP$_{_\textsf{R}}$} & \multicolumn{2}{c}{R$_{_\textsf{R}}$-GP$_{_\textsf{CI}}(10)$} & \multicolumn{2}{c}{PFA$_{_\textsf{CI}}$} & \multicolumn{2}{c}{PFA$_{_\textsf{R}}$} & \multicolumn{2}{c}{R$_{_\textsf{R}}$-PFA$_{_\textsf{CI}}(10)$} & \multicolumn{2}{c}{S-PbP} & \multicolumn{2}{c}{PbP} \\\cmidrule(lr){3-4}\cmidrule(lr){5-6}\cmidrule(lr){7-8}\cmidrule(lr){9-10}\cmidrule(lr){11-12}\cmidrule(lr){13-14}\cmidrule(lr){15-16}\cmidrule(lr){17-18}
	&  & \multicolumn{1}{c}{Avg.} & \multicolumn{1}{c}{Max} & \multicolumn{1}{c}{Avg.} & \multicolumn{1}{c}{Max} & \multicolumn{1}{c}{Avg.} & \multicolumn{1}{c}{Max} & \multicolumn{1}{c}{Avg.} & \multicolumn{1}{c}{Max} & \multicolumn{1}{c}{Avg.} & \multicolumn{1}{c}{Max} & \multicolumn{1}{c}{Avg.} & \multicolumn{1}{c}{Max} & \multicolumn{1}{c}{Avg.} & \multicolumn{1}{c}{Max} & \multicolumn{1}{c}{Avg.} & \multicolumn{1}{c}{Max} \\\cmidrule(lr){1-18}
	0.8 & 6                 & $<$0.1 & $<$0.1 & 0.1 & 0.6 & 2.7 & 8.0       & $<$0.1 & $<$0.1 & 0.1 & 0.7 & 4.0 & 9.4 & 0.8 & 1.9 & 2.3 & 6.1 \\\vspace{0.5em}
        & 12                & $<$0.1 & $<$0.1 & 0.1 & 0.3 & 14.4 & 32.5     & $<$0.1 & $<$0.1 & 0.1 & 0.3 & 14.6 & 30.4 & 1.6 & 3.7 & 8.9 & 23.8 \\
    1.5 & 10                & $<$0.1 & $<$0.1 & 0.1 & 0.5 & 13.5 & 40.0     & $<$0.1 & $<$0.1 & 0.1 & 0.4 & 19.3 & 39.7 & 2.6 & 5.7 & 11.2 & 27.0 \\\vspace{0.5em}
        & 20                & $<$0.1 & $<$0.1 & 0.1 & 0.2 & 80.6 & 158.4    & $<$0.1 & $<$0.1 & 0.1 & 0.3 & 81.8 & 154.1 & 5.1 & 11.3 & 51.8 & 124.5 \\
\multicolumn{2}{c}{Overall} & $<$0.1 & $<$0.1 & 0.1 & 0.6 & 11.0 & 158.4    & $<$0.1 & $<$0.1 & 0.1 & 0.7 & 11.9 & 154.1 & 1.3 & 11.3 & 7.2 & 124.5\\
	\hline
\end{tabular}
\end{adjustbox}
    \end{table}

\RefC{2.9}{As observed in Table \ref{TabOnlineTLarge}, the average decision times of policy S-PbP grow linearly with $K$. This dependence is explained by the fact that, at each decision epoch, each additional vehicle requires the solution of an extra single-knapsack model. On the other hand, the average response time of policy PbP grows more than linearly with $K$. In the PbP, each additional vehicle requires solving a new multiple-knapsack model at each epoch, but also increases the complexity of solving models \eqref{model:mka} for the remaining vehicles (see constraints \eqref{eq:kpctr}). The computational tractability of each policy also depends on the response time requirements posed by the underlying application. However, it is clear that policy S-PbP would scale better in instances with request rate $\Lambda\geq1.5$ and many more than $K=20$ vehicles.}

\begin{figure}[!htb]
\caption{Accuracy of the Multiple-/Single-Knapsack Potential Approximation.\label{FigPotentialCurves}}
{\includegraphics[trim={100 40 90 30},clip,width=\textwidth]{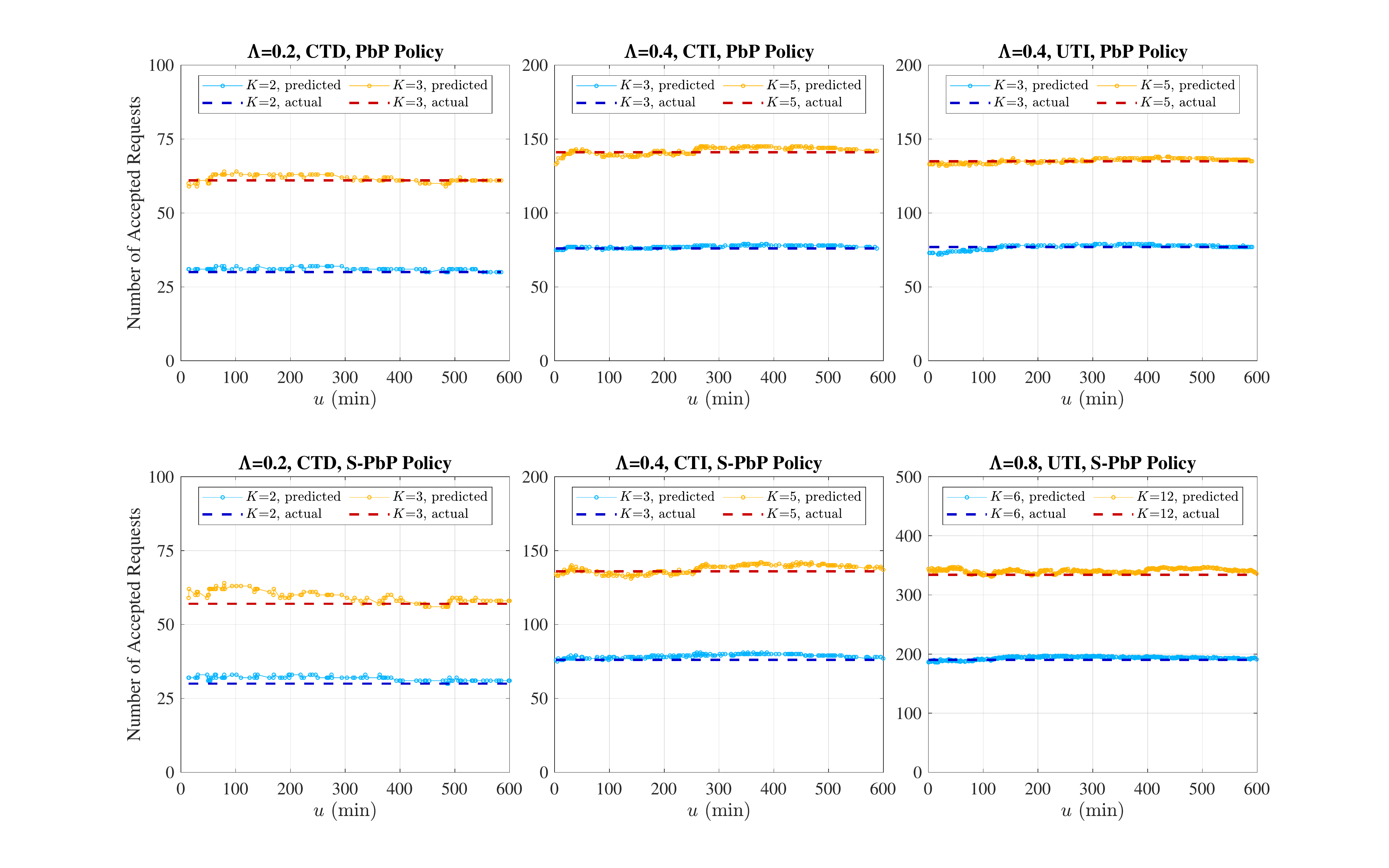}}
{}
\end{figure}

The excellent performance of PbP and S-PbP is attributed to accurate predictions on the fleet's future service capability. Figures \ref{FigPotentialCurves} shows the expected total number of accepted requests as approximated by models \eqref{model:mka} and \eqref{model:ka} throughout the service period for various request scenarios. Quite remarkably, both models are able to predict with significant accuracy, already at the beginning of the service period (i.e., before any stochasticity is disclosed) and regardless of the request distribution and instance size, the total number of accepted requests during the entire period. In combination with the decision rules applied in the optimal scheduling policy, these predictions lead to high-quality dynamic decisions that can be computed in real-time.

\begin{figure}[!htb]
\centering
\caption{Acceptance Profiles of Six Policies in Three Request Scenarios.\label{FigAcceptanceProfiles}}
\includegraphics[trim={95 10 90 20},clip,width=\textwidth]{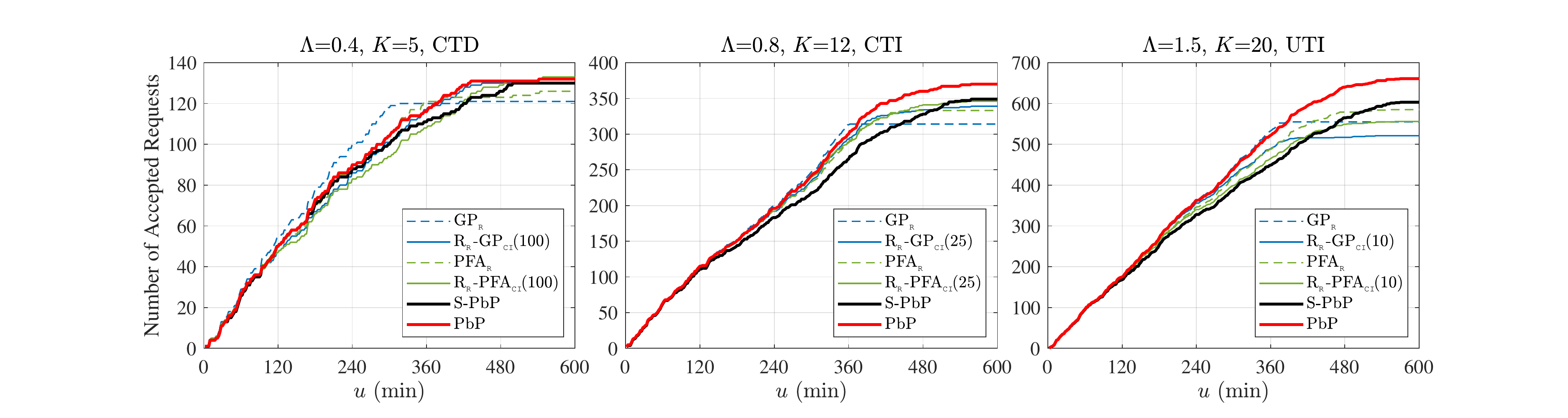}
{}
\end{figure}

Figure \ref{FigAcceptanceProfiles} exhibits the cumulative number of accepted dynamic requests over time (acceptance profiles) of several policies across three request scenarios. Even though data from only 18 simulations are plotted, these profiles' general shape represents well what is generally observed. While the greedy policy GP$_{_\textsf{R}}$ exhausts the fleet's service capability relatively early, the PFA and rollout policies with $H=25$ sample paths or more prescribe better decisions and achieve higher acceptance rates. Capacity is most efficiently consumed by potential-based policies, which accept requests at a more steady rate until near the end of the service period.

\section{Conclusions} \label{SecConclusions}
We studied a vehicle routing problem where a service provider must determine, in real-time, whether to accept dynamic pickup (or service) requests and how to adjust service routes to accommodate accepted requests. The real-time requirement limits the choice of methods effectively applied to solve large instances of the DVRPSR. In particular, traditional ADP methods such as rollout algorithms resort to expensive online simulations that compromise request response times. This paper introduced knapsack-based approximations of the reward-to-go, upon which efficient online scheduling policies and offline planning algorithms are developed. The approximation models follow simple and interpretable steps justified by observed characteristics of dynamic routes, namely, late depot arrival and request order preservation. The excellent performance of the proposed methods was demonstrated by computational experiments on very large DVRPSR instances.

Several interesting research perspectives are resulting from this work. First, the studied DVRPSR is a canonical problem that can be extended in many directions. Given that a simple linear model approximates the reward-to-go, we are optimistic that the method can be adapted and remain effective to handle side-constraints (e.g., time windows and vehicle capacity), \RefC{2.3}{as well as other sources of dynamism such as stochastic and time-dependent travel times}. \RefB{2}{The concept of effective speed is applicable to other problem variants, for example, the DVRP with stochastic pick-up and delivery requests. An interesting research direction is to adapt the knapsack approximations to pick-up and delivery variants, which yield very large DVRP instances as observed in applications such as meal delivery.} Moreover, the action space of the DVRPSR can be extended to allow inter-route reoptimization (i.e., reassignment of accepted requests among vehicles), in addition to only intra-route request reordering. Lastly, given the dramatic improvement in acceptance rates enabled by anticipative offline plans, we believe that it is also promising to investigate offline planning algorithms based on the multiple-knapsack approximation model.

\section*{Acknowledgement}
This project has received funding from the European Union's Horizon 2020 research and innovation programme under the Marie Sk\l{}odowska-Curie grant agreement No 754462. Computational experiments were performed on the Dutch national e-infrastructure with the support of SURF Cooperative. Geographical data for Vienna are copyrighted to OpenStreetMap contributors and available from \cite{OSM2020}.

\bibliography{New_Refs}

\begin{thebibliography}{41}
\expandafter\ifx\csname natexlab\endcsname\relax\def\natexlab#1{#1}\fi
\providecommand{\url}[1]{\texttt{#1}}
\providecommand{\href}[2]{#2}
\providecommand{\path}[1]{#1}
\providecommand{\DOIprefix}{doi:}
\providecommand{\ArXivprefix}{arXiv:}
\providecommand{\URLprefix}{URL: }
\providecommand{\Pubmedprefix}{pmid:}
\providecommand{\doi}[1]{\href{http://dx.doi.org/#1}{\path{#1}}}
\providecommand{\Pubmed}[1]{\href{pmid:#1}{\path{#1}}}
\providecommand{\bibinfo}[2]{#2}
\ifx\xfnm\relax \def\xfnm[#1]{\unskip,\space#1}\fi
\bibitem[{Azi et~al.(2012)Azi, Gendreau \& Potvin}]{azi2012dynamic}
\bibinfo{author}{Azi, N.}, \bibinfo{author}{Gendreau, M.}, \&
  \bibinfo{author}{Potvin, J.-Y.} (\bibinfo{year}{2012}).
\newblock \bibinfo{title}{A dynamic vehicle routing problem with multiple
  delivery routes}.
\newblock {\it \bibinfo{journal}{Annals of Operations Research}\/},  {\it
  \bibinfo{volume}{199}\/}, \bibinfo{pages}{103--112}.
\bibitem[{Baldacci et~al.(2011)Baldacci, Mingozzi \& Roberti}]{baldacci2011new}
\bibinfo{author}{Baldacci, R.}, \bibinfo{author}{Mingozzi, A.}, \&
  \bibinfo{author}{Roberti, R.} (\bibinfo{year}{2011}).
\newblock \bibinfo{title}{New route relaxation and pricing strategies for the
  vehicle routing problem}.
\newblock {\it \bibinfo{journal}{Operations Research}\/},  {\it
  \bibinfo{volume}{59}\/}, \bibinfo{pages}{1269--1283}.
\bibitem[{Barnhart et~al.(1998)Barnhart, Johnson, Nemhauser, Savelsbergh \&
  Vance}]{barnhart1998branch}
\bibinfo{author}{Barnhart, C.}, \bibinfo{author}{Johnson, E.~L.},
  \bibinfo{author}{Nemhauser, G.~L.}, \bibinfo{author}{Savelsbergh, M.~W.}, \&
  \bibinfo{author}{Vance, P.~H.} (\bibinfo{year}{1998}).
\newblock \bibinfo{title}{Branch-and-price: {C}olumn generation for solving
  huge integer programs}.
\newblock {\it \bibinfo{journal}{Operations Research}\/},  {\it
  \bibinfo{volume}{46}\/}, \bibinfo{pages}{316--329}.
\bibitem[{Bent \& Van~Hentenryck(2004)}]{bent2004scenario}
\bibinfo{author}{Bent, R.~W.}, \& \bibinfo{author}{Van~Hentenryck, P.}
  (\bibinfo{year}{2004}).
\newblock \bibinfo{title}{Scenario-based planning for partially dynamic vehicle
  routing with stochastic customers}.
\newblock {\it \bibinfo{journal}{Operations Research}\/},  {\it
  \bibinfo{volume}{52}\/}, \bibinfo{pages}{977--987}.
\bibitem[{Bertsekas \& Tsitsiklis(1996)}]{bertsekas1996neuro}
\bibinfo{author}{Bertsekas, D.~P.}, \& \bibinfo{author}{Tsitsiklis, J.~N.}
  (\bibinfo{year}{1996}).
\newblock {\it \bibinfo{title}{Neuro-dynamic programming}\/}.
\newblock \bibinfo{publisher}{Athena Scientific}.
\bibitem[{Bertsimas \& Van~Ryzin(1991)}]{bertsimas1991stochastic}
\bibinfo{author}{Bertsimas, D.~J.}, \& \bibinfo{author}{Van~Ryzin, G.}
  (\bibinfo{year}{1991}).
\newblock \bibinfo{title}{A stochastic and dynamic vehicle routing problem in
  the euclidean plane}.
\newblock {\it \bibinfo{journal}{Operations Research}\/},  {\it
  \bibinfo{volume}{39}\/}, \bibinfo{pages}{601--615}.
\bibitem[{Branke et~al.(2005)Branke, Middendorf, Noeth \&
  Dessouky}]{Branke_2005}
\bibinfo{author}{Branke, J.}, \bibinfo{author}{Middendorf, M.},
  \bibinfo{author}{Noeth, G.}, \& \bibinfo{author}{Dessouky, M.}
  (\bibinfo{year}{2005}).
\newblock \bibinfo{title}{Waiting strategies for dynamic vehicle routing}.
\newblock {\it \bibinfo{journal}{Transportation Science}\/},  {\it
  \bibinfo{volume}{39}\/}, \bibinfo{pages}{298--312}.
\bibitem[{Chen \& Xu(2006)}]{chen2006dynamic}
\bibinfo{author}{Chen, Z.-L.}, \& \bibinfo{author}{Xu, H.}
  (\bibinfo{year}{2006}).
\newblock \bibinfo{title}{Dynamic column generation for dynamic vehicle routing
  with time windows}.
\newblock {\it \bibinfo{journal}{Transportation Science}\/},  {\it
  \bibinfo{volume}{40}\/}, \bibinfo{pages}{74--88}.
\bibitem[{Costa et~al.(2019)Costa, Contardo \& Desaulniers}]{costa2019exact}
\bibinfo{author}{Costa, L.}, \bibinfo{author}{Contardo, C.}, \&
  \bibinfo{author}{Desaulniers, G.} (\bibinfo{year}{2019}).
\newblock \bibinfo{title}{Exact branch-price-and-cut algorithms for vehicle
  routing}.
\newblock {\it \bibinfo{journal}{Transportation Science}\/},  {\it
  \bibinfo{volume}{53}\/}, \bibinfo{pages}{946--985}.
\bibitem[{Daganzo(1984{\natexlab{a}})}]{daganzo1984distance}
\bibinfo{author}{Daganzo, C.~F.} (\bibinfo{year}{1984}{\natexlab{a}}).
\newblock \bibinfo{title}{The distance traveled to visit {N} points with a
  maximum of {C} stops per vehicle: An analytic model and an application}.
\newblock {\it \bibinfo{journal}{Transportation Science}\/},  {\it
  \bibinfo{volume}{18}\/}, \bibinfo{pages}{331--350}.
\bibitem[{Daganzo(1984{\natexlab{b}})}]{daganzo1984length}
\bibinfo{author}{Daganzo, C.~F.} (\bibinfo{year}{1984}{\natexlab{b}}).
\newblock \bibinfo{title}{The length of tours in zones of different shapes}.
\newblock {\it \bibinfo{journal}{Transportation Research Part B:
  Methodological}\/},  {\it \bibinfo{volume}{18}\/}, \bibinfo{pages}{135--145}.
\bibitem[{Ferrucci \& Bock(2015)}]{ferrucci2015general}
\bibinfo{author}{Ferrucci, F.}, \& \bibinfo{author}{Bock, S.}
  (\bibinfo{year}{2015}).
\newblock \bibinfo{title}{A general approach for controlling vehicle en-route
  diversions in dynamic vehicle routing problems}.
\newblock {\it \bibinfo{journal}{Transportation Research Part B:
  Methodological}\/},  {\it \bibinfo{volume}{77}\/}, \bibinfo{pages}{76--87}.
\bibitem[{Ferrucci \& Bock(2016)}]{ferrucci2016pro}
\bibinfo{author}{Ferrucci, F.}, \& \bibinfo{author}{Bock, S.}
  (\bibinfo{year}{2016}).
\newblock \bibinfo{title}{Pro-active real-time routing in applications with
  multiple request patterns}.
\newblock {\it \bibinfo{journal}{European Journal of Operational Research}\/},
  {\it \bibinfo{volume}{253}\/}, \bibinfo{pages}{356--371}.
\bibitem[{Gendreau et~al.(1999)Gendreau, Guertin, Potvin \&
  Taillard}]{gendreau1999parallel}
\bibinfo{author}{Gendreau, M.}, \bibinfo{author}{Guertin, F.},
  \bibinfo{author}{Potvin, J.-Y.}, \& \bibinfo{author}{Taillard, E.}
  (\bibinfo{year}{1999}).
\newblock \bibinfo{title}{Parallel tabu search for real-time vehicle routing
  and dispatching}.
\newblock {\it \bibinfo{journal}{Transportation Science}\/},  {\it
  \bibinfo{volume}{33}\/}, \bibinfo{pages}{381--390}.
\bibitem[{van Heeswijk et~al.(2019)van Heeswijk, Mes \&
  Schutten}]{van2019delivery}
\bibinfo{author}{van Heeswijk, W.~J.}, \bibinfo{author}{Mes, M.~R.}, \&
  \bibinfo{author}{Schutten, J.~M.} (\bibinfo{year}{2019}).
\newblock \bibinfo{title}{The delivery dispatching problem with time windows
  for urban consolidation centers}.
\newblock {\it \bibinfo{journal}{Transportation Science}\/},  {\it
  \bibinfo{volume}{53}\/}, \bibinfo{pages}{203--221}.
\bibitem[{Hvattum et~al.(2006)Hvattum, L{\o}kketangen \&
  Laporte}]{hvattum2006solving}
\bibinfo{author}{Hvattum, L.~M.}, \bibinfo{author}{L{\o}kketangen, A.}, \&
  \bibinfo{author}{Laporte, G.} (\bibinfo{year}{2006}).
\newblock \bibinfo{title}{Solving a dynamic and stochastic vehicle routing
  problem with a sample scenario hedging heuristic}.
\newblock {\it \bibinfo{journal}{Transportation Science}\/},  {\it
  \bibinfo{volume}{40}\/}, \bibinfo{pages}{421--438}.
\bibitem[{Ichoua et~al.(2006)Ichoua, Gendreau \& Potvin}]{ichoua2006exploiting}
\bibinfo{author}{Ichoua, S.}, \bibinfo{author}{Gendreau, M.}, \&
  \bibinfo{author}{Potvin, J.-Y.} (\bibinfo{year}{2006}).
\newblock \bibinfo{title}{Exploiting knowledge about future demands for
  real-time vehicle dispatching}.
\newblock {\it \bibinfo{journal}{Transportation Science}\/},  {\it
  \bibinfo{volume}{40}\/}, \bibinfo{pages}{211--225}.
\bibitem[{Ichoua et~al.(2007)Ichoua, Gendreau \& Potvin}]{Ichoua12007}
\bibinfo{author}{Ichoua, S.}, \bibinfo{author}{Gendreau, M.}, \&
  \bibinfo{author}{Potvin, J.-Y.} (\bibinfo{year}{2007}).
\newblock \bibinfo{title}{Planned route optimization for real-time vehicle
  routing}.
\newblock In \bibinfo{editor}{V.~Zeimpekis}, \bibinfo{editor}{C.~D.
  Tarantilis}, \bibinfo{editor}{G.~M. Giaglis}, \& \bibinfo{editor}{I.~Minis}
  (Eds.), {\it \bibinfo{booktitle}{Dynamic Fleet Management: Concepts, Systems,
  Algorithms {\&} Case Studies}\/} (pp. \bibinfo{pages}{1--18}).
\newblock \bibinfo{address}{Boston, MA}: \bibinfo{publisher}{Springer US}.
\bibitem[{Klapp et~al.(2018)Klapp, Erera \& Toriello}]{klapp2018dynamic}
\bibinfo{author}{Klapp, M.~A.}, \bibinfo{author}{Erera, A.~L.}, \&
  \bibinfo{author}{Toriello, A.} (\bibinfo{year}{2018}).
\newblock \bibinfo{title}{The dynamic dispatch waves problem for same-day
  delivery}.
\newblock {\it \bibinfo{journal}{European Journal of Operational Research}\/},
  {\it \bibinfo{volume}{271}\/}, \bibinfo{pages}{519--534}.
\bibitem[{Larsen et~al.(2002)Larsen, Madsen \& Solomon}]{larsen2002partially}
\bibinfo{author}{Larsen, A.}, \bibinfo{author}{Madsen, O.}, \&
  \bibinfo{author}{Solomon, M.} (\bibinfo{year}{2002}).
\newblock \bibinfo{title}{Partially dynamic vehicle routing---models and
  algorithms}.
\newblock {\it \bibinfo{journal}{Journal of the Operational Research
  Society}\/},  {\it \bibinfo{volume}{53}\/}, \bibinfo{pages}{637--646}.
\bibitem[{Mitrovi{\'{c}}-Mini{\'{c}} et~al.(2004)Mitrovi{\'{c}}-Mini{\'{c}},
  Krishnamurti \& Laporte}]{Mitrovi_Mini__2004}
\bibinfo{author}{Mitrovi{\'{c}}-Mini{\'{c}}, S.},
  \bibinfo{author}{Krishnamurti, R.}, \& \bibinfo{author}{Laporte, G.}
  (\bibinfo{year}{2004}).
\newblock \bibinfo{title}{Double-horizon based heuristics for the dynamic
  pickup and delivery problem with time windows}.
\newblock {\it \bibinfo{journal}{Transportation Research Part B:
  Methodological}\/},  {\it \bibinfo{volume}{38}\/}, \bibinfo{pages}{669--685}.
\bibitem[{{OSM}(2020)}]{OSM2020}
\bibinfo{author}{{OSM}} (\bibinfo{year}{2020}).
\newblock \bibinfo{title}{{OpenStreetMap project database, accessed August 7,
  2020, \url{http://planet.openstreetmap.org}}}.
\bibitem[{Pillac et~al.(2013)Pillac, Gendreau, Gu{\'e}ret \&
  Medaglia}]{pillac2013review}
\bibinfo{author}{Pillac, V.}, \bibinfo{author}{Gendreau, M.},
  \bibinfo{author}{Gu{\'e}ret, C.}, \& \bibinfo{author}{Medaglia, A.~L.}
  (\bibinfo{year}{2013}).
\newblock \bibinfo{title}{A review of dynamic vehicle routing problems}.
\newblock {\it \bibinfo{journal}{European Journal of Operational Research}\/},
  {\it \bibinfo{volume}{225}\/}, \bibinfo{pages}{1--11}.
\bibitem[{Powell(2011)}]{powell2011approximate}
\bibinfo{author}{Powell, W.~B.} (\bibinfo{year}{2011}).
\newblock {\it \bibinfo{title}{Approximate Dynamic Programming: Solving the
  curses of dimensionality}\/}.
\newblock Wiley Series in Probability and Statistics (\bibinfo{edition}{2nd}
  ed.).
\newblock \bibinfo{address}{Hoboken, New Jersey}: \bibinfo{publisher}{John
  Wiley \& Sons, Inc.}
\bibitem[{Powell \& Meisel(2015)}]{powell2015tutorial}
\bibinfo{author}{Powell, W.~B.}, \& \bibinfo{author}{Meisel, S.}
  (\bibinfo{year}{2015}).
\newblock \bibinfo{title}{Tutorial on stochastic optimization in energy---part
  ii: An energy storage illustration}.
\newblock {\it \bibinfo{journal}{IEEE Transactions on Power Systems}\/},  {\it
  \bibinfo{volume}{31}\/}, \bibinfo{pages}{1468--1475}.
\bibitem[{Powell et~al.(2012)Powell, Simao \&
  Bouzaiene-Ayari}]{powell2012approximate}
\bibinfo{author}{Powell, W.~B.}, \bibinfo{author}{Simao, H.~P.}, \&
  \bibinfo{author}{Bouzaiene-Ayari, B.} (\bibinfo{year}{2012}).
\newblock \bibinfo{title}{Approximate dynamic programming in transportation and
  logistics: a unified framework}.
\newblock {\it \bibinfo{journal}{EURO Journal on Transportation and
  Logistics}\/},  {\it \bibinfo{volume}{1}\/}, \bibinfo{pages}{237--284}.
\bibitem[{Psaraftis et~al.(2016)Psaraftis, Wen \&
  Kontovas}]{psaraftis2016dynamic}
\bibinfo{author}{Psaraftis, H.~N.}, \bibinfo{author}{Wen, M.}, \&
  \bibinfo{author}{Kontovas, C.~A.} (\bibinfo{year}{2016}).
\newblock \bibinfo{title}{Dynamic vehicle routing problems: Three decades and
  counting}.
\newblock {\it \bibinfo{journal}{Networks}\/},  {\it \bibinfo{volume}{67}\/},
  \bibinfo{pages}{3--31}.
\bibitem[{Ritzinger et~al.(2016)Ritzinger, Puchinger \&
  Hartl}]{ritzinger2016survey}
\bibinfo{author}{Ritzinger, U.}, \bibinfo{author}{Puchinger, J.}, \&
  \bibinfo{author}{Hartl, R.~F.} (\bibinfo{year}{2016}).
\newblock \bibinfo{title}{A survey on dynamic and stochastic vehicle routing
  problems}.
\newblock {\it \bibinfo{journal}{International Journal of Production
  Research}\/},  {\it \bibinfo{volume}{54}\/}, \bibinfo{pages}{215--231}.
\bibitem[{Soeffker et~al.(2022)Soeffker, Ulmer \&
  Mattfeld}]{soeffker2021stochastic}
\bibinfo{author}{Soeffker, N.}, \bibinfo{author}{Ulmer, M.~W.}, \&
  \bibinfo{author}{Mattfeld, D.~C.} (\bibinfo{year}{2022}).
\newblock \bibinfo{title}{Stochastic dynamic vehicle routing in the light of
  prescriptive analytics: A review}.
\newblock {\it \bibinfo{journal}{European Journal of Operational Research}\/},
  {\it \bibinfo{volume}{298}\/}, \bibinfo{pages}{801--820}.
\bibitem[{Tassiulas(1996)}]{tassiulas1996adaptive}
\bibinfo{author}{Tassiulas, L.} (\bibinfo{year}{1996}).
\newblock \bibinfo{title}{Adaptive routing on the plane}.
\newblock {\it \bibinfo{journal}{Operations Research}\/},  {\it
  \bibinfo{volume}{44}\/}, \bibinfo{pages}{823--832}.
\bibitem[{Thomas(2007)}]{thomas2007waiting}
\bibinfo{author}{Thomas, B.~W.} (\bibinfo{year}{2007}).
\newblock \bibinfo{title}{Waiting strategies for anticipating service requests
  from known customer locations}.
\newblock {\it \bibinfo{journal}{Transportation Science}\/},  {\it
  \bibinfo{volume}{41}\/}, \bibinfo{pages}{319--331}.
\bibitem[{Toth \& Vigo(2002)}]{toth2002vehicle}
\bibinfo{author}{Toth, P.}, \& \bibinfo{author}{Vigo, D.}
  (\bibinfo{year}{2002}).
\newblock {\it \bibinfo{title}{The vehicle routing problem}\/}.
\newblock \bibinfo{publisher}{SIAM}.
\bibitem[{Ulmer(2020)}]{ulmer2020dynamic}
\bibinfo{author}{Ulmer, M.~W.} (\bibinfo{year}{2020}).
\newblock \bibinfo{title}{Dynamic pricing and routing for same-day delivery}.
\newblock {\it \bibinfo{journal}{Transportation Science}\/},  {\it
  \bibinfo{volume}{54}\/}, \bibinfo{pages}{1016--1033}.
\bibitem[{Ulmer et~al.(2019)Ulmer, Goodson, Mattfeld \&
  Hennig}]{ulmer2019offline}
\bibinfo{author}{Ulmer, M.~W.}, \bibinfo{author}{Goodson, J.~C.},
  \bibinfo{author}{Mattfeld, D.~C.}, \& \bibinfo{author}{Hennig, M.}
  (\bibinfo{year}{2019}).
\newblock \bibinfo{title}{Offline--online approximate dynamic programming for
  dynamic vehicle routing with stochastic requests}.
\newblock {\it \bibinfo{journal}{Transportation Science}\/},  {\it
  \bibinfo{volume}{53}\/}, \bibinfo{pages}{185--202}.
\bibitem[{Ulmer et~al.(2020)Ulmer, Goodson, Mattfeld \&
  Thomas}]{ulmer2020modeling}
\bibinfo{author}{Ulmer, M.~W.}, \bibinfo{author}{Goodson, J.~C.},
  \bibinfo{author}{Mattfeld, D.~C.}, \& \bibinfo{author}{Thomas, B.~W.}
  (\bibinfo{year}{2020}).
\newblock \bibinfo{title}{On modeling stochastic dynamic vehicle routing
  problems}.
\newblock {\it \bibinfo{journal}{EURO Journal on Transportation and
  Logistics}\/},  {\it \bibinfo{volume}{9}\/}, \bibinfo{pages}{100008}.
\bibitem[{Ulmer et~al.(2018{\natexlab{a}})Ulmer, Mattfeld \&
  K{\"o}ster}]{ulmer2018budgeting}
\bibinfo{author}{Ulmer, M.~W.}, \bibinfo{author}{Mattfeld, D.~C.}, \&
  \bibinfo{author}{K{\"o}ster, F.} (\bibinfo{year}{2018}{\natexlab{a}}).
\newblock \bibinfo{title}{Budgeting time for dynamic vehicle routing with
  stochastic customer requests}.
\newblock {\it \bibinfo{journal}{Transportation Science}\/},  {\it
  \bibinfo{volume}{52}\/}, \bibinfo{pages}{20--37}.
\bibitem[{Ulmer et~al.(2018{\natexlab{b}})Ulmer, Soeffker \&
  Mattfeld}]{ulmer2018value}
\bibinfo{author}{Ulmer, M.~W.}, \bibinfo{author}{Soeffker, N.}, \&
  \bibinfo{author}{Mattfeld, D.~C.} (\bibinfo{year}{2018}{\natexlab{b}}).
\newblock \bibinfo{title}{Value function approximation for dynamic multi-period
  vehicle routing}.
\newblock {\it \bibinfo{journal}{European Journal of Operational Research}\/},
  {\it \bibinfo{volume}{269}\/}, \bibinfo{pages}{883--899}.
\bibitem[{Ulmer \& Streng(2019)}]{ulmer2019same}
\bibinfo{author}{Ulmer, M.~W.}, \& \bibinfo{author}{Streng, S.}
  (\bibinfo{year}{2019}).
\newblock \bibinfo{title}{Same-day delivery with pickup stations and autonomous
  vehicles}.
\newblock {\it \bibinfo{journal}{Computers \& Operations Research}\/},  {\it
  \bibinfo{volume}{108}\/}, \bibinfo{pages}{1--19}.
\bibitem[{Ulmer \& Thomas(2018)}]{ulmer2018same}
\bibinfo{author}{Ulmer, M.~W.}, \& \bibinfo{author}{Thomas, B.~W.}
  (\bibinfo{year}{2018}).
\newblock \bibinfo{title}{Same-day delivery with heterogeneous fleets of drones
  and vehicles}.
\newblock {\it \bibinfo{journal}{Networks}\/},  {\it \bibinfo{volume}{72}\/},
  \bibinfo{pages}{475--505}.
\bibitem[{Ulmer \& Thomas(2020)}]{ulmer2020meso}
\bibinfo{author}{Ulmer, M.~W.}, \& \bibinfo{author}{Thomas, B.~W.}
  (\bibinfo{year}{2020}).
\newblock \bibinfo{title}{Meso-parametric value function approximation for
  dynamic customer acceptances in delivery routing}.
\newblock {\it \bibinfo{journal}{European Journal of Operational Research}\/},
  {\it \bibinfo{volume}{285}\/}, \bibinfo{pages}{183--195}.
\bibitem[{Voccia et~al.(2019)Voccia, Campbell \& Thomas}]{voccia2019same}
\bibinfo{author}{Voccia, S.~A.}, \bibinfo{author}{Campbell, A.~M.}, \&
  \bibinfo{author}{Thomas, B.~W.} (\bibinfo{year}{2019}).
\newblock \bibinfo{title}{The same-day delivery problem for online purchases}.
\newblock {\it \bibinfo{journal}{Transportation Science}\/},  {\it
  \bibinfo{volume}{53}\/}, \bibinfo{pages}{167--184}.

\end{thebibliography}

\appendix

\section{Summary of Notations}\label{sec:notations}
\setcounter{table}{0}
\RefA{4}{Tables \ref{tab:notation_pro} and \ref{tab:notation_alg} summarize the problem notations and algorithm notations defined in Section~\ref{SecProblem} and Section~\ref{SecOnline}, respectively.}

\begin{table}[!ht]
    \footnotesize
    \centering
    \caption{Problem Notations}\label{tab:notation_pro}
    \begin{tabular}{ll}
    \hline
         Notation & Description \\
    \hline
         \multicolumn{2}{l}{Global parameters}\\
    \hline
         $\mathcal{G}$                                                   & Graph (street network)                                  \\
         $\mathcal{V}$                                                   & Set of nodes (road intersections)                       \\
         $\mathcal{A}$                                                   & Set of arcs (road segments)                             \\
         $d_{ij}$                                                        & Length of arc $(i,j)\in\mathcal{A}$                     \\
         $\overline{s}$                                                  & Vehicle speed                                           \\
         $t_{ij}=d_{ij}/\overline{s}$                                    & Travel time on arc $(i,j)\in\mathcal{A}$                \\
         $t(i,j)$                                                        & Fastest travel time between nodes $i,j\in\mathcal{V}$   \\
         $[0,U]$                                                         & Service period                                          \\
         $K$                                                             & Number of vehicles                                      \\
         $r=(u,i,d)$                                                     & Customer request                                        \\
         $u\in[0,U]$                                                     & Request arrival time                                    \\
         $i\in\mathcal{V}\setminus\{0\}$                                 & Customer node                                           \\
         $d\in\mathbb{R}_{>0}$                                          & Service time of request                                 \\
         $\mathcal{S}$                                                   & Set of static requests                                  \\
         $\mathcal{D}(u)$                                                & Ordered set of dynamic requests generated up to instant $u$     \\
         $\lambda_{i}(u)$                                                & Request arrival rate at node $i$ and instant $u$        \\
         $\Lambda(u)=\sum_{i\in\mathcal{V}\setminus\{0\}}\lambda_{i}(u)$ & Overall request arrival rate at instant $u$             \\
    \hline
        \multicolumn{2}{l}{State variables}\\
    \hline
        $T=|\mathcal{D}(U)|$                                        & Number of dynamic requests         \\
        $\theta=((v_{0},\delta_{0}),\ldots,(v_{f},\delta_{f}))$     & Route                                                                 \\
        $v_{i}\in\mathcal{V}$, $i\in\{1,\ldots,f-1\}$               & The $i$-th node along route $\theta$                    \\
        $\delta_{i}$, $i\in\{1,\ldots,f-1\}$                        & Set of requests to serve when arriving at node $v_{i}$                     \\
        $\Delta(\theta)=\cup_{i=1}^{f-1}\delta_{i}$                 & Set of all requests served by route $\theta$                          \\ 
        $\Delta_{u}(\theta)=\{(u',i,d)\in\Delta(\theta):u'\geq u\}$ & Ordered set of requests served by route $\theta$ during $[u,U]$             \\
        $b(\tau,\theta)$                                            & Budget of route $\theta$ when it starts at instant $\tau$               \\
        $V_{k}=(\tau_{k},\theta_{k})$                               & State of vehicle $k$ when it starts route $\theta_{k}$ at instant $\tau_{k}$    \\
        $S_{t}=\{V_{1},\ldots,V_{K},r_t\}$                          & State variable at decision epoch $t=\{1,\ldots,T\}$                   \\
    \hline
        \multicolumn{2}{l}{Decision variables}\\
    \hline
        $\mathcal{F}$                           & Feasible decision set at initial state $S_{0}$              \\
        $\mathbf{y}\in\mathcal{F}$              & Offline decision (set of planned routes)                    \\
        $\mathcal{X}_t$                         & Feasible decision set at state $S_{t}$, $t\neq0$            \\
        $x_{t}\in\mathcal{X}_t$                 & Online decision at state $S_{t}$, $t\neq0$                  \\
        $S_{t+1}=S^{M}(S_{t},x_{t},r_{t+1})$    & State transition function                                   \\
    \hline
        \multicolumn{2}{l}{Reward \& objective function}\\
    \hline
        $R(S_{t},x_{t})$                    & Reward of executing decision $x_{t}$ at state $S_{t}$             \\
        $\Pi$                               & Set of feasible policies                                          \\
        $X^{\pi}(S_{t})\in\mathcal{X}_t$    & Decision prescribed by policy $\pi\in\Pi$ at state $S_{t}$        \\
        $\Phi_{\pi}(S_{t})$                 & Potential of state $S_{t}$ when controlled by policy $\pi\in\Pi$  \\
    \hline
    \end{tabular}
\end{table}

\begin{table}[!ht]
    \footnotesize
    \centering
    \caption{Algorithm Notations}
    \label{tab:notation_alg}
    \begin{tabular}{ll}
    \hline
         Notation & Description \\
    \hline
         $\rho$                                                                                         & Routing policy                                                                        \\
         $\rho(\theta,r)$                                                                               & New route resulting from adding request $r$ to route $\theta$ by policy $\rho$        \\
         $\pi^*\in\Pi$                                                                                  & Optimal scheduling policy                                                               \\
         $x_{t}^{-}\in\mathcal{X}_{t}$                                                                  & ``Reject'' decision with respect to request $r_{t}$                                   \\
         $x_{t}^{+}\in\mathcal{X}_{t}\setminus\{x_{t}^{-}\}$                                            & ``Accept'' decision with respect to request $r_{t}$                                   \\
         $\mathcal{X}_{t}^{k}\subset\mathcal{X}_{t}$                                                    & Set of decisions where $r_{t}$ is accepted and assigned to vehicle $k$                \\
         $x_{t}^{k}\in\mathcal{X}_{t}^{k}$                                                              & Decision of accepting $r_{t}$ and assigning it to vehicle $k$                         \\
         $\check{s}_{u_{t}}(V_{k})$                                                                     & Effective speed of a vehicle with state $V_{k}$ estimated at instant $u_{t}$          \\
         $d_{u_{t}}(V_{k})$                                                                             & Remaining distance of route $\theta_{k}$ at instant $u_{t}$                           \\
         $\mathcal{V}_{u_{t}}(V_{k},u')\subseteq\mathcal{V}$                                            & Set of nodes visited by $\theta_{k}$ during $[u',U]$, estimated at instant $u_{t}$    \\
         $c_{u_{t}}(V_{k},r')$                                                                          & Cost of assigning future request $r'$ to vehicle $k$, estimated at instant $u_{t}$    \\
         $\omega=\{r^{\omega}_{1},\ldots,r^{\omega}_{T_{\omega}}\}$                                     & Sample path of the request arrival process $\{\Lambda(u)\}_{u\in[0,U]}$               \\
         $H$                                                                                            & Number of sample paths                                                                \\
         $\Omega=\{\omega_{1},\ldots,\omega_{H}\}$                                                      & Set of sample paths                                                                   \\
         $\overline{K}=\{k\in\{1,\ldots,K\}:V_{k}\neq\emptyset\}$                                       & Set of non-idle vehicles                                                              \\
         $\phi^{\omega}_{\pi^*}(S_{t+1}|x_{t})$                                                         & Multiple-knapsack approximation of potential in sample path $\omega$                  \\
         $p^{\omega}_{u_{t}}(V_{k})$                                                                    & Single-knapsack approximation of vehicle $k$'s potential in sample path $\omega$      \\
         $\alpha^{\omega}$                                                                              & Compensation ratio of single-knapsack approximation                                   \\ 
         $\hat{\Phi}_{\pi^{*}}(S_{t+1}|x_{t})\approx\mathbb{E}[\Phi_{\pi^{*}}(S_{t+1})|x_{t}]$          & Multiple-/Single-knapsack approximation of the expected potential                     \\
         $\hat{\Phi}_{_\textsf{PFA}}(S_{t+1}|x_{t})\approx\mathbb{E}[\Phi_{\pi^{*}}(S_{t+1})|x_{t}]$    & PFA approximation of the expected potential                                           \\
    \hline
    \end{tabular}
    
\end{table}

\section{Auxiliary Algorithms}
\subsection{Predicting Vehicle Locations}\label{sec:algpredict}
Consider a state $S_{t}=(V_{1},\ldots,V_{K},r_{t})$, $r_{t}=(u_{t},i_{t},d_{t})$. The following procedure is used to predict, at the current instant $u_{t}$, the location of a non-idle vehicle $k$, $V_{k}=(\tau_{k},\theta_{k})$, at a future instant $u'>u_{t}$, and to determine, also at the current instant $u_{t}$, the set $\mathcal{V}_{u_{t}}(V_{k},u')$ of nodes along $\theta_{k}$ to be traversed after $u'$:

\noindent\emph{Step 1: Predicting the location of vehicle $k$ at a future instant $u'$:}
\begin{enumerate}[label=(\roman*)]
\item \RefA{4}{Let $d_{u_{t}}(V_{k})$ be the remaining distance to be traveled along route $\theta_{k}$ at instant $u_{t}$. Then, the effective speed $\check{s}_{u_{t}}(V_{k})$ is given by:}
\begin{equation}\nonumber
\check{s}_{u_{t}}(V_{k})=\frac{d_{u_{t}}(V_{k})\overline{s}}{d_{u_{t}}(V_{k})+b(\tau_{k},\theta_{k})\overline{s}}.
\end{equation}
\item Let $\theta_{k}(t)=((v_{t},\delta_{t}),\ldots,(v_{f},\delta_{f}))$ be the substring of $\theta_{k}$ such that $v_{t}$ is the current node at which vehicle $k$ is located (if the vehicle is currently in service) or the next node to be reached by vehicle $k$ (if the vehicle is currently in transit).
\item Let $\tilde{u}_{s}$ be the elapsed (service) time since vehicle $k$ arrived at node $v_{t}$, if the vehicle is in service, $\tilde{u}_{s}=0$ otherwise. Further, let $\tilde{u}_{t}$ be the remaining (transit) time under the effective speed $\check{s}_{u_{t}}(V_{k})$ before vehicle $k$ reaches node $v_{t}$, if the vehicle is in transit, $\tilde{u}_{t}=0$ otherwise.
\item Taking into account the initial conditions set by $\tilde{u}_{s}$ and $\tilde{u}_{t}$, the remaining route $\theta_{k}(t)$ is simulated under effective speed $\check{s}_{u_{t}}(V_{k})$ for an amount of time $u'-u_{t}$.
\end{enumerate}

\noindent\emph{Step 2: Determining $\mathcal{V}_{u_{t}}(V_{k},u')$:}
\begin{enumerate}[label=(\roman*)]\setcounter{enumi}{3}
\item Let $\theta_{k}'(t)=((v_{t}',\delta_{t}'),\ldots,(v_{f},\delta_{f}))$ be the substring of $\theta_{k}(t)$ such that $v_{t}'$ is the (predicted) first unvisited node along route $\theta_{k}$ at instant $u'$.
\item The procedure returns $\mathcal{V}_{u_{t}}(V_{k},u')=\{v:(v,\delta)\in\theta_{k}'(t)\}$.
\end{enumerate}

\subsection{Column Generation for the DCVRP}\label{sec:colgen}
Model \eqref{model:spm} under cost function \eqref{eq:myo} corresponds to a DCVRP defined on graph $\mathcal{G}'=(\mathcal{V}', \mathcal{A}')$. Let $\tilde{\Theta}$ be the set of all feasible DCVRP routes, where each route $\tilde{\theta}\in\tilde{\Theta}$ is a non-empty elementary sequence of nodes $\tilde{\theta}=(\tilde{v}_{1},\ldots,\tilde{v}_{f})$, $\tilde{v}_{1},\ldots,\tilde{v}_{f}\in\mathcal{V}'\setminus\{0\}$. The cost of a route $\tilde{\theta}=(\tilde{v}_{1},\ldots,\tilde{v}_{f})$ is given by $C_{\tilde{\theta}}=\sum_{i=2}^{f}c_{v_{i-1}{v_{i}}}$, where $c_{ij}$ is the cost of arc $(i,j)\in\mathcal{A}'$ as specified in Section \ref{sec:budplanner}. Introducing binary variables $z_{\tilde{\theta}}$, $\tilde{\theta}\in\tilde{\Theta}$, to indicate the routes in the solution, the DCVRP is rewritten as follows:
\begin{align}
&\min&&\sum_{\tilde{\theta}\in\tilde{\Theta}}C_{\tilde{\theta}}z_{\tilde{\theta}},\nonumber\\
&\text{s.t.}&&\sum_{\tilde{\theta}\in\tilde{\Theta}}[i\in\tilde{\theta}]z_{\tilde{\theta}}=1,&i\in\mathcal{V}'\setminus\{0\},\nonumber\\
&&&\sum_{\tilde{\theta}\in\tilde{\Theta}}z_{\tilde{\theta}}\leq K.\label{model:dcvrp}\tag{DCVRP-SPM}
\end{align}

The linear relaxation of \eqref{model:dcvrp} cannot be solved directly because the number of variables $|\tilde{\Theta}|$ grows exponentially with $|\mathcal{V}'|$. So, we apply column generation and, in the process, generate a large pool of routes to subsequently use in the offline planning heuristic. Within column generation, the restricted master problem (RMP) refers to the linear relaxation of \eqref{model:dcvrp} when only a subset $\tilde{\Theta}_{_\textsf{R}}\subset\tilde{\Theta}$ of variables is considered. A profitable variable to add to $\tilde{\Theta}_{_\textsf{R}}$ is identified by solving the pricing problem, which calls for a column $\tilde{\theta}$ with negative reduced cost:
\begin{equation}\label{dcvrp:pp}
C_{\tilde{\theta}}-\sum_{i\in\mathcal{V}'\setminus\{0\}}[i\in\tilde{\theta}]\beta_{i}-\beta_{0}<0,
\end{equation}
where $\beta_{i}$ ($i\in\mathcal{V}'\setminus\{0\}$) and $\beta_{0}$ are, respectively, the dual values associated with the partitioning and the maximum number of vehicles constraints in \eqref{model:dcvrp}, obtained by solving the RMP. Once profitable variables are identified, they are added to $\tilde{\Theta}_{_\textsf{R}}$ and the procedure iterates. Convergence is attained when no further variable satisfies \eqref{dcvrp:pp}.

The pricing problem is solved with a labeling procedure, as common in algorithms for several VRP variants \citep{costa2019exact}. Each label corresponds to a partial path from node $0$ to a node $j\in\mathcal{V}'\setminus\{0\}$, and stores the reduced cost accumulated so far along the partial path. To ensure tractability of the labeling procedure, we allow non-elementary routes according to \emph{ng-}route relaxation \citep{baldacci2011new}, and control the combinatorial growth of labels with dominance rules.

\section{PFA Policy}\label{sec:PFA}
Given a state $S_{t}=(V_{1},\ldots,V_{K},r_{t})$, $r_{t}=(u_{t},i_{t},d_{t})$, it is clear that $r_t$ should be accepted only if this acceptance leads to an increase in the expected total reward, as stated by decision rule \eqref{eq:dr1}. Hence, we design a PFA policy with an analytical function $\hat{\Phi}_{_\textsf{PFA}}(S_{t+1}|x_t)$, which approximates the expected reward-to-go $\mathbb{E}[\Phi_{\pi^*}(S_{t+1})|x_t]$ based on the total time budget:
\begin{equation}\label{eq:pfa}
    \hat{\Phi}_{_\textsf{PFA}}(S_{t+1}|x_t)=\left\{\begin{aligned}
        & \gamma\sum_{V_{k}\neq\emptyset}b(\tau_k,\theta_k),&&\text{if }x_t=x_{t}^{-}\in\mathcal{X}_{t},\\
        &\gamma\left(b(\tau_k,\rho(\theta_k,r_t))+\sum_{i\neq k,V_{i}\neq\emptyset}b(\tau_i,\theta_i)\right),&&\text{if }x_t\in\mathcal{X}_{t}^{k},
    \end{aligned}\right.
\end{equation}
where $x_{t}^{-}$ is the `reject' decision with respect to $r_t$, $\mathcal{X}_{t}^{k}$ is the set of `accept' decisions where $r_t$ is assigned to vehicle $k$, and coefficient \RefC{2.6}{$\gamma>0$} is a relative weight that trades off the immediate reward and the reward-to-go. The larger $\gamma$ is, the more likely request $r_t$ is rejected so as to preserve the fleet's future service capability.

Parameter $\gamma$ is tuned in advance via offline simulations. As the optimal value of $\gamma$ depends not only on instance settings (graph type, fleet size, request rate and distribution, etc.) but also on the initial route plan, the tuning process of $\gamma$ takes place after the initial routes have been fixed by an offline planner and before the beginning of the service period. During the tuning process, we simulate PFA policies with different values of $\gamma$ over $H$ complete sample paths. The values of $\gamma$ simulated are selected based on a variable step size search procedure, \RefC{2.6}{which assumes that the total reward is concave on $\gamma$, as observed empirically.} Finally, we select the value of $\gamma$ which leads to the highest simulated total reward.

\section{Detailed Computational Results}\label{sec:detailed results}
\RefA{8}{This section presents the detailed computational results with respect to the three different request distributions: UTI, CTI, and CTD. The results shown in Tables \ref{TabOfflineARDetail}, \ref{TabOnlineARSmallDetail}, \ref{TabOnlineDTSmallDetail}, \ref{TabOnlineARLargeDetail} and \ref{TabOnlineTLargeDetail} correspond to those shown in Tables \ref{TabOfflineAR}, \ref{TabOnlineARSmall}, \ref{TabOnlineDTSmall}, \ref{TabOnlineARLarge} and \ref{TabOnlineTLarge} in Section \ref{SecComputation}, respectively. Moreover, Figures \ref{FigEvaluationRCI} and \ref{FigEvaluationRHBig} compare the performances of online scheduling policies in small ($\Lambda\in\{0.2,0.4\}$) and large ($\Lambda\in\{0.8,1.5\}$) instances, respectively.}

\setcounter{table}{0}
\setcounter{figure}{0}

\begin{table}[!htb]
\begin{center}
    \renewcommand{\baselinestretch}{1.2}
    \caption{Comparison of Offline Planning Algorithms: Acceptance Rate (\%) and Relative Improvement (\%)\label{TabOfflineARDetail}}
\footnotesize
\begin{threeparttable}
\addtolength\tabcolsep{-0.15em}
\begin{tabular}{cccrrrrrrrrrrrrrrrr}
	\hline
	\multicolumn{3}{c}{Instance} & & \multicolumn{3}{c}{GP$_{_\textsf{CI}}$} &  & \multicolumn{3}{c}{GP$_{_\textsf{R}}$} &  & \multicolumn{3}{c}{R$_{_\textsf{CI}}$-GP$_{_\textsf{CI}}$(50)} &  & \multicolumn{3}{c}{R$_{_\textsf{R}}$-GP$_{_\textsf{CI}}$(50)} \\\cmidrule(lr){1-3}\cmidrule(lr){5-7}\cmidrule(lr){9-11}\cmidrule(lr){13-15}\cmidrule(lr){17-19}
	$\Lambda$ & $K$ & Dist. &  & \multicolumn{1}{c}{MY} & \multicolumn{1}{c}{PB} & \multicolumn{1}{c}{Impr.} &  & \multicolumn{1}{c}{MY} & \multicolumn{1}{c}{PB} & \multicolumn{1}{c}{Impr.} &  & \multicolumn{1}{c}{MY} & \multicolumn{1}{c}{PB} & \multicolumn{1}{c}{Impr.} &  & \multicolumn{1}{c}{MY} & \multicolumn{1}{c}{PB} & \multicolumn{1}{c}{Impr.} \\\cmidrule(lr){1-19}
	0.2 & 2 & UTI & & 19.9 & 19.9 & 0.0 &   & 20.3 & 20.3 & 0.0 &   & 24.2 & 24.4 & +0.6* &     & 24.4 & 24.4 & $-$0.1* \\
	&  & CTI &      & 20.0 & 20.0 & 0.0 &   & 19.8 & 19.8 & 0.0 &   & 24.2 & 24.4 & +0.7* &     & 24.5 & 24.5 & 0.0 \\\vspace{0.5em}
	&  & CTD &      & 18.6 & 18.6 & 0.0 &   & 19.5 & 19.5 & 0.0 &   & 24.1 & 24.1 & $-$0.1* &   & 24.1 & 24.1 & $-$0.3* \\
	& 3 & UTI &     & 36.6 & 42.1 & +15.1 & & 37.6 & 44.9 & +19.2 & & 47.4 & 49.4 & +4.3 &      & 47.6 & 51.0 & +7.2 \\
	&  & CTI &      & 38.4 & 43.4 & +12.8 & & 39.1 & 44.7 & +14.2 & & 49.8 & 50.6 & +1.7 &      & 49.8 & 52.4 & +5.3 \\\vspace{0.5em}
	&  & CTD &      & 42.2 & 45.2 & +7.0 &  & 42.0 & 46.8 & +11.4 & & 50.8 & 51.6 & +1.6 &      & 51.4 & 52.5 & +2.2 \\
	0.4 & 3 & UTI & & 20.5 & 24.2 & +17.7 & & 20.8 & 25.2 & +20.9 & & 26.5 & 28.6 & +8.0 &      & 26.6 & 29.5 & +11.1 \\
	&  & CTI &      & 21.2 & 21.9 & +3.1* & & 21.1 & 23.4 & +11.0 & & 26.2 & 27.8 & +6.3 &      & 26.6 & 28.5 & +7.0 \\\vspace{0.5em}
	&  & CTD &      & 18.2 & 22.9 & +25.8 & & 19.1 & 22.2 & +16.6 & & 25.2 & 27.6 & +9.4 &      & 25.2 & 27.8 & +10.4 \\
	& 5 & UTI &     & 42.8 & 50.0 & +16.7 & & 43.1 & 52.5 & +21.6 & & 51.4 & 56.0 & +8.9 &      & 51.8 & 56.3 & +8.8 \\
	&  & CTI &      & 43.7 & 49.9 & +14.4 & & 43.4 & 50.1 & +15.5 & & 52.4 & 56.0 & +6.7 &      & 52.5 & 56.5 & +7.6 \\\vspace{0.5em}
	&  & CTD &      & 41.3 & 47.8 & +15.8 & & 42.8 & 47.9 & +11.9 & & 50.5 & 53.5 & +6.0 &      & 50.6 & 54.2 & +7.2 \\
	\multicolumn{3}{c}{Overall} &  & 30.3 & 33.8 & +11.7 &  & 30.7 & 34.8 & +13.2 &  & 37.7     & 39.5 & +4.7 &  & 37.9 & 40.1 & +5.9 \\
	\cmidrule(lr){1-19}
	\multicolumn{3}{c}{Instance} & & \multicolumn{3}{c}{PFA$_{_\textsf{CI}}$} &  & \multicolumn{3}{c}{PFA$_{_\textsf{R}}$} &  & \multicolumn{3}{c}{S-PbP} &  & \multicolumn{3}{c}{PbP} \\\cmidrule(lr){1-3}\cmidrule(lr){5-7}\cmidrule(lr){9-11}\cmidrule(lr){13-15}\cmidrule(lr){17-19}
	$\Lambda$ & $K$ & Dist. &  & \multicolumn{1}{c}{MY} & \multicolumn{1}{c}{PB} & \multicolumn{1}{c}{Impr.} &  & \multicolumn{1}{c}{MY} & \multicolumn{1}{c}{PB} & \multicolumn{1}{c}{Impr.} &  & \multicolumn{1}{c}{MY} & \multicolumn{1}{c}{PB} & \multicolumn{1}{c}{Impr.} &  & \multicolumn{1}{c}{MY} & \multicolumn{1}{c}{PB} & \multicolumn{1}{c}{Impr.} \\\cmidrule(lr){1-19}
	0.2 & 2 & UTI & & 26.9 & 26.8 & $-$0.3* &   & 26.9 & 26.9 & +0.1* &     & 27.8 & 27.7 & $-$0.2* &   & 27.8 & 27.7 & $-$0.4 \\
	&  & CTI &      & 26.0 & 26.0 & $-$0.2* &   & 26.4 & 26.5 & +0.3* &     & 27.6 & 27.7 & +0.3* &     & 27.5 & 27.7 & +0.6 \\\vspace{0.5em}
	&  & CTD &      & 27.0 & 27.0 & $-$0.1* &   & 27.2 & 27.1 & $-$0.3*&    & 27.8 & 27.8 & 0.0 &       & 27.8 & 27.8 & 0.0 \\
	& 3 & UTI &     & 40.3 & 48.6 & +20.6 &     & 43.1 & 51.0 & +18.4 &     & 51.0 & 53.3 & +4.5 &      & 51.4 & 53.7 & +4.5 \\
	&  & CTI &      & 41.5 & 49.5 & +19.3 &     & 42.5 & 50.9 & +19.8 &     & 52.6 & 52.3 & $-$0.5* &   & 52.7 & 54.2 & +2.8 \\\vspace{0.5em}
	&  & CTD &      & 43.0 & 47.4 & +10.3 &     & 42.7 & 48.8 & +14.3 &     & 52.4 & 52.4 & 0.0 &       & 53.2 & 55.0 & +3.2 \\
	0.4 & 3 & UTI & & 24.1 & 30.3 & +25.9 &     & 24.6 & 31.6 & +28.2 &     & 29.4 & 32.0 & +8.5 &      & 29.9 & 32.5 & +8.8 \\
	&  & CTI &      & 24.0 & 30.0 & +24.7 &     & 24.1 & 29.7 & +23.3 &     & 29.4 & 32.1 & +9.4 &      & 29.4 & 32.4 & +10.1 \\\vspace{0.5em}
	&  & CTD &      & 23.7 & 29.1 & +22.7 &     & 23.9 & 28.7 & +20.2 &     & 27.4 & 30.1 & +9.9 &      & 27.4 & 30.5 & +11.3 \\
	& 5 & UTI &     & 46.3 & 54.9 & +18.7 &     & 45.8 & 54.2 & +18.2 &     & 53.1 & 55.8 & +5.2 &      & 54.8 & 59.0 & +7.6 \\
	&  & CTI &      & 47.0 & 55.5 & +18.1 &     & 47.5 & 55.8 & +17.6 &     & 53.9 & 55.8 & +3.5 &      & 55.3 & 58.5 & +5.8 \\\vspace{0.5em}
	&  & CTD &      & 44.4 & 52.7 & +18.6 &     & 45.5 & 53.3 & +17.3 &     & 53.1 & 55.6 & +4.9 &      & 53.5 & 56.5 & +5.7 \\
	\multicolumn{3}{c}{Overall} &  & 34.5 & 39.8 & +15.3 &  & 35.0 & 40.4 & +15.3 &  & 40.5 & 41.9 & +3.6 &  & 40.9 & 43.0 & +5.0\\
	\hline
\end{tabular}
\footnotesize\textit{Notes.} Dist., request distribution; Impr., relative improvement of PB over MY: (PB$-$MY)/MY$\times\%$.\\
{* indicates the improvements that are NOT statistically significant  (i.e., $p>0.05$ in paired two-sample \textit{t}-tests).}
\end{threeparttable}
\end{center}
    
\end{table}

\begin{table}[!htb]
\renewcommand{\baselinestretch}{1.2}
    \caption{Comparison of Online Scheduling Policies ($\Lambda\in\{0.2,0.4\}$): Percentage Improvement in the Acceptance Rate Compared to GP$_{_\textsf{CI}}$\label{TabOnlineARSmallDetail}}
    \footnotesize
\addtolength\tabcolsep{-0.1em}
\begin{adjustbox}{center}
\begin{threeparttable}
\begin{tabular}{cccrrrrrrrrrrr}
	\hline
	\multicolumn{3}{c}{Instance} & \multicolumn{11}{c}{Online scheduling policy} \\\cmidrule(lr){1-3}\cmidrule(lr){4-14}
	\multirow{2}{*}{$\Lambda$} & \multirow{2}{*}{$K$} & \multirow{2}{*}{Dist.} & \multicolumn{1}{c}{\multirow{2}{*}{GP$_{_\textsf{R}}$}} & \multicolumn{3}{c}{R$_{_\textsf{R}}$-GP$_{_\textsf{CI}}(H)$} & \multicolumn{1}{c}{\multirow{2}{*}{PFA$_{_\textsf{CI}}$}} & \multicolumn{1}{c}{\multirow{2}{*}{PFA$_{_\textsf{R}}$}} & \multicolumn{3}{c}{R$_{_\textsf{R}}$-PFA$_{_\textsf{CI}}(H)$} & \multicolumn{1}{c}{\multirow{2}{*}{S-PbP}} & \multicolumn{1}{c}{\multirow{2}{*}{PbP}} \\\cmidrule(lr){5-7}\cmidrule(lr){10-12}
	& & & & \multicolumn{1}{c}{$H$=25} & \multicolumn{1}{c}{$H$=50} & \multicolumn{1}{c}{$H$=100} & & & \multicolumn{1}{c}{$H$=25} & \multicolumn{1}{c}{$H$=50} & \multicolumn{1}{c}{$H$=100} & & \\\cmidrule(lr){1-14}
	0.2 & 2 & UTI   & +2.0* & +21.2 & +22.7 & +23.4 & +34.6 & +35.1 & +36.5 & +36.9 & +37.3 & +39.2 & +39.3 \\
	&  & CTI        & $-$0.9* & +21.2 & +22.8 & +23.4 & +30.0 & +32.7 & +36.1 & +37.6 & +37.6 & +38.9 & +38.5 \\\vspace{0.5em}
	&  & CTD        & +4.7* & +28.3 & +29.6 & +30.6 & +45.1 & +45.9 & +46.7 & +48.3 & +48.3 & +49.6 & +49.9 \\
	& 3 & UTI       & +6.5 & +19.7 & +21.0 & +22.3 & +15.3 & +21.2 & +23.7 & +24.9 & +26.6 & +26.5 & +27.6 \\
	&  & CTI        & +3.0 & +18.7 & +20.8 & +22.5 & +14.2 & +17.3 & +23.0 & +24.4 & +25.4 & +20.7 & +25.0 \\\vspace{0.5em}
	&  & CTD        & +3.5* & +13.8 & +16.3 & +17.8 & +5.0 & +8.1 & +16.8 & +18.3 & +19.1 & +15.9 & +21.7 \\
	0.4 & 3 & UTI   & +4.1* & +18.8 & +22.1 & +23.8 & +25.4 & +30.7 & +30.1 & +32.3 & +32.5 & +32.2 & +34.4 \\
	&  & CTI        & +7.2 & +28.3 & +30.2 & +31.2 & +37.2 & +35.9 & +45.5 & +47.3 & +48.7 & +46.9 & +48.3 \\\vspace{0.5em}
	&  & CTD        & $-$2.8* & +18.6 & +21.7 & +23.6 & +27.5 & +25.8 & +29.1 & +31.0 & +32.0 & +31.9 & +33.5 \\
	& 5 & UTI       & +5.0 & +10.3 & +12.7 & +14.9 & +10.0 & +8.4 & +13.9 & +16.3 & +17.6 & +11.8 & +18.0 \\
	&  & CTI        & +0.3* & +10.2 & +13.2 & +15.5 & +11.0 & +11.8 & +14.1 & +16.6 & +17.6 & +11.7 & +17.1 \\\vspace{0.5em}
	&  & CTD        & +0.1* & +10.3 & +13.3 & +15.1 & +10.1 & +11.5 & +14.9 & +17.0 & +18.4 & +16.3 & +18.1 \\
	\multicolumn{3}{c}{Overall} & +2.8 & +16.4 & +18.7 & +20.4 & +17.7 & +19.4 & +23.6 & +25.4 & +26.4 & +23.9 & +27.1 \\
	\hline
\end{tabular}
{* indicates the improvements that are NOT statistically significant  (i.e., $p>0.05$ in paired two-sample \textit{t}-tests).}
\end{threeparttable}
\end{adjustbox}
    \end{table}
    
\begin{table}[!htb]
\renewcommand{\baselinestretch}{1.2}
    \caption{Comparison of Online Scheduling Policies ($\Lambda\in\{0.2,0.4\}$): Maximum Decision Time (seconds)\label{TabOnlineDTSmallDetail}}
    \footnotesize
\addtolength\tabcolsep{-0.2em}
\begin{adjustbox}{center}
\begin{tabular}{cccrrrrrrrrrrrr}
	\hline
	\multicolumn{3}{c}{Instance} & \multicolumn{12}{c}{Online scheduling policy} \\\cmidrule(lr){1-3}\cmidrule(lr){4-15}
	\multirow{2}{*}{$\Lambda$} & \multirow{2}{*}{$K$} & \multirow{2}{*}{Dist.} & \multicolumn{1}{c}{\multirow{2}{*}{GP$_{_\textsf{CI}}$}} & \multicolumn{1}{c}{\multirow{2}{*}{GP$_{_\textsf{R}}$}} & \multicolumn{3}{c}{R$_{_\textsf{R}}$-GP$_{_\textsf{CI}}(H)$} & \multicolumn{1}{c}{\multirow{2}{*}{PFA$_{_\textsf{CI}}$}} & \multicolumn{1}{c}{\multirow{2}{*}{PFA$_{_\textsf{R}}$}} & \multicolumn{3}{c}{R$_{_\textsf{R}}$-PFA$_{_\textsf{CI}}(H)$} & \multicolumn{1}{c}{\multirow{2}{*}{S-PbP}} & \multicolumn{1}{c}{\multirow{2}{*}{PbP}} \\\cmidrule(lr){6-8}\cmidrule(lr){11-13}
	& & & & & \multicolumn{1}{c}{$H$=25} & \multicolumn{1}{c}{$H$=50} & \multicolumn{1}{c}{$H$=100} & & & \multicolumn{1}{c}{$H$=25} & \multicolumn{1}{c}{$H$=50} & \multicolumn{1}{c}{$H$=100} & & \\\cmidrule(lr){1-15}
	0.2 & 2 & UTI & $<$0.1 & 0.3 & 0.9 & 2.0 & 3.1 & $<$0.1 & 0.2 & 1.0 & 2.1 & 3.8 & 0.3 & 0.4 \\
	&  & CTI & $<$0.1 & 0.3 & 1.0 & 1.7 & 3.1 & $<$0.1 & 0.3 & 0.9 & 2.0 & 3.9 & 0.4 & 0.4 \\\vspace{0.5em}
	&  & CTD & $<$0.1 & 0.2 & 0.9 & 1.6 & 3.1 & $<$0.1 & 0.2 & 0.9 & 1.7 & 3.6 & 0.4 & 0.6 \\
	& 3 & UTI & $<$0.1 & 0.1 & 1.9 & 3.1 & 6.3 & $<$0.1 & 0.1 & 1.8 & 3.5 & 6.3 & 0.4 & 0.6 \\
	&  & CTI & $<$0.1 & 0.2 & 1.6 & 3.0 & 6.0 & $<$0.1 & 0.2 & 1.6 & 3.1 & 7.0 & 0.3 & 0.5 \\\vspace{0.5em}
	&  & CTD & $<$0.1 & 0.1 & 1.6 & 3.7 & 6.1 & $<$0.1 & 0.2 & 1.9 & 3.2 & 7.2 & 0.3 & 0.5 \\
	0.4 & 3 & UTI & $<$0.1 & 0.4 & 3.1 & 6.1 & 14.1 & $<$0.1 & 0.2 & 3.2 & 7.3 & 14.5 & 0.6 & 1.0 \\
	&  & CTI & $<$0.1 & 0.4 & 3.1 & 6.1 & 14.1 & $<$0.1 & 0.3 & 3.1 & 7.2 & 14.3 & 0.6 & 1.0 \\\vspace{0.5em}
	&  & CTD & $<$0.1 & 0.3 & 3.1 & 6.1 & 12.0 & $<$0.1 & 0.2 & 3.1 & 6.3 & 14.2 & 0.6 & 1.1 \\
	& 5 & UTI & $<$0.1 & 0.1 & 7.4 & 14.6 & 29.2 & $<$0.1 & 0.1 & 7.5 & 14.9 & 34.3 & 0.7 & 2.0 \\
	&  & CTI & $<$0.1 & 0.2 & 7.3 & 14.6 & 29.2 & $<$0.1 & 0.2 & 8.7 & 17.1 & 34.5 & 0.8 & 2.2 \\\vspace{0.5em}
	&  & CTD & $<$0.1 & 0.3 & 7.2 & 14.4 & 33.2 & $<$0.1 & 0.3 & 8.6 & 16.8 & 33.7 & 0.8 & 2.4 \\
	\multicolumn{3}{c}{Overall} & $<$0.1 & 0.4 & 7.4 & 14.6 & 33.2 & $<$0.1 & 0.3 & 8.7 & 17.1 & 34.5 & 0.8 & 2.4 \\
	\hline
\end{tabular}
\end{adjustbox}
    \end{table}
    
\begin{figure}[!htb]
\centering
	\caption{Comparison of Online Scheduling Policies R$_{_\textsf{CI}}$-GP$_{_\textsf{CI}}(H)$ and R$_{_\textsf{R}}$-GP$_{_\textsf{CI}}(H)$ ($\Lambda\in\{0.2,0.4\}$).\label{FigEvaluationRCI}}
	{\includegraphics[width=0.6\textwidth]{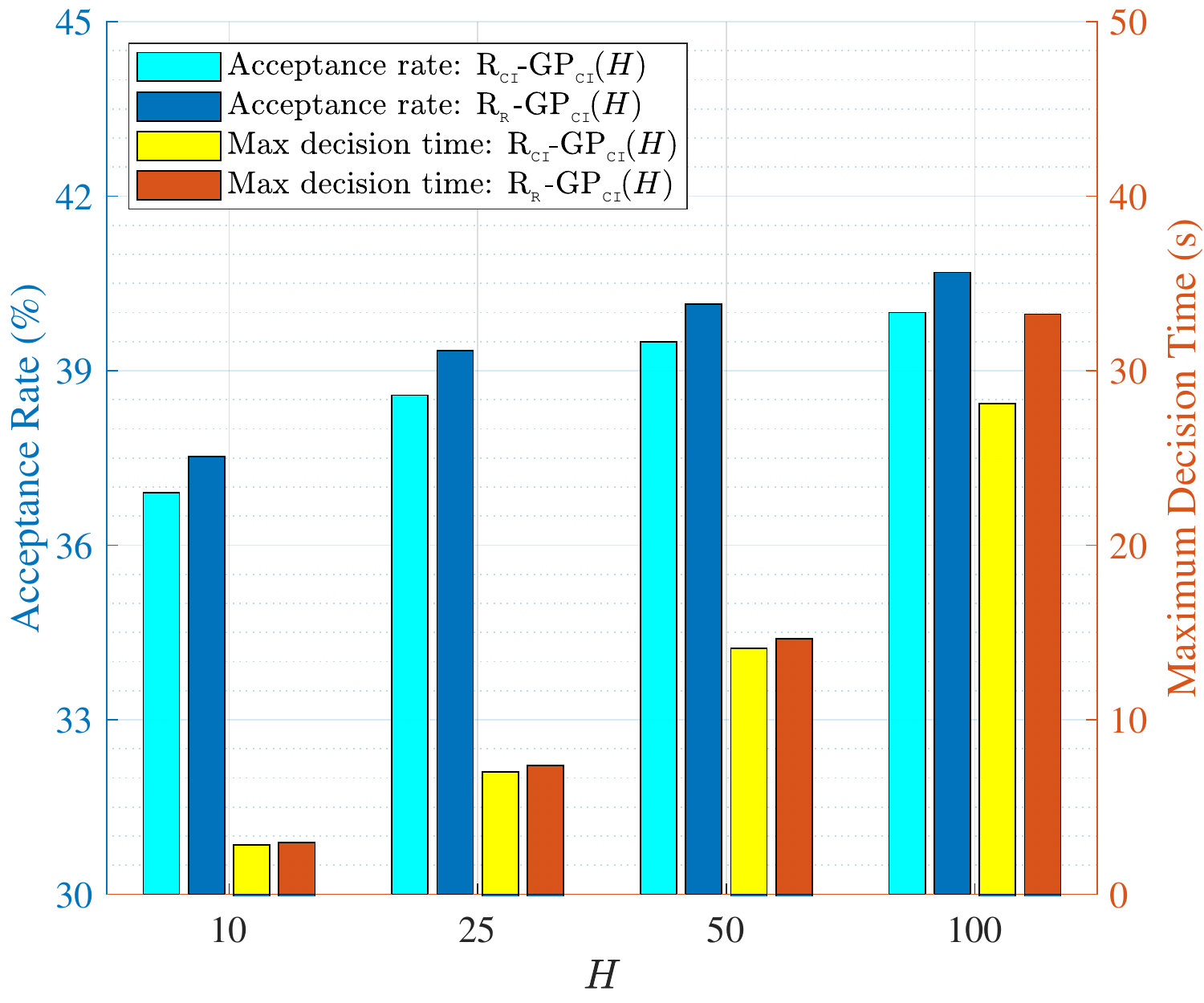}}
	{}
\end{figure}

\begin{figure}[!htb]
\centering
\caption{Comparison of Online Scheduling Policies ($\Lambda\in\{0.8,1.5\}$).\label{FigEvaluationRHBig}}
{\includegraphics[width=0.6\textwidth]{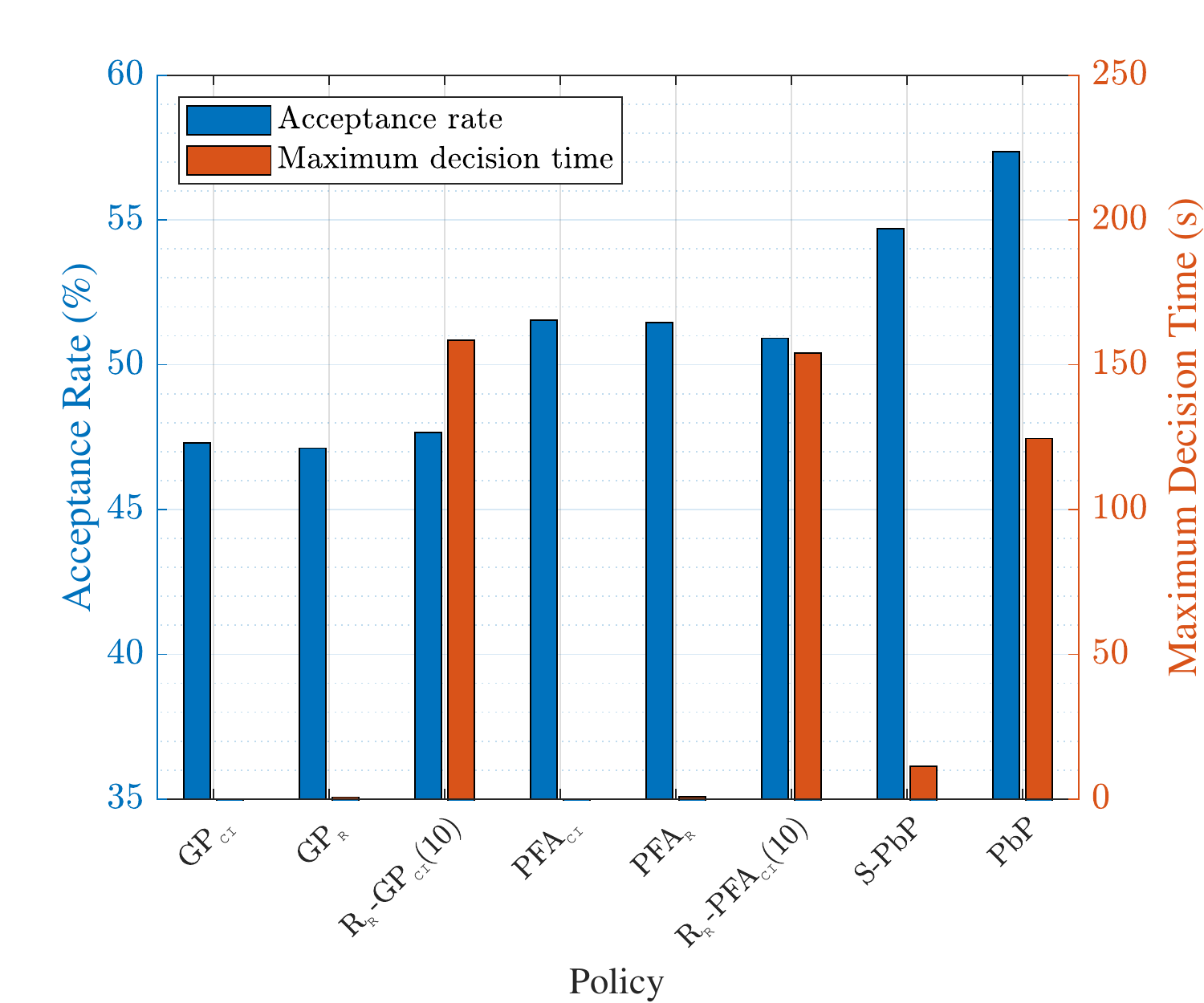}}
{}
\end{figure}

\begin{table}[!htb]
\renewcommand{\baselinestretch}{1.2}
    \caption{Comparison of Online Scheduling Policies ($\Lambda\in\{0.8,1.5\}$): Acceptance Rate (\%) and Percentage Improvement Compared to GP$_{_\textsf{CI}}$\label{TabOnlineARLargeDetail}}
    \footnotesize
    \addtolength\tabcolsep{-0.2em}
\begin{adjustbox}{center}
\begin{threeparttable}
\begin{tabular}{ccccccccccc}
	\hline
	\multicolumn{3}{c}{Instance} & \multicolumn{8}{c}{Online scheduling policy} \\\cmidrule(lr){1-3}\cmidrule(lr){4-11}
	$\Lambda$ & $K$ & Dist. & \multicolumn{1}{c}{GP$_{_\textsf{CI}}$} & \multicolumn{1}{c}{GP$_{_\textsf{R}}$} & \multicolumn{1}{c}{R$_{_\textsf{R}}$-GP$_{_\textsf{CI}}(10)$} & \multicolumn{1}{c}{PFA$_{_\textsf{CI}}$} & \multicolumn{1}{c}{PFA$_{_\textsf{R}}$} & \multicolumn{1}{c}{R$_{_\textsf{R}}$-PFA$_{_\textsf{CI}}(10)$} & \multicolumn{1}{c}{S-PbP} & \multicolumn{1}{c}{PbP} \\\cmidrule(lr){1-11}
	0.8 & 6 & UTI 	            & 30.9 & 31.8(+2.9*) 	& 32.3(+4.7) 	& 36.8(+19.4) 				& 37.2(+20.5) 				& 35.9(+16.5) 				& 38.9(+26.2) 				& 39.7(+28.6) \\
				&  & CTI 		& 30.6 & 30.0($-$2.0*) 	& 33.3(+8.8) 	& 37.3(+21.6) 				& 36.4(+18.8) 				& 37.8(+23.4) 				& 40.2(+31.3) 				& 41.3(+34.9) \\\vspace{0.5em}
				&  & CTD 		& 30.5 & 30.7(+0.7*) 	& 33.0(+8.0) 	& 36.7(+20.4) 				& 36.0(+18.1) 				& 37.0(+21.3) 				& 39.8(+30.5) 				& 40.3(+32.2) \\
				& 12 & UTI 		& 62.9 & 61.8($-$1.8*) 	& 60.9($-$3.3) & 64.6\hspace{0.5em}(+2.7) 	& 64.6\hspace{0.5em}(+2.7) 	& 63.5\hspace{0.5em}(+0.9) 	& 68.3\hspace{0.5em}(+8.5) 	& 72.8(+15.7) \\
				&  & CTI 		& 64.8 & 63.5($-$2.0*) 	& 63.8($-$1.5) & 66.7\hspace{0.5em}(+2.9) 	& 66.8\hspace{0.5em}(+3.1) 	& 65.7\hspace{0.5em}(+1.4) 	& 69.0\hspace{0.5em}(+6.6) 	& 74.6(+15.2) \\\vspace{0.5em}
				&  & CTD 		& 63.9 & 63.8($-$0.2*) 	& 62.6($-$2.1) & 66.3\hspace{0.5em}(+3.8) 	& 67.0\hspace{0.5em}(+4.8) 	& 65.0\hspace{0.5em}(+1.7) 	& 71.0(+11.1) 				& 74.1(+16.0) \\\vspace{0.5em}
				1.5 & 10 & UTI 	& 33.4 & 35.7(+7.0) 	& 34.6(+3.8) & 39.7(+19.1) 				& 39.9(+19.7) 				& 39.9(+19.5) 				& 43.1(+29.0) 				& 44.2(+32.5) \\\vspace{0.5em}
				& 20 & UTI 		& 62.4 & 64.5(+3.4) 	& 61.6($-$1.4) & 67.1\hspace{0.5em}(+7.4) 	& 66.7\hspace{0.5em}(+6.8) 	& 65.0\hspace{0.5em}(+4.1) 	& 70.8(+13.4) 				& 76.9(+23.2) \\
				\multicolumn{3}{c}{Overall} & 47.3 & 47.1($-$0.4*) & 47.7(+0.8) & 51.5\hspace{0.5em}(+8.9) & 51.5\hspace{0.5em}(+8.8) & 50.9\hspace{0.5em}(+7.6) & 54.7(+15.6) 			& 57.4(+21.2) \\
	\hline
\end{tabular}
{* indicates the improvements that are NOT statistically significant  (i.e., $p>0.05$ in paired two-sample \textit{t}-tests).}
\end{threeparttable}
\end{adjustbox}
    \end{table}

\begin{table}[!htb]
\renewcommand{\baselinestretch}{1.2}
    \caption{Comparison of Online Scheduling Policies ($\Lambda\in\{0.8,1.5\}$): Decision Times (seconds)\label{TabOnlineTLargeDetail}}
    \footnotesize
    \addtolength\tabcolsep{-0.3em}
\begin{adjustbox}{center}
\begin{tabular}{cccrrrrrrrrrrrrrrrr}
	\hline
	\multicolumn{3}{c}{Instance} & \multicolumn{16}{c}{Online scheduling policy} \\\cmidrule(lr){1-3}\cmidrule(lr){4-19}
	\multirow{2}{*}{$\Lambda$} & \multirow{2}{*}{$K$} & \multirow{2}{*}{Dist.} & \multicolumn{2}{c}{GP$_{_\textsf{CI}}$} & \multicolumn{2}{c}{GP$_{_\textsf{R}}$} & \multicolumn{2}{c}{R$_{_\textsf{R}}$-GP$_{_\textsf{CI}}(10)$} & \multicolumn{2}{c}{PFA$_{_\textsf{CI}}$} & \multicolumn{2}{c}{PFA$_{_\textsf{R}}$} & \multicolumn{2}{c}{R$_{_\textsf{R}}$-PFA$_{_\textsf{CI}}(10)$} & \multicolumn{2}{c}{S-PbP} & \multicolumn{2}{c}{PbP} \\\cmidrule(lr){4-5}\cmidrule(lr){6-7}\cmidrule(lr){8-9}\cmidrule(lr){10-11}\cmidrule(lr){12-13}\cmidrule(lr){14-15}\cmidrule(lr){16-17}\cmidrule(lr){18-19}
	&  &  & \multicolumn{1}{c}{Avg.} & \multicolumn{1}{c}{Max} & \multicolumn{1}{c}{Avg.} & \multicolumn{1}{c}{Max} & \multicolumn{1}{c}{Avg.} & \multicolumn{1}{c}{Max} & \multicolumn{1}{c}{Avg.} & \multicolumn{1}{c}{Max} & \multicolumn{1}{c}{Avg.} & \multicolumn{1}{c}{Max} & \multicolumn{1}{c}{Avg.} & \multicolumn{1}{c}{Max} & \multicolumn{1}{c}{Avg.} & \multicolumn{1}{c}{Max} & \multicolumn{1}{c}{Avg.} & \multicolumn{1}{c}{Max} \\\cmidrule(lr){1-19}
	0.8 & 6 & UTI & $<$0.1 & $<$0.1 & 0.1 & 0.4 & 2.7 & 8.0 & $<$0.1 & $<$0.1 & 0.1 & 0.4 & 3.9 & 9.4 & 0.8 & 1.8 & 2.3 & 6.1 \\
	&  & CTI & $<$0.1 & $<$0.1 & 0.1 & 0.6 & 2.7 & 8.0 & $<$0.1 & $<$0.1 & 0.1 & 0.7 & 4.1 & 9.4 & 0.8 & 1.8 & 2.2 & 5.5 \\\vspace{0.5em}
	&  & CTD & $<$0.1 & $<$0.1 & 0.1 & 0.6 & 2.7 & 7.8 & $<$0.1 & $<$0.1 & 0.1 & 0.7 & 4.1 & 9.2 & 0.9 & 1.9 & 2.4 & 6.0 \\
	& 12 & UTI & $<$0.1 & $<$0.1 & 0.1 & 0.2 & 13.9 & 29.3 & $<$0.1 & $<$0.1 & 0.1 & 0.2 & 14.6 & 29.1 & 1.5 & 3.4 & 9.2 & 23.8 \\
	&  & CTI & $<$0.1 & $<$0.1 & 0.1 & 0.3 & 13.9 & 28.1 & $<$0.1 & $<$0.1 & 0.1 & 0.3 & 14.8 & 30.4 & 1.6 & 3.7 & 9.1 & 23.3 \\\vspace{0.5em}
	&  & CTD & $<$0.1 & $<$0.1 & 0.1 & 0.2 & 15.5 & 32.5 & $<$0.1 & $<$0.1 & 0.1 & 0.2 & 14.4 & 28.0 & 1.6 & 3.6 & 8.5 & 20.7 \\\vspace{0.5em}
	1.5 & 10 & UTI & $<$0.1 & $<$0.1 & 0.1 & 0.5 & 13.5 & 40.0 & $<$0.1 & $<$0.1 & 0.1 & 0.4 & 19.3 & 39.7 & 2.6 & 5.7 & 11.2 & 27.0 \\\vspace{0.5em}
	& 20 & UTI & $<$0.1 & $<$0.1 & 0.1 & 0.2 & 80.6 & 158.4 & $<$0.1 & $<$0.1 & 0.1 & 0.3 & 81.8 & 154.1 & 5.1 & 11.3 & 51.8 & 124.5 \\
	\multicolumn{3}{c}{Overall} & $<$0.1 & $<$0.1 & 0.1 & 0.6 & 11.0 & 158.4 & $<$0.1 & $<$0.1 & 0.1 & 0.7 & 11.9 & 154.1 & 1.3 & 11.3 & 7.2 & 124.5 \\
	\hline
\end{tabular}
\end{adjustbox}
    \end{table}

\end{document}